\documentclass{SCIYA2023enOL}
\online

\newtheorem{assumptions}[theorem]{Assumptions}
\usepackage{caption}
\usepackage{subfig}
\usepackage{graphicx}
\usepackage{makecell}
\begin{document}

\ensubject{fdsfd}

\ArticleType{}
\Year{}
\Month{}%
\Vol{}
\No{}
\BeginPage{1} %
\DOI{}

\title[SMC with Gaussian Mixture Approximation]{Sequential Monte Carlo with Gaussian Mixture Approximation for Infinite-Dimensional Statistical Inverse Problems}
{Sequential Monte Carlo with Gaussian Mixture Approximation for Infinite-Dimensional Statistical Inverse Problems}

\author[1]{Haoyu Lu}{{luhaoyu@stu.xjtu.edu.cn}}
\author[1,$\ast$]{Junxiong Jia}{{jjx323@xjtu.edu.cn}}
\author[1,2]{Deyu Meng}{dymeng@mail.xjtu.edu.cn}

\AuthorMark{H. Lu, J. Jia, and D. Meng}

\AuthorCitation{H. Lu, J. Jia, and D. Meng}

\address[1]{School of Mathematics and Statistics,
Xi'an Jiaotong University,
Xi'an {\rm 710049},
China}
\address[2]{Pazhou Laboratory (Huangpu)}

\abstract{
By formulating the inverse problem of partial differential equations (PDEs) as a statistical inference problem, the Bayesian approach provides a general framework for quantifying uncertainties. In the inverse problem of PDEs, parameters are defined on an infinite-dimensional function space, and the PDEs induce a computationally intensive likelihood function. Additionally, sparse data tends to lead to a multi-modal posterior. These features make it difficult to apply existing sequential Monte Carlo (SMC) algorithms. To overcome these difficulties, we propose new conditions for the likelihood functions, construct a Gaussian mixture based preconditioned Crank-Nicolson transition kernel, and demonstrate the universal approximation property of the infinite-dimensional Gaussian mixture probability measure. By combining these three novel tools, we propose a new SMC algorithm with Gaussian mixture approximation, together with an easy-to-use reduced version. For this new algorithm, we obtain a convergence theorem that allows Gaussian priors, illustrating that the sequential particle filter actually reproduces the true posterior distribution. Furthermore, the proposed new algorithm is rigorously defined on the infinite-dimensional function space, naturally exhibiting the discretization-invariant property. Numerical experiments demonstrate that the reduced version has a strong ability to probe the multi-modality of the posterior, significantly reduces the computational burden, and numerically exhibits the discretization-invariant property (important for large-scale problems).
}



\maketitle

\section{Introduction}\label{Sec1}
In fields such as seismic exploration and medical imaging, noisy measurements of solutions to partial differential equations (PDEs) are typically available, and there is often interest in the parameters within these PDEs. This is known as an inverse problem, which has experienced tremendous growth over the past few decades \cite{Benning2018}. Furthermore, we aim not only to estimate the unknown parameters but also to quantify the uncertainty of these estimates, such as variance and confidence intervals, necessitating Bayesian modeling for the inverse problem. The Bayesian framework's viability for inverse problems was extensively reviewed in \cite{Kaipio2005}, which demonstrated the potential for Markov chain Monte Carlo methods (MCMC, \cite{Beskos2015SC}).

Inverse problems of PDEs have generally been defined on infinite-dimensional spaces, which are not compatible with the well-studied finite-dimensional Bayesian inference approach. To overcome this obstacle, the “Bayesianize-then-discretize” strategy is employed. Bayes' formula and algorithms are initially constructed in infinite-dimensional space, and once the infinite-dimensional algorithm is established, a finite-dimensional approximation is carried out \cite{Dashti2017}. Despite the advantages of the Bayesian framework, the non-linear forward problems often lead to complex posterior distributions, making it infeasible to analytically derive the posterior distributions.

The MCMC method, serving as a vital component within the Bayesian framework, provides the means to navigate and sample from these intricate posterior distributions. In the infinite-dimensional setting, the preconditioned Crank-Nicolson (pCN) and pCN Langevin transition kernel can be derived from the discretization of the Langevin equation \cite{Andrew_review,Beskos2008,Hairer2007}. However, the Bayesian inference for computationally intensive models can be limited by the expense of MCMC sampling. In theory, MCMC methods can sample from the posterior distribution, but in practice, the computational cost of MCMC in PDE inverse problems can be very high. Additionally, the burn-in process and the difficulty of parallelizing traditional MCMC \cite{Ye2024} further contribute to the computational expense.

A parallelizable alternative to MCMC without the need for burn-in is sequential Monte Carlo (SMC, \cite{Dai2022}), a class of algorithms used to approximate distributions of interest. SMC methods combine particle filtering, importance sampling, tempering ideas \cite{Latz2021}, and MCMC sampling. Early SMC methods were developed in a finite-dimensional setting \cite{Moral2006}. Later, in \cite{Kantas2014}, they were applied to high-dimensional PDE inverse problems, and in \cite{Beskos2015SC}, the convergence of SMC in the infinite-dimensional setting was proven, achieving a convergence rate that is independent of the discrete dimension.

Besides, MCMC can also be used as a transition kernel in SMC. Therefore, SMC and MCMC share many common features, such as the need to develop the theory in the infinite-dimensional function space to establish a SMC sampling method for unknown function parameters. Additionally, the computational cost of SMC is high due to the use of MCMC within SMC. To save time costs, the following strategies can be employed:
\begin{itemize}
	\item Add likelihood information into the MCMC transition kernel to accelerate MCMC mixing \cite{Andrew_review}, which will be discussed in detail below.
	\item Use a surrogate or approximate likelihood to avoid the expensive computation cost if possible \cite{Bon2021}.
	\item Combine other useful tools with MCMC, such as variational inference \cite{Naesseth2017}.
	\item Adjust the structure of SMC. For example, the multi-level-SMC method was proposed in \cite{Beskos2017} and further discussed in \cite{Beskos2018, Latz2018}.
\end{itemize}

The first strategy has an extremely wide range of applications. The choices of transition kernels in SMC methods are highly flexible. For example, the basic random walk and its infinite-dimensional improvement, pCN, are used in SMC as described in \cite{Kantas2014} and \cite{Beskos2015SC}. The dimension-independent likelihood-informed transition kernel is utilized in \cite{Cui2024}. Another SMC method that leverages the gradient information of the likelihood function is presented in \cite{Rosato2024}. Finally, a series of studies \cite{Gunawan2020, Buchholz2021, Burda2023} employs Hamiltonian Monte Carlo as the transition kernel to enhance sampling efficiency.

In this paper, we focus on improving SMC in infinite-dimensional spaces. To further reduce computational costs, we consider improvements for the MCMC transition steps. SMC methods are particularly attractive for a diverse range of inference challenges due to their ability to effectively manage multi-modal distributions. We aim to preserve these desirable approximation properties while performing the approximation. For handling multi-modal distributions, the Gaussian mixture is an appropriate choice. In fact, it has been proven that a Gaussian mixture can approximate any distribution in the finite-dimensional case \cite{Nguyen2020}. The Gaussian mixture approximation, employed in particle filters \cite{Raihan2018}, is a commonly used method for approximation. Gaussian mixture approximations have been applied in SMC filtering. In \cite{Sehyun2019}, the authors constructed the prior samples into a Dirac measure, regarded it as a Gaussian mixture distribution with zero variance, and then linearized the model to update the parameters of the Gaussian mixture, ultimately obtaining an approximation of the posterior. 

Unlike these works, we will combine the first two strategies by incorporating a likelihood-informed Gaussian mixture into the pCN proposal, resulting in the pCN-GM method. Employing pCN-GM as the mutation step in SMC leads to the SMC-pCN-GM algorithm. 
Furthermore, by replacing the mutation step with a Gaussian mixture sampler, we obtain a reduced version, termed SMC-GM.

The main difficulties are as follows: First, existing SMC convergence theorems require the potential function to be bounded \cite{Beskos2015SC}, yet commonly used Gaussian priors and many forward problems, such as Darcy flow, cannot apply this theorem. Second, traditional MCMC simulators using Metropolis-Hastings do not necessarily behave well when dealing with sampling from multi-modal distributions \cite{Daviet2018}. Additionally, in the infinite-dimensional setting, there are singularity issues between measures, such as changes of the mean of a Gaussian measure, which almost surely lead to singularities \cite{Bogachev1998}. Third, while mixture Gaussian densities can approximate density functions well in finite-dimensional settings, can infinite-dimensional Gaussian mixture measures also approximate the posterior measure, and if so, in what sense of distance? 

In response to these issues, overall, we have made the following contributions:
\begin{enumerate}
	\item We propose a weaker condition and prove a convergence theorem of SMC, thereby providing a theoretical foundation for sampling from a posterior measure derived from a Gaussian prior using SMC.
	\item We derive a generalized pCN algorithm based on Gaussian mixtures (pCN-GM) from the Crank-Nicolson discretization of the Langevin system that describes multiple particles with various perturbations. This prevents the occurrence of singularities, i.e., the pCN-GM method is well-defined under the infinite-dimensional setting. This algorithm can be viewed as an extension of both pCN and independent sampling with Gaussian mixtures \cite{Feng2018}. 
	\item Using pCN-GM in the mutation step of SMC yields the SMC-pCN-GM method. Based on the fact that finite-dimensional Gaussian mixture density functions can approximate any density function, along with the continuity of the likelihood function and the properties of infinite-dimensional Gaussian measures, we prove the denseness of Gaussian mixture measures in infinite-dimensional spaces. This motivates a reduced version of the SMC-pCN-GM method.
\end{enumerate}

This paper is organized as follows. In Section \ref{Sec.SMC-GM} 
we propose the SMC-pCN-GM method, together with its reduced version, SMC-GM. Both methods are SMC algorithms with Gaussian mixture approximation. Subsection \ref{Subsec.SMC} provides a brief introduction to the SMC method. In Subsection  \ref{Subsec.ApproximationTheorySMC}, we prove the convergence theorem for SMC under weak conditions. In Subsection \ref{Subsec.pCN-GMM}  we propose the pCN-GM kernel, a new Markov transition kernel,  and prove that it is well-defined in the infinite-dimensional space. In Subsection \ref{Subsec.GM}  we present the universal approximation property of Gaussian mixture measures in the infinite-dimensional setting, and give the approximated SMC-GM method. Finally, Section \ref{Sec.Numerics} applies the developed SMC-GM method to three typical inverse problems.
All proofs are provided in \ref{SM-AllProofs}.  All the programming codes are accessible at \url{https://github.com/jjx323/IPBayesML-SMC-GM}.

\section{SMC with Gaussian mixture approximation}\label{Sec.SMC-GM}
In this section, we adopt the SMC methods defined in infinite-dimensional function spaces \cite{Beskos2015SC} and implement three critical improvements. First, we provide a convergence theorem for SMC under weaker conditions. Second, we combine the Gaussian mixture with pCN to introduce a new transition kernel, known as the pCN-GM algorithm, which, when utilized as the transition kernel in SMC, forms the SMC-pCN-GM method. Lastly, we present an approximation to SMC-pCN-GM as a reduced version, the SMC-GM algorithm.

\subsection{Algorithm overview}\label{Subsec.SMC}
Let \( \mathcal{H} \) be a separable Hilbert space, \( \mathcal{B}(\mathcal{H}) \) be the Borel \( \sigma \)-algebra, and \( N_d \) be a positive integer. Denote the Gaussian measure on the measurable space \((\mathcal{H},\mathcal{B}(\mathcal{H}))\) with mean \(m\) and covariance \(\mathcal{C}\) as \(\mathcal{N}(m,\mathcal{C})\), and Gaussian measure on \(\mathbb{R}^{N_d}\)  with zero mean and covariance \(C\)  as \(N(0,C)\). The forward problem is defined by
\begin{align}
\boldsymbol{d} = \mathcal{F}(u) + \epsilon,\quad \epsilon \sim N(0, \sigma^2I),
\end{align}
where \( u \in \mathcal{H} \) represents the parameter of interest, the mapping \( \mathcal{F} : \mathcal{H} \rightarrow \mathbb{R}^{N_d} \) is the forward operator determined by a PDE, \(\sigma^2\in \mathbb{R}^+\) is the noise variance, \(I\in \mathbb{R}^{N_d\times N_d}\) is the identity matrix, \( \epsilon \) represents the observation error, and \( \boldsymbol{d} \in \mathbb{R}^{N_d} \) is the measurement data. The negative log-likelihood function, or potential function, is given by
\begin{equation}
\Phi(u)=\frac{1}{2\sigma^2}\|\mathcal{F}(u)-\boldsymbol{d}\|^2,
\end{equation}
where the norm is defined as \(\|a\|:=\sqrt{a^Ta}\). 

Within the Bayesian framework, we must select a prior measure and sample from the posterior distribution to estimate the statistics of \(u\). Let \(\mathcal{N} (0,\mathcal{C})\) be the prior measure, denoted as \(\mu_0\). Assume that the covariance operator \( \mathcal{C}\) is positive, self-adjoint, and of trace class, that is, it satisfies the condition
\begin{align}
\operatorname{tr}(\mathcal{C}) := \sum_{i=1}^\infty \lambda_i < \infty,
\end{align}
where \( \lambda_i \) denotes the \( i \)-th eigenvalue of \( \mathcal{C} \). With these assumptions, \( \mathcal{N}(0, \mathcal{C}) \) is well-defined on \( (\mathcal{H}, \mathcal{B}(\mathcal{H})) \). The Bayes' theorem \cite{Dashti2017} can be applied to determine the posterior
\begin{align}
\frac{d\mu^\mathrm{d}}{d\mu_0}(u) = \frac{1}{Z} \exp\left\{-\Phi(u)\right\},\label{Bayes}
\end{align}
where \( Z \) is the normalization constant, ensuring that the posterior distribution is properly normalized. The inverse problem is to estimate statistics of \(u\) using the data \(\boldsymbol{d}\). In the Bayesian framework, the estimation is obtained by sampling from the posterior \(\mu^\mathrm{d}\).

Now we provide a brief introduction to the SMC method, for more details, see \cite{Moral2006,Beskos2015SC}. 
Recall that \( \mu_0 \) is the prior measure on a separable Hilbert space \( \mathcal{H} \), equipped with the Borel \( \sigma \)-algebra \( \mathcal{B}(\mathcal{H}) \). We aim to sample from the posterior measure \( \mu^\mathrm{d} \) defined in (\ref{Bayes}). The SMC methods share the general structure of particle filters, which uses a \( N \)-particle Dirac measure
\begin{align}
\mu^N = \sum_{n=1}^N w^{(n)} \delta_{v^{(n)}}
\label{Dirac approximation}
\end{align}
to approximate the posterior measure \( \mu^\mathrm{d} \), where \(\{w^{(n)}\}_{n=1}^N\) are \(N\) positive numbers satisfying \(\sum_{n=1}^Nw^{(n)}=1\), representing weights. And the objective of the SMC method is to find particles \( \{v^{(n)}\}_{n=1}^N \) such that \( \mu^N \) is close to \( \mu^\mathrm{d} \).

\begin{remark}
	Note that finding the appropriate weights here can be reduced to finding the appropriate samples, as samples with large weights will be replicated multiple times during the SMC process \cite{Kantas2014}.
\end{remark}

The SMC methods use importance sampling to sample from the posterior. However, since it is only suitable for posteriors that are similar to the prior, it would be difficult to sample directly from \( \mu^\mathrm{d} \). Thus, the SMC methods employ the concept of intermediate measures for stratification in the field of inverse problems for PDEs \cite{Beskos2015SC}. Given an integer \( J \), then for \( 1 \leq j \leq J \), we can define a sequence of measures \( \mu_j \ll \mu_0 \) by the Radon-Nikodym derivative:
\begin{align*}
\Phi_i(u) = h_i \Phi(u), \quad \frac{d\mu_j}{d\mu_0}(u) \propto \exp\left(-\sum_{i=1}^j \Phi_i(u)\right),
\end{align*}
where \( \{h_j\}_{j=1}^J \) is a positive sequence satisfying \( \sum_{j=1}^J h_j = 1 \). The method for determining \(h_i\) is described in \ref{SM-NumericalDetails}. Now we have a measure sequences \(\{\mu_j\}_{j=0}^J\), and the importance sampling can be used to sample from \(\mu_{j+1}\) based on samples from \(\mu_j\).

The entire SMC algorithm consists of three iterative steps, the re-weighting step, the resampling step and the mutation step (See Algorithm \ref{alg SMC}). First, importance sampling is used to obtain samples from the current layer to the next, which is specifically reflected in the updating of weights. Second, to discard low-probability points and sample more high-probability points, resampling can be incorporated during the iteration, which may reduce sample diversity. Third, a \(\mu_j\)-invariant transition kernel is used to enhance sample diversity. Ultimately, the SMC methods achieve an approximation of the posterior distribution.

In the mutation step, we can use the MCMC method as a \(\mu_j\)-invariant transition kernel. The simplest MCMC transition kernel are random walk and the pCN method, which are used in SMC in \cite{Beskos2015SC} and \cite{Kantas2014}, respectively. A general SMC frame is outlined in Algorithm \ref{alg SMC}.
\begin{algorithm}
	\caption{Sequential Monte Carlo}
	\label{alg SMC}
	\begin{algorithmic}[1]
		\STATE{\textbf{Input}: Specify \(\mu_0\). Draw \(\{v_0^{(n)}\}_{n=1}^N\sim \mu_0.\)}
		\STATE{
			
			\textbf{For} \(j=0,1,\cdots ,J-1\), \textbf{do}:
			
			}
		\STATE{\qquad Set \(w_j^{(n)}=1/N,n=1,2,\cdots, N\) and define \(\mu_j^N=\sum_{n=1}^Nw^{(n)}_{j}\delta_{v^{(n)}_{j}}\) ;}
		
		\STATE{\qquad (Re-weighting) Update \(w_{j+1}^{(n)}\) by:}
		\begin{align*}
		\hat{w}_{j+1}^{(n)}=e^{-h_{j+1}\Phi(v_{j}^{(n)})}w_j^{(n)},\qquad w_{j+1}^{(n)}=\frac{\hat{w}_{j+1}^{(n)}}{\sum_{n=1}^N\hat{w}_{j+1}^{(n)}}
		\end{align*}
		\STATE{\qquad (Resampling) Update \(\hat{v}_{j+1}^{(n)}\) by resampling from  \(\sum_{n=1}^Nw^{(n)}_{j+1}
			\delta_{v^{(n)}_{j}}\).}
        \STATE{\qquad (Mutation) Draw \(\{v_{j+1}^{(n)}\}_{n=1}^N\) from a \(\mu_{j+1}\)-preserved MCMC transition kernel;}
		\STATE{\textbf{Output}: \(\mu_J^N\).}
	\end{algorithmic}
\end{algorithm}

To clarify our improvement process, we list the difficulties to be overcome as follows:
\begin{itemize}
	\item Difficulty 1: We will use a transition kernel based on an improved pCN, which, like pCN, applies a Gaussian prior. However, the Bayesian inverse problem constructed from this prior and the Darcy flow forward problem does not meet the conditions required by existing SMC convergence theorems \cite{Beskos2015SC}.
	\item Difficulty 2: In infinite-dimensional spaces, the proposal of the well-posed pCN algorithm involves mixing the current state with prior samples. Such a proposal is not ideal for complex multi-modal posteriors. In fact, traditional MCMC simulators using Metropolis-Hastings do not necessarily behave well when dealing with sampling from multi-modal distributions \cite{Daviet2018}.
	\item Difficulty 3: Finally, SMC can involve a large number of likelihood evaluations, which is expensive for models with computationally intensive likelihood functions \cite{Bon2021}.
\end{itemize}

In general, we aim to propose an efficient SMC for infinite-dimensional spaces capable of handling complex multi-modal posteriors. In the remainder of this section, we will address these challenges. In Section~\ref{Subsec.ApproximationTheorySMC}, we propose a weaker condition than Equation (19) in \cite{Beskos2015SC} and prove an analogous SMC convergence theorem under this condition. This theorem addresses Difficulty 1, providing a theoretical foundation for using transition kernels like pCN and our improved version in SMC. In Section~\ref{Subsec.pCN-GMM}, we introduce the Gaussian mixture-based pCN method, demonstrating its well-definedness. The incorporation of a Gaussian mixture addresses Difficulty 2. Finally, in Section~\ref{Subsec.GM}, we demonstrate the denseness of the mixed Gaussian measure in infinite-dimensional spaces and propose the SMC-GM algorithm, an approximation of SMC-pCN-GM. The high efficiency and multi-modal posterior sampling capability of SMC-GM will be demonstrated in Section \ref{Sec.Numerics}, addressing Difficulty 3.

To clarify the various SMC methods discussed in this paper, each corresponding to a different transition kernel, we summarize them as follows:
\begin{itemize}
	\item SMC-pCN: The SMC method utilizing the pCN method as the transition kernel.
	\item SMC-pCN-GM: The SMC method utilizing the pCN-GM method as the transition kernel. Here pCN-GM method is an improved version of pCN detailed in Section \ref{Subsec.pCN-GMM}.
	\item SMC-GM: The SMC method that replaces the transition kernel with a Gaussian mixture sampler, as introduced in Section \ref{Subsec.GM}.
	\item SMC-RW: The SMC method employing the random walk method as the transition kernel, used in Section \ref{example2} for numerical comparison.
\end{itemize}
For readers' convenience, we detail the aforementioned transition kernels (pCN, pCN-GM, GM, RW) in \ref{SM-4Kernels}. To adapt an SMC-XX method, simply replace the mutation step in Algorithm \ref{alg SMC} with the XX method. 
The algorithms proposed in this paper are SMC methods with Gaussian mixture approximations, including SMC-pCN-GM and its reduced version SMC-GM.

\subsection{Approximation theory of SMC}\label{Subsec.ApproximationTheorySMC}
In this section, we establish the convergence theory for the SMC algorithm with the following steps: First, we regard the entire algorithm as a mapping from the prior measure to the approximate posterior measure, and decompose this entire mapping into a series of mappings. Second, we estimate the error introduced by each mapping. Finally, we estimate the error of the entire mapping process, leading to the convergence theorem. In this subsection, we make the following assumptions:
\begin{assumptions}\label{A1}
	The potential functions $ \Phi_j $ satisfy the following conditions:
	
	\noindent (a) They are lower bounded, i.e., there exists \(\kappa_1>0\) such that for all \(j\in\{1,2,\cdots,\)\(J\}\), we have\begin{align*}\exp\left(-\Phi_j(u)\right)\leq \kappa_1^{-1}.\end{align*}
	(b) They are finite almost everywhere, i.e., for all \( j\in\{1,2,\cdots,J\}\) we have
 \begin{align*}\Phi_j(u)<\infty, \quad \mu_0-a.e.\end{align*}
\end{assumptions}
\begin{remark}
	It has been proved that the SMC method has an approximation theory in    \cite{Beskos2015SC}, where the potential function \(\Phi(u)\) is assumed to have an upper bound. However, for the Gaussian prior and the Darcy flow forward problem, this condition is not met. To address this limitation, we propose new assumptions (b)  to replace the upper boundedness, which can be verified under the conditions discussed above (the Darcy flow forward problem and the Gaussian prior). Assumptions \ref{A1}.(a) is easily satisfied because \(\Phi(u)=\frac{1}{2\sigma^2}\|\mathcal{F}u-\boldsymbol{d}\|^2\geq 0\). With Assumptions \ref{A1}.(b), we can prove that $ \{\mu_j\}_{j=1}^J  $ are mutually equivalent and, additionally,
	\begin{align}
	\frac{d\mu_{j}}{d\mu_{j-1}}(u)=\frac{1}{Z_{j}}\exp(-\Phi_j(u)),
	\end{align}
    {\color{black}where 
    \begin{align}\label{defZj}
    Z_j:=\int_\mathcal{H}\exp\left(-\Phi_j(u)\right)\mu_{j-1}(du).
    \end{align}
    In what follows, each $Z_j$ is defined analogously.} Hence, the distance between the measures \( \mu_j \) and \( \mu_{j-1} \) can be closer than that between \( \mu_0 \) and \( \mu_J \), which is advantageous for sampling in SMC.
\end{remark}

The SMC method consists of three steps: mutation, re-weighting, and resampling. First, We denote the transition kernel in the mutation step with \(P_j\), which preserves $ \mu_j $:
\begin{align}
(P_j\mu_j)(dv)=&\int_{w\in \mathcal{H}}P_j(w,dv){\color{black}\mu_j}(dw).\label{P_j}
\end{align}
 We then denote the re-weighting step as
\begin{align}
(L_j\mu)(dv)=&\frac{\exp(-\Phi_j(v))\mu(dv)}{\int_\mathcal{H}\exp(-\Phi_j(v))\mu(dv)}.\label{L_j}
\end{align}
Here we use \(L_j\) because this step corresponds the multiplication of the \(j\)-th potential function. With this notation we can express the true posterior distribution as
\begin{equation}
\mu_J=L_JL_{J-1}\cdots L_2L_1\mu_0.\label{truepost}
\end{equation}
The third operator is the sampling operator, because we use a weighted sum of Dirac measure to approximate the posterior:
\begin{equation}
(S^N \mu)(dv)=\frac{1}{N}\sum^N_{n=1}\delta_{v^{(n)}}(dv),\quad v^{(n)}\sim \mu,\ i.i.d ..\label{S^N}
\end{equation}
Using these notations, the measure of the next layer in SMC can be expressed in terms of the measure of the previous layer as follows:\color{black}
\begin{equation}
\mu_{j+1}^N=S^NP_{j+1}L_{j+1}\mu_j^N,\quad 0\leq j\leq J-1,\quad  \mu_0^N=S^N\mu_0.\label{decidemuJN}
\end{equation}\color{black}
To measure the error between \(\mu_J^N\) determined by (\ref{decidemuJN}) and the posterior in (\ref{truepost}), we need to consider the error introduced by these operators, respectively. To determine the distance between measures, we use the same metric in \cite{Dashti2017}:
\begin{equation}
d(\mu,\nu)=\sup_{|f|_\infty\leq 1}\sqrt{\mathbb{E}|\mu(f)-\nu(f)|^2},\label{random TV distance}
\end{equation}
where $ \mu(f):=\int_\mathcal{H}f(v)\mu(dv) $  represents the average value of the function \(f\) with respect to the measure \(\mu\), and \(\mathbb{E}\) is the expectation with respect to the samples within \(\mu\) and \(\nu\).  
\begin{remark}
	Note that this distance is designed for the measures associated with random samples, and if the measures is deterministic, the distance (\ref{random TV distance}) degenerates into
	\begin{equation}
	\sup_{|f|_\infty\leq 1}|\mu(f)-\nu(f)|=d_{\text{\tiny TV}}(\mu,\nu).\label{distance}
	\end{equation}
	which is the total variation distance. For more details, see \ref{SM-TV}.
\end{remark}

Let \(\mathcal{P}(\mathcal{H})\) denote the set of all probability measures on \(\mathcal{H}\).
 We can estimate the error induced by the above three operators \(S^N,P_j\), and \(L_j\) defined in \eqref{S^N},\eqref{P_j}, and \eqref{L_j}:
\begin{lemma}\label{Lemma1-3}
	We have the following results:
	
\noindent	(a) The sampling operator $ S^N $ satisfies\begin{align}
	\sup_{\mu\in \mathcal{P}(\mathcal{H})}d(S^N\mu,\mu)\leq \frac{1}{\sqrt{N}}.
	\end{align}
	(b) For any \(\nu,\nu'\in \mathcal{P}(\mathcal{H})\), the Markov kernel $ P_j $ satisfies\begin{align}
	d(P_j\nu,P_j\nu')\leq d(\nu,\nu').
	\end{align}
	(c) Under Assumptions \ref{A1},  for any \(\mu\in \mathcal{P}(\mathcal{H})\), we have
	\begin{align}
	d(L_j\mu_{j-1},L_j\mu)\leq 2\kappa_1^{-1}Z_{j}^{-1}d(\mu_{j-1},\mu).
	\end{align}
\end{lemma}

\begin{remark}
Proof of Lemma \ref{Lemma1-3}.(a),(b), see \cite{Dashti2017}, and proof of Lemma \ref{Lemma1-3}.(c), see \ref{SM-ProofLemmaC}. 
\end{remark}

Finally we give the approximation theory of SMC:
\begin{theorem}
If Assumptions \ref{A1} holds, then we have \label{SMC approximation theorem}
\begin{align*}
d(\mu_J^N, \mu_J) \le \frac{1}{\sqrt{N}} \left(
1 + \sum_{j=1}^{J} \frac{2^j}{\kappa_1^j} \prod_{k=1}^{j}\frac{1}{Z_{J+1-k}}
\right).
\end{align*} 
\end{theorem}
Theorem \ref{SMC approximation theorem} ensures the rationality of SMC for posterior sampling derived from a Gaussian prior. The proof is provided in \ref{SM-ConvergenceSMC}.  
In the next subsection, we will refine the pCN algorithm, which is also designed for posteriors derived from a Gaussian prior. Without Theorem \ref{SMC approximation theorem}, using pCN as a transition kernel would lack theoretical justification.

\subsection{pCN algorithm based on Gaussian Mixtures}\label{Subsec.pCN-GMM}
In this subsection, we construct a \(\mu_j\)-preserved transition kernel associated with Gaussian mixture measures, that is, the pCN algorithm based on Gaussian Mixtures (pCN-GM), which can be considered a generalization of the pCN algorithm. Therefore, we first briefly introduce the Metropolis-Hastings algorithm and the pCN algorithm, and then present our new kernel.

The Metropolis-Hastings algorithm is a standard sampling method that can be used to sample from the posterior distribution, thereby solving inverse problems within the Bayesian framework. Additionally, it can serve as the transition kernel in the mutation step of SMC methods. The Metropolis-Hastings algorithm consists of two iterative steps: first, a transition is made from the current state \(u\) according to a proposal kernel \(Q(u,dv)\), and second, the transition is accepted with a certain probability \(a(u,v)\). Alternatively, we can represent the entire process with a formula:
\begin{align}
P(u,dv)=Q(u,dv)a(u,v)+\delta_u(dv)\int_\mathcal{H}(1-a(u,w))Q(u,dw).\label{Q-P}
\end{align}

The key components of this kernel are the proposal \( Q(u, dv) \) and the acceptance rate function \( a(u, v) \). To sample from the target measure \(\rho\), which is absolutely contiounous with respect to prior \(\mu_0\), we must select \( Q(u, dv) \) and  \( a(u, v) \) such that \( P(u, dv) \) defined by (\ref{Q-P}) preserves \( \rho \) as follows:
\begin{equation}
\int_{u \in \mathcal{H}} P(u, dv) \rho(du) = \rho(dv).
\end{equation}
According to Theorem 21 in \cite{Dashti2017}, we can choose \( a(u, v) \) as follows:
\begin{equation}
a_{\text{\tiny MH}} = \min\left\{1, \frac{\rho(dv) Q(v, du)}{\rho(du) Q(u, dv)}\right\},\label{eq:Q-a}
\end{equation}
so that \( P(u, dv) \) preserves \(\rho \), provided that \( \rho(du) Q(u, dv) \) and \( \rho(dv) Q(v, du) \) are equivalent as measures on \( (\mathcal{H} \times \mathcal{H}, B(\mathcal{H}) \times B(\mathcal{H})) \). We will refer to a Metropolis-Hastings algorithm as well-defined if this condition of absolute continuity is satisfied \cite{Cui2016}. Hence, we only need to select the proposal kernel \( Q(u, dv) \) and demonstrate the equivalence between \( \rho(du) Q(u, dv) \) and \( \rho(dv) Q(v, du) \).

It should be noted that only carefully designed Metropolis-Hastings methods have interpretations in infinite dimensions. In fact, most Metropolis-Hastings methods defined in finite dimensions will not make sense in the infinite-dimensional limit because the acceptance probability is defined as the Radon-Nikodym derivative between two measures, and measures in infinite dimensions have a tendency to be mutually singular \cite{Dashti2017}. 

One carefully designed algorithm is the preconditioned Crank-Nicolson (pCN) method, which has the corresponding transition kernel
\begin{equation}
Q_{\text{pCN}}(u, dv) =m+ \sqrt{1 - \beta^2} (u-m) + \beta  \mathcal{N}_{0,\mathcal{C}}(dv),\label{pCNproposal}
\end{equation}
where $\mathcal{N}_{0,\mathcal{C}}$ denotes the Gaussian measure $\mathcal{N}(0,\mathcal{C})$. Note that \eqref{pCNproposal} is a drifted version \cite{Pinski2015} of the pCN proposal. The special case without drift (i.e., when \(m=0\)) is more frequently discussed \cite{Andrew_review,Beskos2015SC}.
This kernel allows for deriving a simple acceptance rate function
\begin{align}
a_\text{pCN}(u,v)=\min\left\{1,\exp(\Phi(u)-\Phi(v))\right\},\quad  \frac{d\rho}{d\mu_0}(u)\propto e^{-\Phi(u)}.\label{a in pCN}
\end{align}
As an MCMC sampling method, the irreducibility of the pCN algorithm is of fundamental importance, as it is closely related to ergodicity, the law of large numbers, and central limit theorems. 
For MCMC methods arising from discretizations of continuous-time dynamical systems, irreducibility may sometimes depend on the specific discretization scheme. A notable example is Hamiltonian Monte Carlo~\cite{Alain2020}.
Due to intrinsic properties of infinite-dimensional Hilbert spaces, the pCN algorithm generally loses irreducibility. 
Nevertheless, ergodicity, the law of large numbers, and central limit theorems still hold. 
Since a detailed discussion of these issues is not the main focus of this paper,  we refer the reader to  \ref{SM-4Kernels} for more details.

The pCN method will serve as the transition kernel for SMC in numerical experiments, constituting an accurate posterior sampling algorithm for comparison. On the other hand, to prove the convergence of our algorithm (in Section \ref{Subsec.pCN-GMM}), we generalize the pCN method to the case with a Gaussian mixture proposal, and prove that it is also well-defined in an infinite-dimensional space.

Recall that each iteration of pCN consists of two steps: proposal defined in (\ref{pCNproposal}), and acceptance according to (\ref{a in pCN}). Recall that the Gaussian prior \(\mathcal{N}(0,\mathcal{C})\) has the eigenpairs \((\lambda_i,\phi_i)\). Assume that the covariance operators \( \{\mathcal{C}_j\}_{j=1}^M \) in the Gaussian mixture have the same eigenfunctions to the prior:
\begin{align}
\mathcal{C}_j\phi_i=\lambda_{ji}\phi_i,\ 1\leq j\leq M.
\end{align}
Moreover, we use a Gaussian mixture proposal
\begin{align}
Q(u,dv)=\sum_{i=1}^Mw_i m_i+\sqrt{1-\beta^2}\sum_{i=1}^Mw_i(u-m_i)+\beta\sum_{i=1}^Mw_i \mathcal{N}_{0,\mathcal{C}_i}(dv),\label{GMM proposal}
\end{align}
where \(\beta\in(0,1]\) is the step size and \(\{w_i\}_{i=1}^M\) are weights satisfying \(\sum_{i=1}^M w_i=1.\)
\begin{remark}\label{remark of pCN-GM}
	This proposal has been carefully constructed, and it will yield a well-defined Metropolis-Hastings algorithm at the end of this subsection. We can obtain this proposal from the Crank-Nicolson discretization of a Langevin system, see \ref{SM-Langevin2pCNGM}. This proposal degenerates to the pCN proposal with a non-zero mean prior \cite{Pinski2015} when \(M=1\), and to the proposal in the Gaussian mixture independence sampler \cite{Feng2018} when \(\beta=1\). We will estimate the parameters of a Gaussian mixture from samples by the method from \cite{Feng2018}, which is established in function space. And the resulting Gaussian mixture is also equivalent to the prior measure \(\mu_0\).
\end{remark}

For simplicity, rewrite the proposal (\ref{GMM proposal}) into the following form:
\begin{align}
Q(u,dv)=\gamma u+(1-\gamma)\sum_{j=1}^Mw_j m_j+\beta \sum_{j=1}^M w_j\mathcal{N}_{0,\mathcal{C}_j}(dv),\label{pCN-GM proposal}
\end{align}
where \(\gamma:=\sqrt{1-\beta^2}\). And the accept rate function can be computed by
\begin{align}\label{auv of pCN-GMM}
a(u,v)=\min\left\{1,\frac{\rho(dv)Q(v,du)}{\rho(du)Q(u,dv)}\right\}.
\end{align}
Here the accept rate function \(a(u,v)\) is well-defined in infinite-dimensional Hilbert space if and only if the Radon-Nikodym derivative exists, that is, \(\rho(dv)Q(v,du)\) is absolutely continuous with respect to \(\rho(du)Q(u,dv)\), as measures on \(\mathcal{H}\times \mathcal{H}\). We will first illustrate the well-definiteness of the measure \(\mu_{0}(dv)Q(v,du)\) and then show their equivalence.
\begin{theorem}\label{mu_0dvPvdu is well defined}
	Choose the prior measure \(\mu_0=\mathcal{N}(0,\mathcal{C})\), and \(Q(u,dv)\) is defined by \eqref{GMM proposal}. Assume that Gaussian components in \eqref{GMM proposal} are equivalent to the prior \(\mu_0\).
	Consider measures \(\mu_{0}(dv)Q(v,du)\) and \(\mu_{0}(du)Q(u,dv)\).
	
	\noindent (a). They are both \(M\)-component Gaussian mixture measures defined on \(\mathcal{H}\times\mathcal
	H\):
	\begin{align}
	\mu_0(du)Q(u,dv)=&\sum_{j=1}^Mw_j\rho_{j}(du, dv),\quad \rho_{j}=\mathcal{N}\left(
	m_{j},\mathcal{V}_{j}
	\right),\label{muduPudv}\\
	\mu_0(dv)Q(v,du)=&\sum_{j=1}^Mw_j\rho_{j}'(du, dv),\quad \rho_{j}'=\mathcal{N}\left(
	m_{j}',\mathcal{V}_{j}'
	\right).\label{mudvPvdu}
	\end{align}
	where the mean functions and covariance operators are as follows:
	\begin{align}
	m_{j}=&[0,(1-\gamma)) m_j],
	&\mathcal{V}_{j}=\left[\begin{array}{cc}
	\mathcal{C} & \gamma \mathcal{C} \\
	\gamma \mathcal{C} & \beta^2\mathcal{C}_j+\gamma^2\mathcal{C}
	\end{array}\right],\\
	m_{j}'=&[(1-\gamma) m_j,0],
	&\mathcal{V}_{j}'=\left[\begin{array}{cc}
	\beta^2\mathcal{C}_j+\gamma^2\mathcal{C}& \gamma \mathcal{C} \\
	\gamma \mathcal{C} & \mathcal{C}
	\end{array}\right].
	\end{align}
	(b). The covariance operators \(\mathcal{V}_{j},\mathcal{V}_{j}'\) are positive, self-adjoint, and of trace class. Thus  \(\mu_{0}(dv)Q(v,du)\) and \(\mu_{0}(du)Q(u,dv)\) are both well defined Gaussian mixture measures on \((\mathcal{H}\times\mathcal
	H,\mathcal{B}(\mathcal{H})\times\mathcal{B}(\mathcal{H}))\). Every component of them is a well defined Gaussian measure.
\end{theorem}
The proof of theorem \ref{mu_0dvPvdu is well defined} is given in \ref{SM-mu_0dvPvdu is well defined}.
Furthermore we can prove the equivalence between \(\mu_{0}(dv)Q(v,du)\) and \(\mu_{0}(du)Q(u, dv).\)
\begin{theorem}\label{HS}
	Under the assumptions of theorem 
	\ref{mu_0dvPvdu is well defined} we have that the operator \(\mathcal{V}_{j}^{-1/2}\mathcal{V}_{j}'\mathcal{V}_{j}^{-1/2}-\mathcal{I}\) is a Hilbert-Schmidt operator, and \(\mu_0(du)\)
	\(Q(u,dv)\), \(\mu_0(dv)Q(v,du)\) are equivalent. In fact, all of their components are equivalent. Furthermore, \(\rho(du)Q(u,dv)\) and \(\rho(dv)Q(v,du)\) are equivalent. 
\end{theorem}

\begin{remark}
	The equivalence between the two measures in this theorem is crucial (see the proof in \ref{SM-well-defined pCN-GM}), directly affecting whether the algorithm possesses dimension independence. A simple example is the comparison between the pCN and the random walk transition kernel:
	\begin{align*}Q_{\text{\tiny pCN}}(u,dv)=\sqrt{1-\beta^2}u+\beta\mathcal{N}_{0,\mathcal{C}}(dv),\quad Q_{\text{\tiny RW}}(u,dv)=u+\beta\mathcal{N}_{0,\mathcal{C}}(dv),\end{align*}
	where just a slight difference in coefficients can lead to significant changes in the algorithm's acceptance rate as the discrete dimension increases \cite{Andrew_review}.
	Moreover, if we require a proposal of the form given in Equation \eqref{pCN-GM proposal} without the assumption \(\gamma = \sqrt{1 - \beta^2}\) to yield a well-defined acceptance rate function, then the only option is \(\beta^2 + \gamma^2 = 1\). Discarding the meaningless negative solutions, the generalization we provide is unique. The proof of  uniqueness is given in \ref{SM-UniquenessOfpCNGM}.
\end{remark}

The last part of the pCN-GM method is to compute the accept rate function by (\ref{auv of pCN-GMM}). It follows from 
\begin{align}\label{auv of pCN-GMM2}
\begin{split}
\frac{\rho(dv)Q(v,du)}{\rho(du)Q(u,dv)}
=&\frac{\rho(dv)Q(v,du)}{\mu_{0}(dv)Q(v,du)}
\frac{\mu_{0}(dv)Q(v,du)}{\mu_{0}(du)Q(u,dv)}
\frac{\mu_{0}(du)Q(u,dv)}{\rho(du)Q(u,dv)},\\
=&\frac{d\rho/d\mu_0(v)}{d\rho/d\mu_0(u)}\frac{\mu_{0}(dv)Q(v,du)}{\mu_{0}(du)Q(u,dv)}
\end{split}
\end{align}
that we only need to compute the Radon-Nikodym derivative \(\frac{\mu(du)Q(u,dv)}{\mu(dv)Q(v,du)}\). According to Theorem  \ref{HS} this Radon-Nikodym derivative exists. With Theorem \ref{mu_0dvPvdu is well defined} we can write \(\mu(du)Q(u,dv)\) and \(\mu(dv)Q(v,du)\) into Gaussian mixture forms:
\begin{align}
\mu_0(du)Q(u,dv)=\sum_{j=1}^M w_{j}\rho_{j}(du, dv),\quad
\mu_0(dv)Q(v,du)=\sum_{j=1}^M w_{j}\rho_{j}'(du, dv).
\end{align}
All Gaussian measure components are equivalent, hence the Radon-Nikodym derivative is
\begin{align}
\frac{\mu_{0}(dv)Q(v,du)}{\mu_{0}(du)Q(u,dv)}=\sum_{j_2=1}^M\frac{w_{j_2}'}{\sum_{j_1=1}^Mw_{j_1} \frac{d\mu_{j_1}}{d\mu_{j_2}'}}.\label{auv of pCN-GMM3}
\end{align}
Using the Radon-Nikodym derivative between Gaussian measures \cite{Kukush2019}, we have
\begin{align}
\frac{d\mu_{j_1}}{d\mu_{j_2}'}[u,v]
=\frac{\exp\left\{
	-\frac{1}{2}\|[u,v]-m_{j_1}\|_{\mathcal{V}_{j_1}}^2
	+\frac{1}{2}\|[u,v]-m_{j_2}'\|_{\mathcal{V}_{j_2}'}^2
	\right\} }{\sqrt{\det(\mathcal{V}_{j_2}'^{-1/2}\mathcal{V}_{j_1}\mathcal{V}_{j_2}'^{-1/2})}},\label{auv of pCN-GMM4}
\end{align}
where we use the symbol \(\|\cdot\|_\mathcal{V}^2=\|\mathcal{V}^{-1/2}\cdot\|^2\). 

\begin{algorithm}
	\caption{pCN based on Gaussian mixtures}
	\label{alg pCN GM}
	\begin{algorithmic}[1]
		\STATE{\textbf{Input}: Specify the prior \(\mu_0\). Draw \(u_0\sim\mu_0.\)}
		\STATE{\textbf{For} \(i=0,1,\cdots ,I\), \textbf{do}:}
		\STATE{\qquad Sample \(v_i\) from the proposal defined in (\ref{GMM proposal});}
		\STATE{\qquad Accept \(v_i\) with probability  \(a(u,v)\) computed by (\ref{auv of pCN-GMM}),(\ref{auv of pCN-GMM2}),(\ref{auv of pCN-GMM3}),(\ref{auv of pCN-GMM4});}
		\STATE{\qquad Otherwise, keep the current state, \( u^{(k+1)} = u^{(k)} \).}
		\STATE{\textbf{Output}:  The final state \(u_{I}\).}
	\end{algorithmic}
\end{algorithm}

Finally we have a new Markov transition kernel, which is well-defined in separable Hilbert space. The significance of this algorithm is reflected in the following three aspects:
\begin{itemize}
	\item It maintains the property of the pCN algorithm being invariant to the posterior; hence, when used as a transition kernel in SMC methods, it can theoretically achieve an approximation of the posterior with arbitrary precision according to Theorem \ref{SMC approximation theorem}.
	\item Disregarding the relationship with SMC, when measurements are sparse, the posterior often exhibits multi-modality. At such times, this algorithm can incorporate multi-modal information into the proposal.
	\item Direct sampling from the Gaussian mixture, is a part of the pCN-GM method under specific parameters. Thus, if we can find a good Gaussian mixture approximation of the posterior, we can use the Gaussian mixture sampler to replace the whole mutation step, which is the intuitive idea of the next subsection.
\end{itemize}

Using Algorithm~\ref{alg pCN GM} as the mutation step in Algorithm~\ref{alg SMC} yields the SMC-pCN-GM algorithm. 
The approximation properties of SMC-pCN-GM are established in Theorem~\ref{SMC approximation theorem}, while its well-definedness is ensured by Theorem~\ref{HS}.

\subsection{Approximation theory of Gaussian mixtures}\label{Subsec.GM}
In SMC, we need to choose a transition kernel in the mutation step, and the simplest choice is to sample from the prior, i.e., to use an independence sampler. This sampler method works well when the likelihood is not too informative, but works poorly if the information in the likelihood is substantial \cite{Dashti2017}. To utilize the information from the likelihood, several proposals are available, for example, the pCN, the pCN Langevin \cite{Andrew_review}, and the Hessian-preconditioned explicit Langevin \cite{Cui2016}, each utilizing information of \( \Phi \), \( D\Phi \), and \( D^2\Phi \) respectively. We need to solve one, two, and four PDEs, respectively, to compute them with the adjoint method \cite{Ghattas2021}, which is computationally intensive. To utilize likelihood information while avoiding the computation of the likelihood function, we use a Gaussian mixture approximation of the posterior. In other words, we retain the first step of pCN-GM, which proposes a sample, and discard the second step of the acceptance process, accepting all proposals directly. The resulting SMC method is SMC-GM, because the mutation step is simply a Gaussian mixture sampler. SMC-GM makes sense only if we can find a good Gaussian mixture approximation to the posterior. In this subsection, we prove the denseness of the Gaussian mixture measure in the posterior measures, provide an error estimate, and propose the SMC-GM method finally.

Note that under the finite-dimensional setting, a Gaussian mixture density function can approximate any density function. Define the set of \(m\)-component location-scale linear combinations of the probability density function \( g \) as
\begin{align*}\mathcal{M}_M^g=\left\{h:h(x)=\sum_{i=1}^M\frac{c_i}{\sigma_i^n}g\left(\frac{x-\mu_i}{\sigma_i}\right),\mu_i\in \mathbb{R}^n,\sigma_i,c_i> 0,1\leq i\leq M,\sum_{i=1}^Mc_i=1\right\}.\end{align*}
Define the set of continuous functions that vanish at infinity by
\begin{align*}\mathcal{C}_0=\left\{f \text{ is continuous}:\forall \epsilon>0,\exists \text{ a compact } K\subset \mathbb{R}^n, \text{ such that } |f|<\epsilon,\text{ for any } x\notin K\right\}.\end{align*}
Then it follows from Theorem 5.(d) in \cite{Nguyen2020} that Gaussian mixtures have the universal approximation property in the finite-dimensional setting. If \( f \) and \( g \) are probability density functions and \( g \in \mathcal{C}_0 \), then for any measurable \( f \), there exists a sequence \( \{ h_M \} \in \mathcal{M}_M^g \) such that
\begin{align*} \lim_{M \to \infty} h_M = f, \text{ a.e.} \end{align*}
We will extend this approximation property to the infinite-dimensional setting. First, we define the Gaussian mixture measure. We say a probability measure is a mixture of Gaussians if it can be expressed as a weighted finite sum of Gaussian measures, i.e.,
\begin{align*}
\mu(du) = \sum_{i=1}^M w_i \mu_i(du),
\end{align*}
where \( w_i > 0 \) and \( \mu_i \) are well-defined Gaussian measures on \( \mathcal{H} \), \( 1 \leq i \leq m \). From this definition, it is directly obtained that \( \sum_{i=1}^m w_i = 1 \). To express our approximation property, we use the distance defined in \eqref{distance}. Now we can prove that Gaussian mixtures are dense in posterior measures, i.e., the universal approximation property of Gaussian mixtures in the infinite-dimensional setting.
\begin{theorem}
	Given a well-defined Gaussian prior \( \mu_0 = \mathcal{N}(0, \mathcal{C}) \), a continuous potential function \( \Phi\), and the posterior measure determined by Bayes' formula
	\begin{align*} \frac{d\mu^\mathrm{d}}{d\mu_0}(u) = \frac{1}{Z} \exp\left(-\Phi(u)\right) .\end{align*}
	Then for every \( \epsilon > 0 \), there exists \( M\in \mathbb{N} \) and \( \{ w_i, m_i, \mathcal{C}_i \}_{i=1}^M \), forming a well-defined Gaussian mixture measure in infinite-dimensional space:
	\begin{align*} \tilde{\mu}(du) = \sum_{i=1}^M w_i \mu_i(du),\quad \mu_i=\mathcal{N}(m_i, \mathcal{C}_i) ,\end{align*}
	such that \( d_{\text{\tiny TV}}(\tilde{\mu}, \mu^\mathrm{d}) < \epsilon \). Moreover, the Gaussian components \(\{\mu_i\}_{i=1}^M\) are all equivalent to the prior \( \mu_0 \).\label{GMM approximate post}
\end{theorem}

\begin{remark}
	The above theorem states that for posterior measures in a separable Hilbert space, Gaussian mixture measures are dense within them; furthermore, this Gaussian mixture is also mutually equivalent to the prior measure. This theorem is the motivation for our use of Gaussian mixtures for approximation. The proof is given in \ref{SM-GMM approximate post}.
\end{remark}

In practice, we estimate the parameters of a Gaussian mixture from samples using the method from \cite{Feng2018}, {\color{black} in which the Gaussian mixture is formulated as a probability measure equivalent to the prior measure $\mu_0$ on an infinite-dimensional function space.} Then we show that sampling from Gaussian mixture is an approximation of pCN-GM. Let \(\beta=1\) in (\ref{GMM proposal}), and we obtain
\begin{align}
Q(u,dv)=&\sum_{i=1}^Iw_i \mathcal{N}_{m_i,\mathcal{C}_i}(dv).\label{GMM beta=1}
\end{align}
Now \(Q(u,dv)\) is independent with respect to \(u\), thus we can drop \(u\) and rewrite it as \(Q(dv)\). The accept rate function is
\begin{align}
a(u,v)=\min\left\{1,\frac{e^{-\Phi(v)}\mu_{0}(dv)Q(du)}{e^{-\Phi(u)}\mu_{0}(du)Q(dv)}\right\}.
\end{align}
Then with the help of Theorem \ref{GMM approximate post} we have the equivalence between \(Q(du)\) and \(\mu_0(du)\), therefore \(a(u,v)\) has a simpler form:
\begin{align}
a(u,v)=\min\left\{1,
\frac{e^{-\Phi(v)}\mu_{0}(dv)}{Q(dv)}
\frac{Q(du)}{e^{-\Phi(u)}\mu_{0}(du)}
\right\}.
\end{align}
If we have a good Gaussian mixture approximation of the posterior, then \(a(u,v)\approx 1\), which means sampling from the Gaussian mixture is an approximation of pCN with a Gaussian mixture proposal shown in Algorithm \ref{alg pCN GM}. More precisely, recall that \(Q(dv)\) is the Gaussian mixture sampler defined by (\ref{GMM beta=1}),  and we have:
\begin{theorem}\label{theorem4}
Recall that \(\mu_j\) is the posterior of layer j. Then we have that the Gaussian mixture sampler  \(Q_{\mu_j^N}\) is an approximation of the transition kernel \(P\), defined by (\ref{Q-P}): 
\begin{align}
d(Q_{\mu_j^N},P\mu_j^N)\leq d(Q_{\mu_j^N},\mu_{j})+d(\mu_j^N,\mu_{j}).  \label{GMneq1}
\end{align}
Moreover there exists a Gaussian mixture measure $Q_{\mu_j^N}$ such that
\begin{align}
	d'(Q_{\mu_j^N},P\mu_j^N)\leq (2+\epsilon)d'(\mu_j^N,\mu_{j}), \text{ for any }\epsilon>0,\label{GMneq2}
\end{align}
where \(d'\) is a weaker distance defined by
\begin{align*}
d'(\mu,\nu):=\sup_{|f|_\infty\le1,\|D^2f\|\leq 1}\sqrt{
	\mathbb E\Big[
	\int fd\mu-\int fd\nu
	\Big]^2}.
\end{align*}
\end{theorem} 
Note that the original error bound still holds, since $d'(\mu,\nu) \leq d(\mu,\nu)$. The proof of Theorem \ref{theorem4} is given in \ref{SM-ErrOfGMM}.
This theorem indicates that, if we have a good Dirac approximation of the posterior, then it is reasonable to directly use a Gaussian mixture sampler to replace the transition kernel \(P\). In general, this approximation introduces a larger error than standard MCMC methods. The significance of this theorem is that the additional error remains controllable, which is also reflected in the numerical results presented in the next section. 
The resulting SMC-GM is concluded in Algorithm \ref{alg SMC-GM}.

\begin{algorithm}
	\caption{Sequential Monte Carlo with Gaussian mixture approximation}
	\label{alg SMC-GM}
	\begin{algorithmic}[1]
		\STATE{\textbf{Input}: Specify \(\mu_0\). Draw \(\{v_0^{(n)}\}_{n=1}^N\sim \mu_0.\)}
		\STATE{\textbf{For} \(j=0,1,\cdots ,J\), \textbf{do}:}
		\STATE{\qquad Set \(w_j^{(n)}=1/N,n=1,2,\cdots, N\) and define \(\mu_j^N=\sum_{n=1}^Nw^{(n)}_{j}\delta_{v^{(n)}_{j}}\) ;}
		
		\STATE{\qquad (Re-weighting) Update \(w_{j+1}^{(n)}\) by:}
		\begin{align*}
		\hat{w}_{j+1}^{(n)}=e^{-h_{j+1}\Phi(v_{j}^{(n)})}w_j^{(n)},\qquad w_{j+1}^{(n)}=\frac{\hat{w}_{j+1}^{(n)}}{\sum_{n=1}^N\hat{w}_{j+1}^{(n)}}
		\end{align*}
		\STATE{\qquad (Resampling) Update \(\hat{v}_{j+1}^{(n)}\) by resampling from  \(\sum_{n=1}^Nw^{(n)}_{j+1}
			\delta_{v^{(n)}_{j}}\).}
        \STATE{\qquad (Mutation) Estimate a Gaussian mixture \(\tilde{\mu}\) by samples, and draw \begin{align*}\{v_{j+1}^{(n)}\}_{n=1}^N\sim \tilde{\mu};\end{align*}}

		\STATE{\textbf{Output}: \(\mu_J^N\).}
	\end{algorithmic}
\end{algorithm}
The difference between SMC and SMC-GM is that in SMC, the mutation step must be a transition kernel preserving \(\mu_j\), the posterior at the \(j\)-th layer, whereas in SMC-GM, this step involves sampling from a Gaussian mixture measure determined by the samples.

\section{Numerical Examples}\label{Sec.Numerics}
In this section, we provide three numerical examples to substantiate the results in Section \ref{Sec.SMC-GM}. Compared with SMC-pCN-GM, the reduced SMC-GM algorithm requires fewer computational resources while maintaining satisfactory performance. Hence, we mainly report the numerical results of the SMC-GM algorithm.
We apply the SMC-GM method to two typical inverse problems, including one linear and one nonlinear problem. For the nonlinear problem, we also consider a sparse measurement situation as the third example. Through these experiments, we demonstrate that the SMC-GM algorithm is faster than SMC-pCN and has strong sampling capabilities for multi-modal posteriors.
In the following, SMC-pCN is used for comparison. Thus we briefly discuss the computational cost of them. If we require \(N\) samples, and the number of intermediate measures is \(J\), then SMC-GM method only need to compute the weights in the resampling step for each particle, and a total of \(NJ\) PDEs. As for SMC-pCN method, we need to sample a Markov chain for every particle in each layer, the length of the chain can be chosen between 5 and 1000 \cite{Beskos2015SC}. In all three numerical examples, we will select Markov chains of length 200, requiring us to solve a total of \(200NJ\) PDEs. Note that the method for determining \(h_i\) is described in section \ref{SM-NumericalDetails} in the supplementart material and we estimate the parameters of a Gaussian mixture from samples using the method from \cite{Feng2018}. All the programming codes are accessible at \url{https://github.com/jjx323/IPBayesML-SMC-GM}. 
SMC is intrinsically parallel, and so are our implementations.

\subsection{A multi-modal example}
As the first example, we construct a multi-modal problem in a 1-dimensional interval \(\Omega = [0,1]\) as described in \cite{Feng2018}. Let the unknown \(u \in \mathcal{H} = L^2(\Omega)\) and assume the prior is a zero-mean Gaussian measure on \(\mathcal{H}\) with covariance operator \(\mathcal{C} = (I - 0.01 \Delta)^{-2}\). We consider a four-modal likelihood function given by
\begin{align}
\exp(-\Phi(u)) = \sum_{i=1}^4 \exp\left(-\frac{1}{2\sigma^2} \| u - f_i(x)  \|^2\right),  \label{example1phi}
\end{align}
and the modals are chosen as follows:
\begin{align}
f_1 = \cos(\pi x), \quad f_2 = -\cos(\pi x), \quad f_3 = \cos(2\pi x), \quad f_4 = \cos(3\pi x).
\end{align}
It can be verified that the likelihood defined in \eqref{example1phi} satisfies the Assumptions 6.1 in \cite{Andrew_review}. It is easy to see that the posterior distribution should have four modes, i.e., \(\{f_i\}_{i=1}^4\).

\begin{figure}[ht]
	\centering
	\captionsetup[subfigure]{skip=-0.5pt}
	\captionsetup[subfigure]{labelformat=empty}
	\captionsetup{aboveskip=0.1pt} 
	
	\subfloat[\hspace{2.2em}(a) Samples using SMC-GM]{
		\includegraphics[ keepaspectratio=true, width=0.46\textwidth, clip=true]{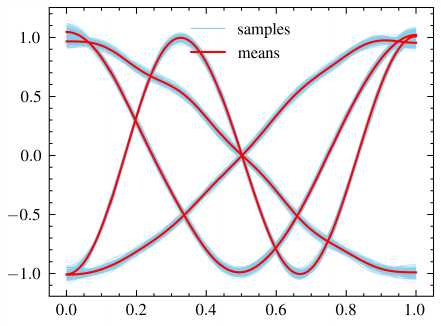}}
	\subfloat[\hspace{2.2em}(b) Samples using SMC-pCN]{
		\includegraphics[ keepaspectratio=true, width=0.46\textwidth, clip=true]{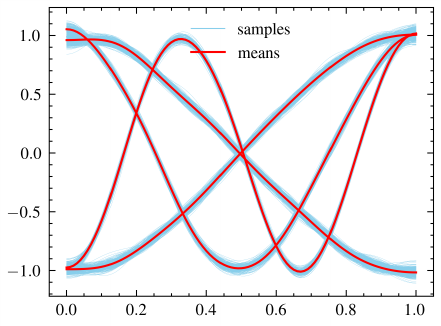}}\\
	\vskip 0.3 cm
	\caption{\emph{\small 
			(a): Posterior samples obtained from the SMC-GM method. (b): Posterior samples obtained from the SMC-pCN method.}}	
	
	\label{fig_ex1_1}
\end{figure}

We drew \(2 \times 10^4\) samples from the posterior of \(u\) using SMC-GM method, and SMC-pCN method for comparison. According to the illustrations given in \cite{Beskos2015SC}, the length of the pCN chain was chosen as \(200\), and we show \(1000\) samples, which are clustered using K-means, in Figure \ref{fig_ex1_1}. The samples are represented by sky-blue dashed lines, and the mean of each cluster is depicted with a red solid line.

We observe no discernible difference between the posterior samples obtained by SMC-pCN and SMC-GM, indicating that our algorithm approximates the posterior well. SMC-GM has been verified for accuracy and is satisfactory. In terms of efficiency, SMC-GM obtained 20,000 samples in 270 seconds, averaging 74 samples per second. In contrast, SMC-pCN required 5774 seconds, averaging 4 samples per second. This clearly demonstrates the significant advantage of SMC-GM in terms of sampling efficiency, as its efficiency is 21 times that of the SMC-pCN method.

\subsection{Elliptic PDE Inverse Problem}\label{example2}
\begin{figure}[ht]
	\centering
	\captionsetup[subfigure]{skip=3pt}
	\captionsetup[subfigure]{labelformat=empty}
	\captionsetup{aboveskip=0.1pt} 
	
	\subfloat[\hspace{-1em}(a) Background truth]{
		\includegraphics[width=0.31\linewidth]{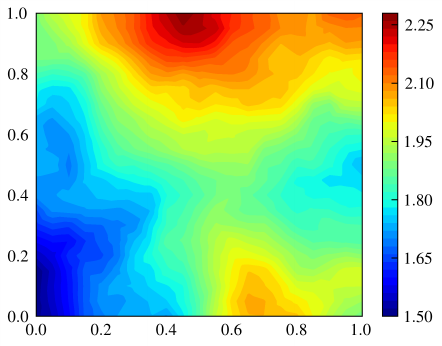}
	}
	\subfloat[\hspace{-0.5em}(b) Mean of SMC-GM]{
		\includegraphics[width=0.3\linewidth]{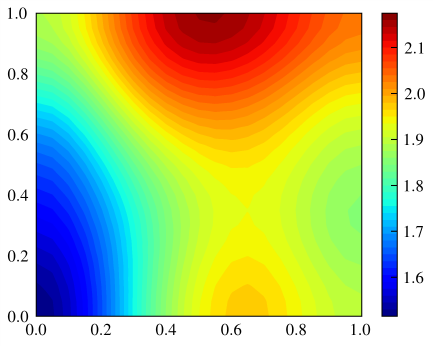}
	}
	\subfloat[\hspace{-1em}(c) Mean of SMC-pCN]{
		\includegraphics[width=0.3\linewidth]{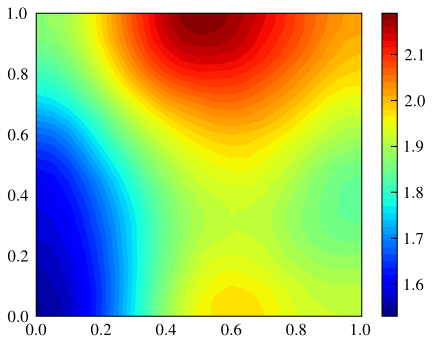}
	}
	\vskip 0.3 cm
	\caption{\emph{\small 
			(a): The background truth of  \(u\). (b): The posterior mean estimated using SMC-GM. (c): The posterior mean estimated using SMC-pCN.}}	
	\label{fig:example1_mean}
\end{figure}

\paragraph{Basic setting}
As the second example, we consider a canonical nonlinear inverse problem involving inference of the diffusion coefficient in the following steady-state Darcy flow equation \cite{Dashti2011}:
\begin{align}\label{example2:PDE}
\begin{split}
-\nabla \cdot (e^{u(x)} w(x)) &= f(x), \quad x \in \Omega, \\
w(x) &= 0, \quad x \in \partial \Omega, 
\end{split}
\end{align}
where \(\Omega = (0,1)^2\), \(f \in H^{-1}(\Omega)\) is the known source function, and \(u \in L^\infty(\Omega)\) denotes the log-permeability, which is our parameter of interest. Let \(\boldsymbol{x}_1, \ldots, \boldsymbol{x}_N\) be the measurement locations. Then the forward operator has the following form:
\begin{align}
\mathcal{F}u = (w(\boldsymbol{x}_1), \ldots, w(\boldsymbol{x}_N))^T.
\end{align}
To illustrate the effectiveness of the SMC-GM method, we compare it with SMC-pCN. For clarity, we list the common parameters of the two algorithms as follows:
\begin{itemize}
	\item Let the domain \(\Omega\) be \((0,1)^2\), and the measurement points \(\{x_i\}_{i=1}^{100}\) are taken at the coordinates \(\{(\frac{9+98i}{900},\frac{9+98j}{900})\}_{i,j=0}^{9}\). To avoid the inverse crime \cite{Kaipio2005}, the data is generated on a fine mesh with the number of grid points equal to \(500 \times 500\). We use a \(20 \times 20\) mesh in the inverse stage.
	\item We use a Gaussian prior \(\mathcal{N}(0,\mathcal{C}_0)\), where the covariance operator is given by \(\mathcal{C}_0 = (I - \Delta)^{-2}\). The source function in \eqref{example2:PDE} is set as a constant function \(f = 1\). We assume the data are produced from an underlying true signal sampled from the prior, which is plotted in Figure \ref{fig:example1_mean}(a).
	\item Assume that 2\% random Gaussian noise \(\epsilon \sim N(0, \Gamma_{\text{noise}})\) is added, where \(\Gamma_{\text{noise}} = \tau^{-2}I\) and \(\tau^{-1} = 0.02 \max_i |(\mathcal{F}u)_i|\). The particle number in SMC is set to \(N = 1000\).
\end{itemize}

Furthermore, SMC-pCN has its unique parameter settings as follows: 
\begin{itemize}
	\item The initial step norm is \(\beta = 0.2\), then we update it using an adaptive strategy (see \ref{SM-NumericalDetails}).
	\item In \cite{Beskos2015SC}, it is recommended that the number of transition steps for the transfer kernel in SMC be between 5 and 1000. Here, the Markov chain length is set to \(l_{\text{pCN}} = 200\), significantly lower than usual pCN sampling processes.
\end{itemize}

\begin{figure}[ht]
	\centering
	\captionsetup[subfigure]{skip=3pt} 
	\captionsetup[subfigure]{labelformat=empty}
	\captionsetup{aboveskip=0.1pt} 
	\subfloat[\hspace{1em}(a) Density of \(u_1\).]{
		\includegraphics[width=0.3\linewidth]{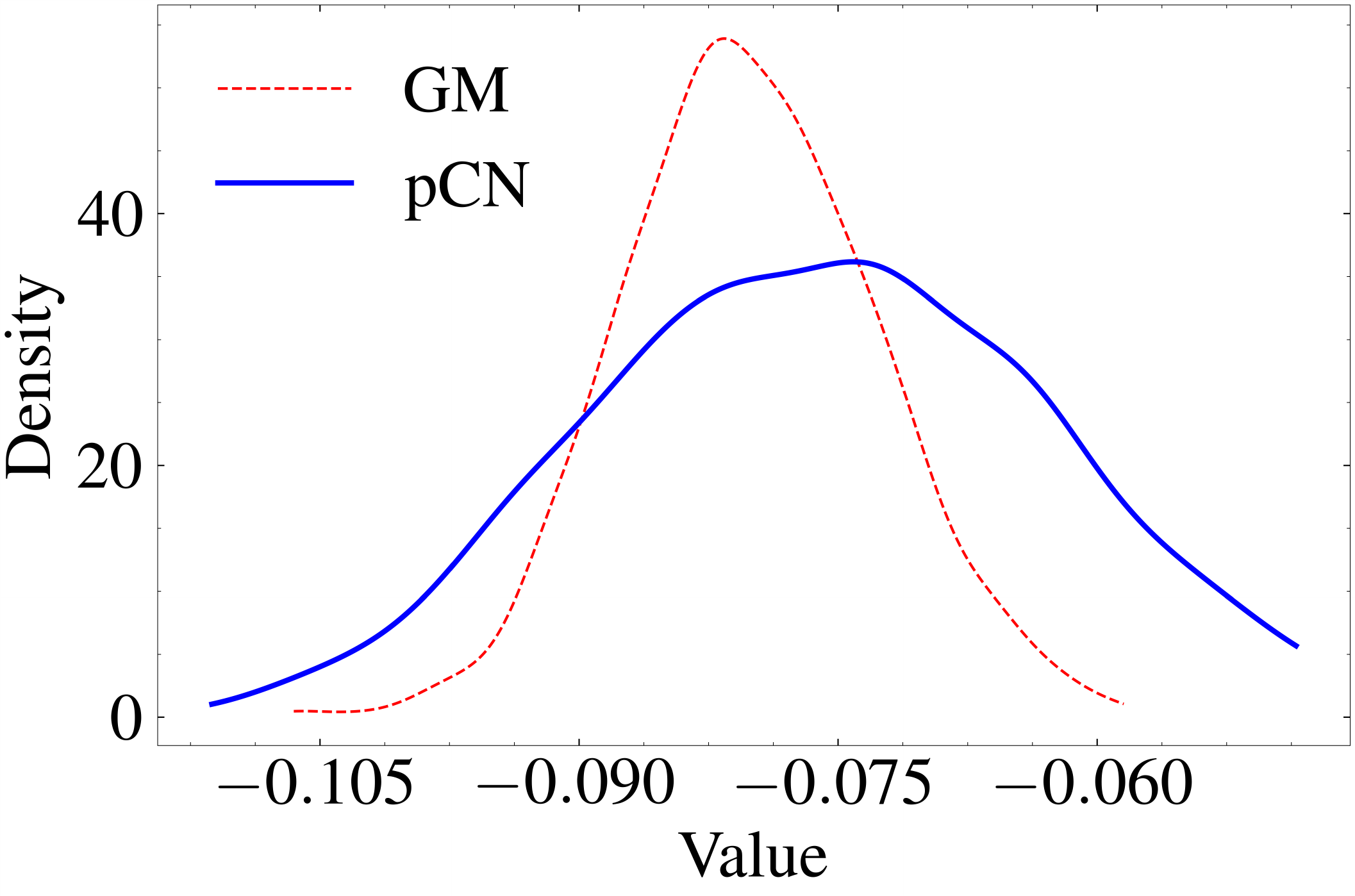}
	}
	\subfloat[\hspace{1em}(b) Density of \(u_4\).]{
		\includegraphics[width=0.3\linewidth]{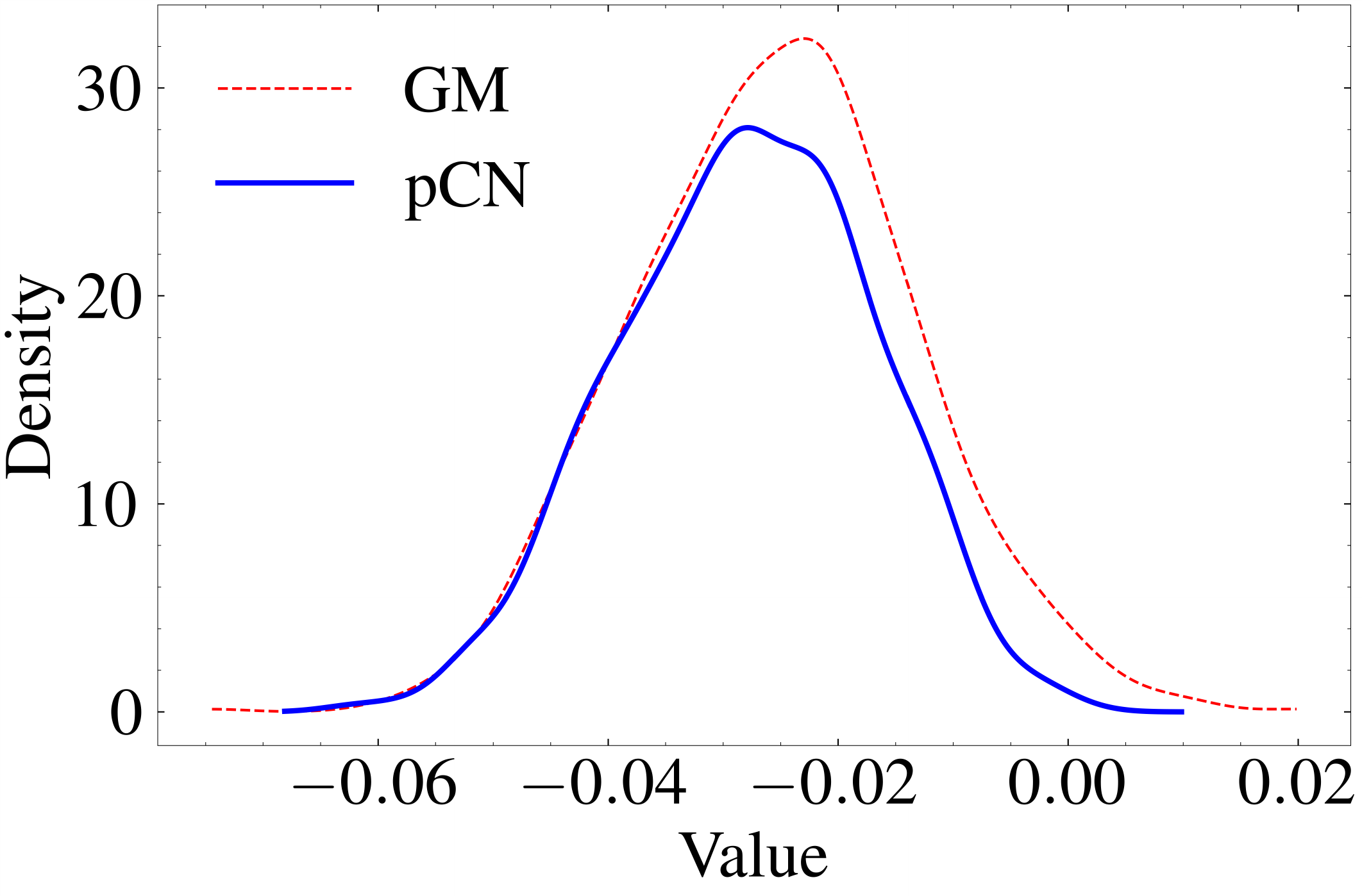}
	}
	\subfloat[\hspace{1em}(c) Density of \(u_7\).]{
		\includegraphics[width=0.3\linewidth]{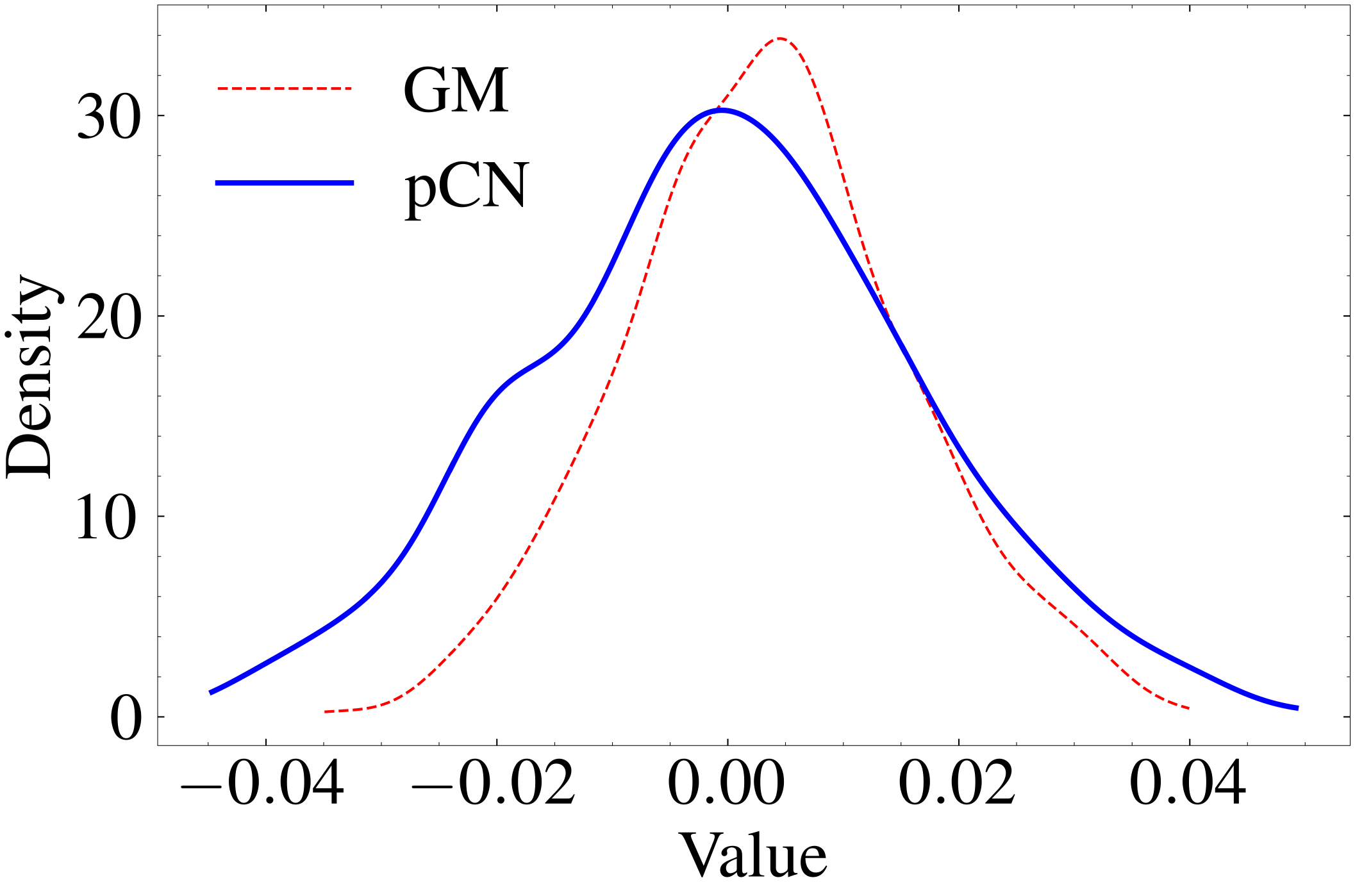}
	}\vspace{10pt}
	
	\subfloat[\hspace{1em}(d) Density of \(u_{10}\).]{
		\includegraphics[width=0.3\linewidth]{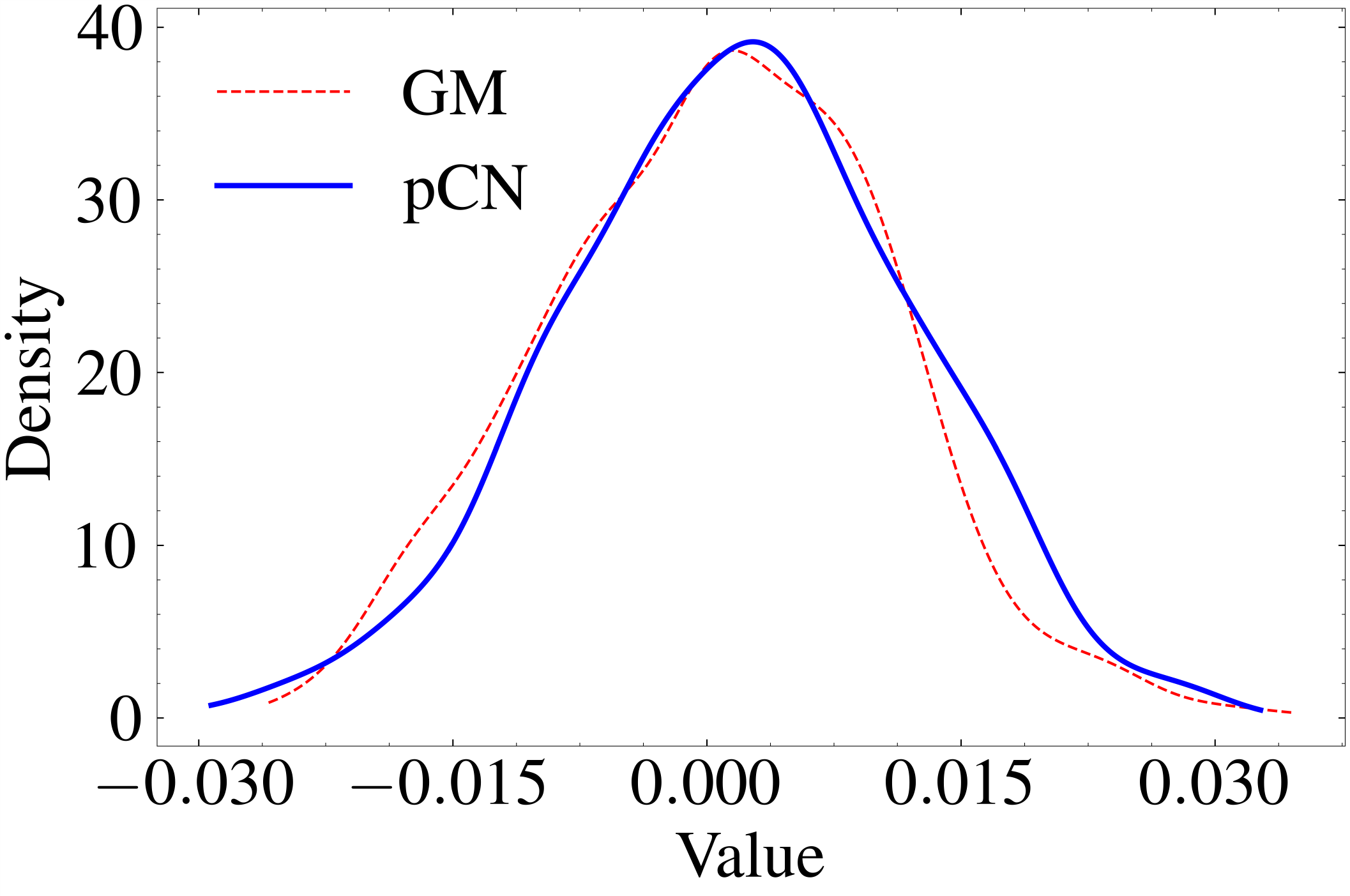}
	}    
	\subfloat[\hspace{1em}(e) Density of \(u_{13}\).]{
		\includegraphics[width=0.3\linewidth]{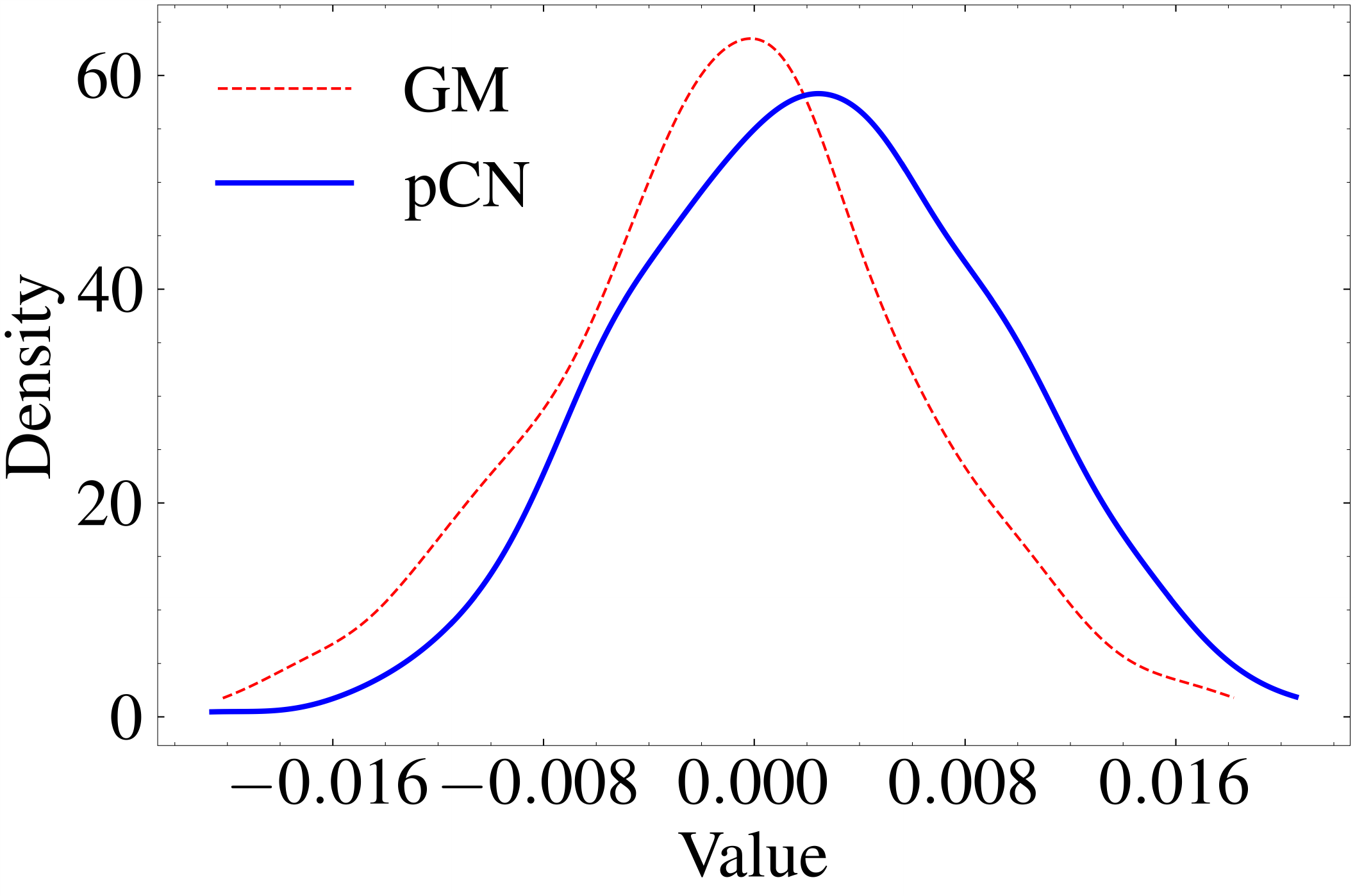}
	}    
	\subfloat[\hspace{1em}(f) Density of \(u_{16}\).]{
		\includegraphics[width=0.3\linewidth]{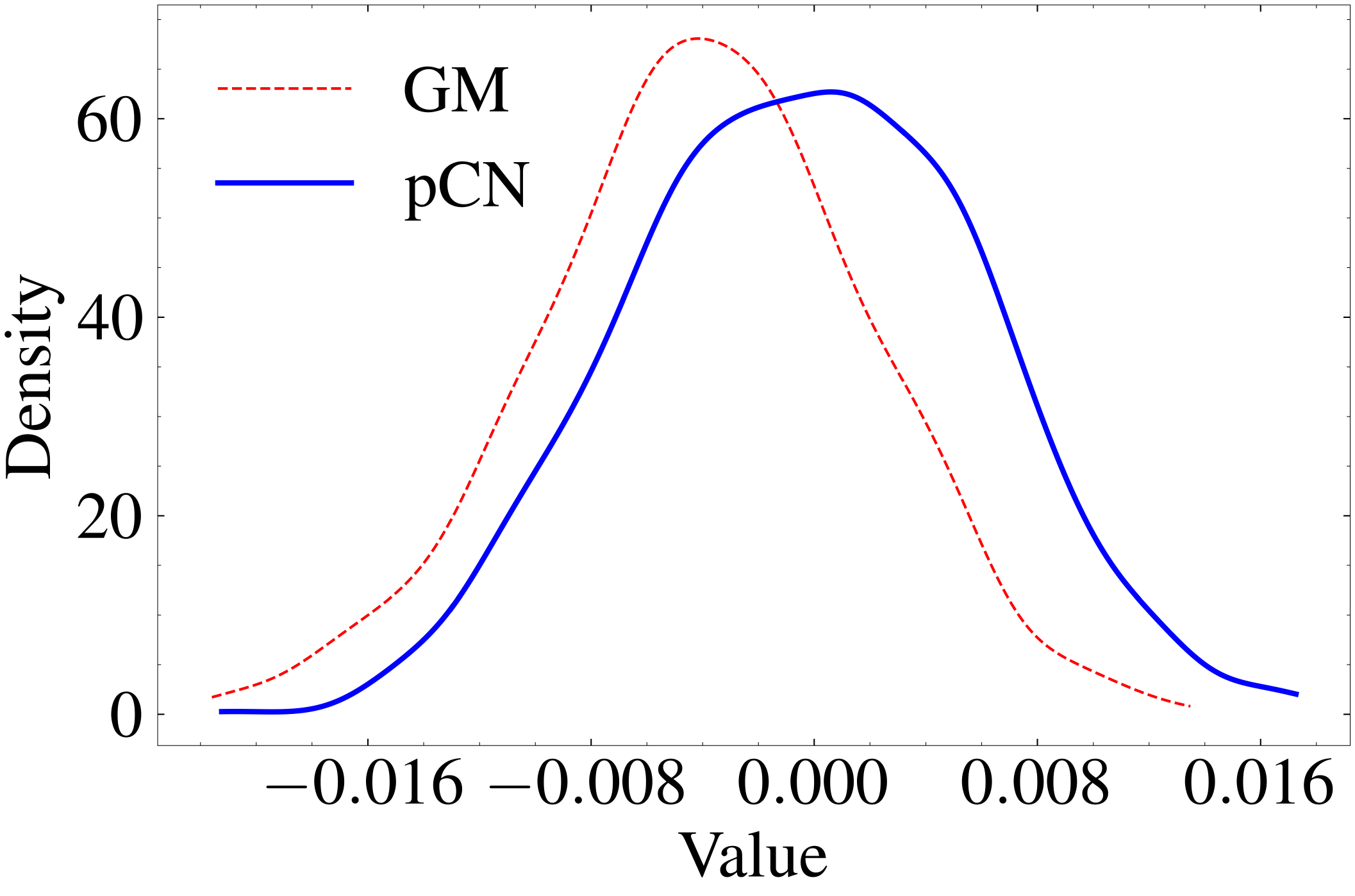}
	}
	\vskip 0.3 cm
	\caption{\emph{\small 
			Posterior densities of the 1st, 4th, 7th, 10th, 13th, 16th Fourier coefficients.
	}}	
	\label{fig:DensityOfCoef}
\end{figure}

\begin{figure}[ht]
	\centering
	\captionsetup[subfigure]{skip=3pt}
	\captionsetup[subfigure]{labelformat=empty}
	\captionsetup{aboveskip=0.1pt} 
	\subfloat[\hspace{3em}(a) Sum of \(h_j\) in SMC-GM]{
		\includegraphics[width=0.4\linewidth]{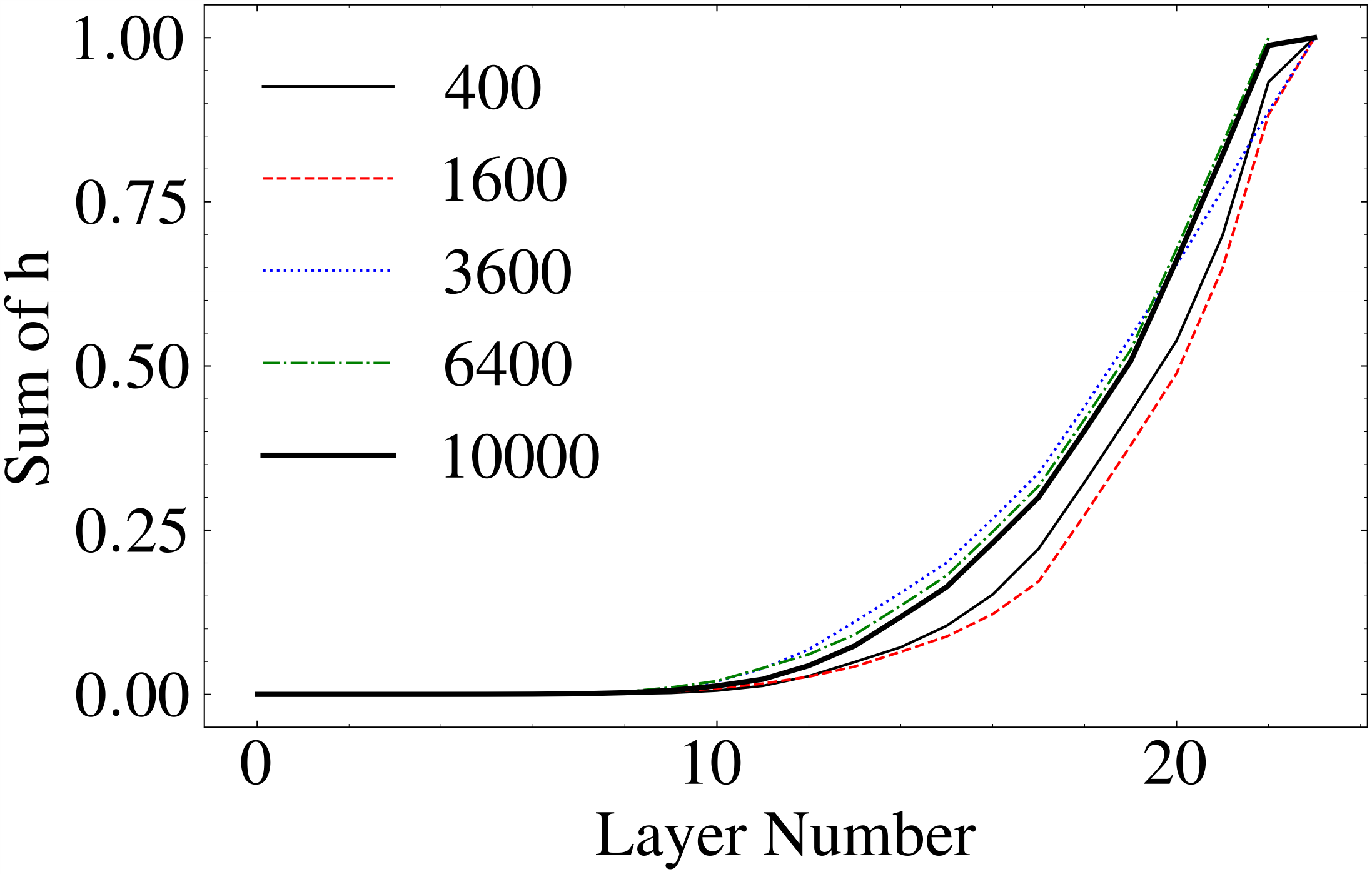}
	}\hspace{2em}
	\subfloat[\hspace{3em}(b) Sum of \(h_j\) in SMC-RW]{
		\includegraphics[width=0.4\linewidth]{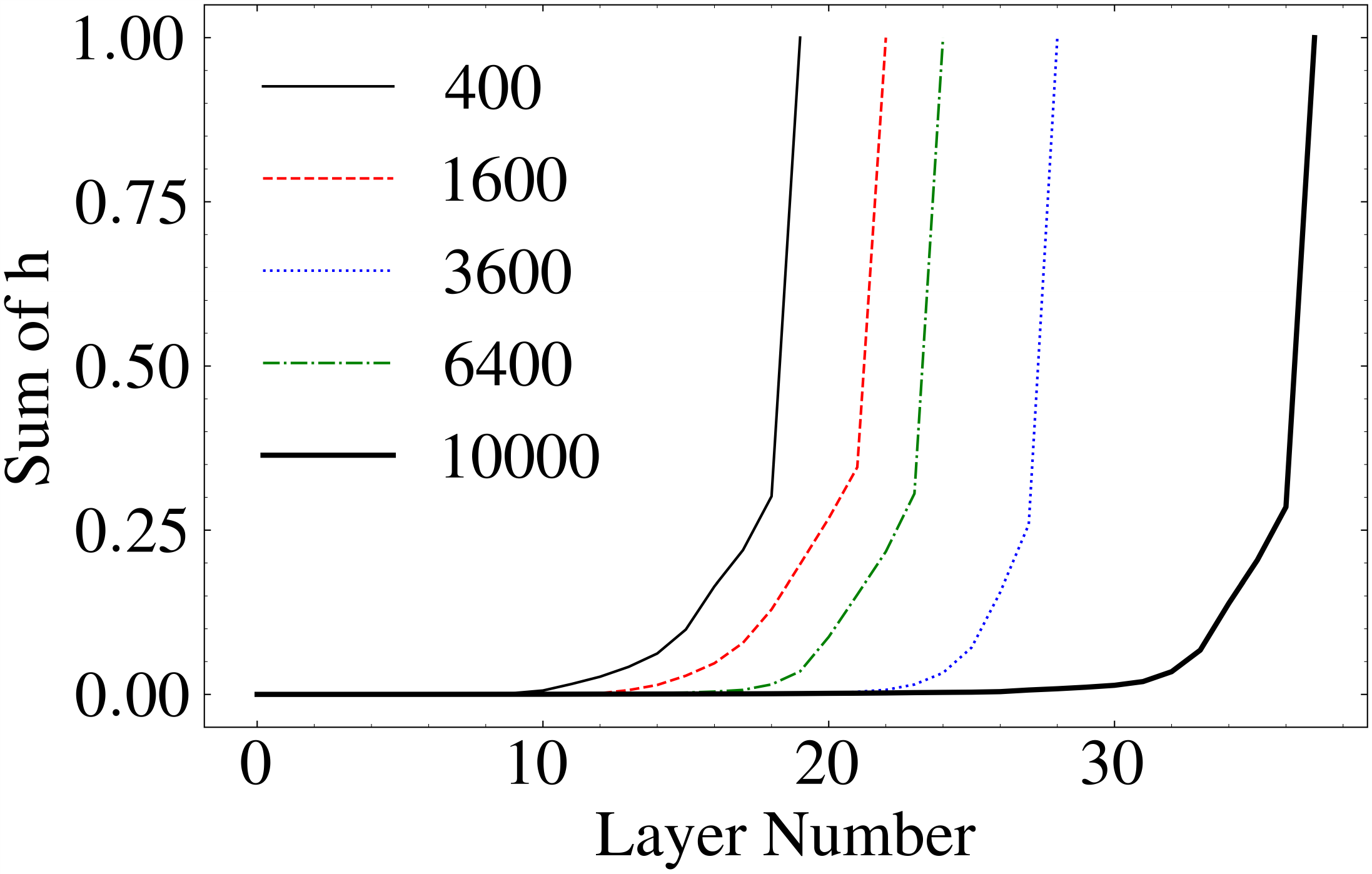}
	}
	\vskip 0.3 cm
	\caption{\emph{\small 
			In the SMC method, we always need to find a sequence \(\{h_j\}_{j=1}^J\) that ranges from 0 to 1  to determine \(\{\mu_j\}_{j=1}^J\), the intermediate measures. (a): \(\{h_j\}_{j=1}^J\) found by the SMC-GM algorithm under different dimensional grid discretizations. (b): \(\{h_j\}_{j=1}^J\) found by SMC-RW under different dimensional grid discretizations.}}	
	\label{fig_MeshIdpdt}
\end{figure}

\begin{figure}[ht]
	\centering
	\captionsetup[subfigure]{skip=3pt} 
	\captionsetup[subfigure]{labelformat=empty}
	\captionsetup{aboveskip=0.1pt} 
	\subfloat[\hspace{-1em}(a) Modal 1 of SMC-GM]{
		\includegraphics[width=0.3\linewidth]{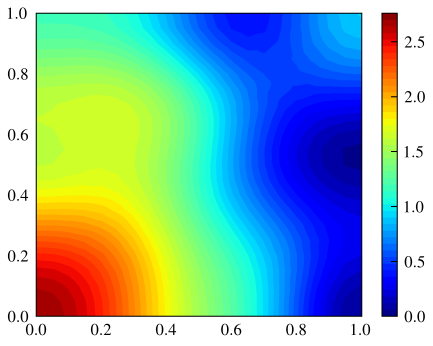}
	}
	\subfloat[\hspace{-1em}(b) Modal 2 of SMC-GM]{
		\includegraphics[width=0.3\linewidth]{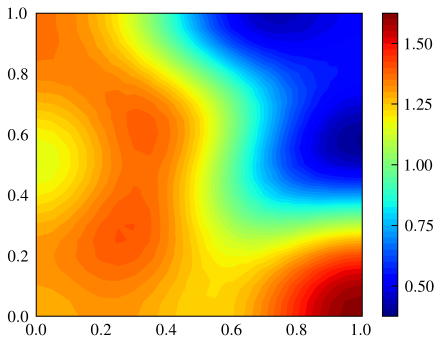}
	}
	\subfloat[\hspace{-1em}(c) Modal 3 of SMC-GM]{
		\includegraphics[width=0.3\linewidth]{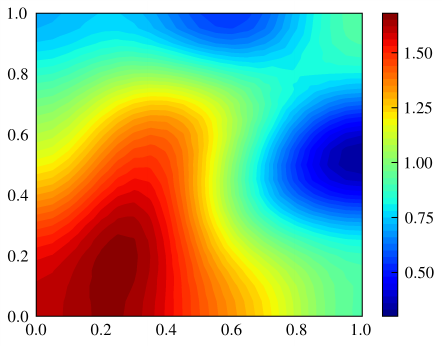}
	}
	
	\subfloat[\hspace{-1em}(d) Modal 4 of SMC-GM]{
		\includegraphics[width=0.3\linewidth]{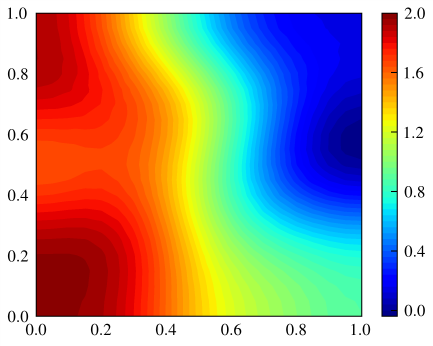}
	}
	\subfloat[\hspace{-1em}(e) Modal 1 of SMC-pCN]{
		\includegraphics[width=0.3\linewidth]{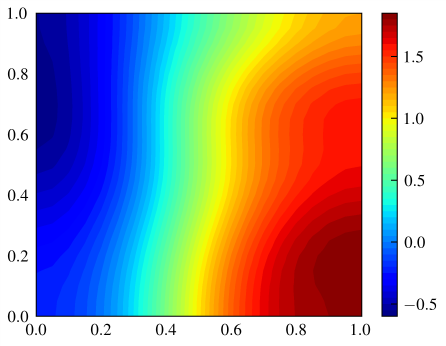}
	}
	\subfloat[\hspace{-1em}(f) Modal 2 of SMC-pCN]{
		\includegraphics[width=0.3\linewidth]{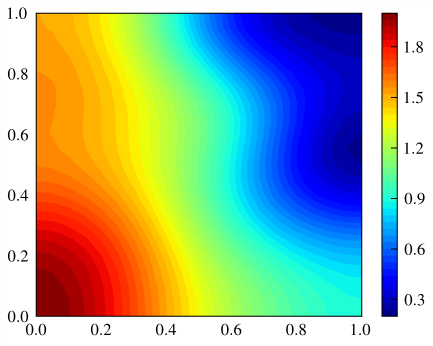}
	}
	\vskip 0.3 cm
	\caption{\emph{\small 
			(a)-(d) display four out of the eight cluster means from the SMC-GM sampling results, while (e) and (f) show two out of the five cluster means from the SMC-pCN sampling results.
	}}	
	\label{fig:Example3Modals}
\end{figure}

\begin{figure}[ht]
	\centering
	\captionsetup[subfigure]{skip=3pt} 
	\captionsetup[subfigure]{labelformat=empty}
	\captionsetup{aboveskip=0.1pt} 
	\subfloat[\hspace{3.2em}(a) Error of SMC-GM]{
		\includegraphics[width=0.4\linewidth]{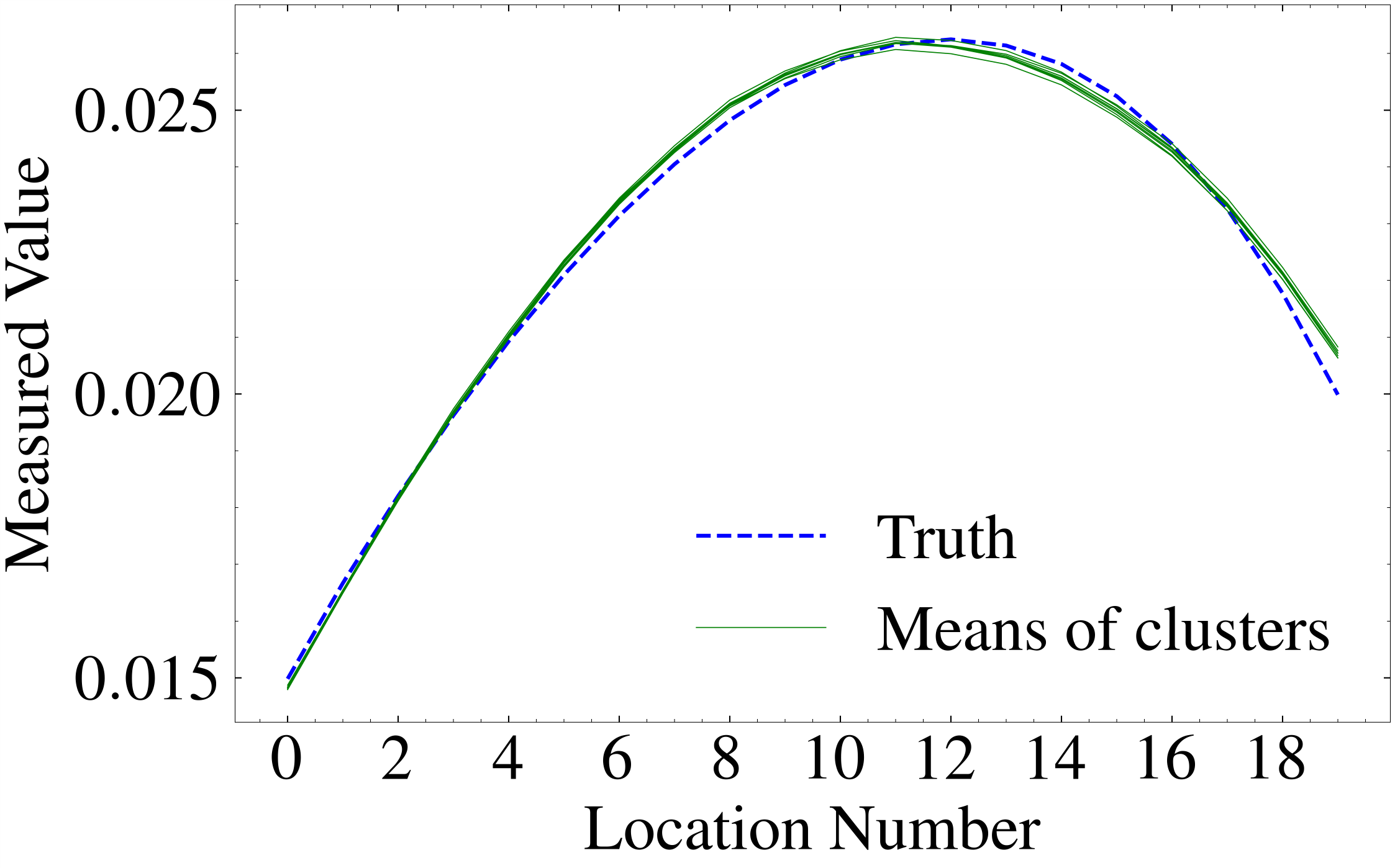}
	}\hspace{2em}
	\subfloat[\hspace{3.2em}(b) Error of SMC-pCN]{
		\includegraphics[width=0.4\linewidth]{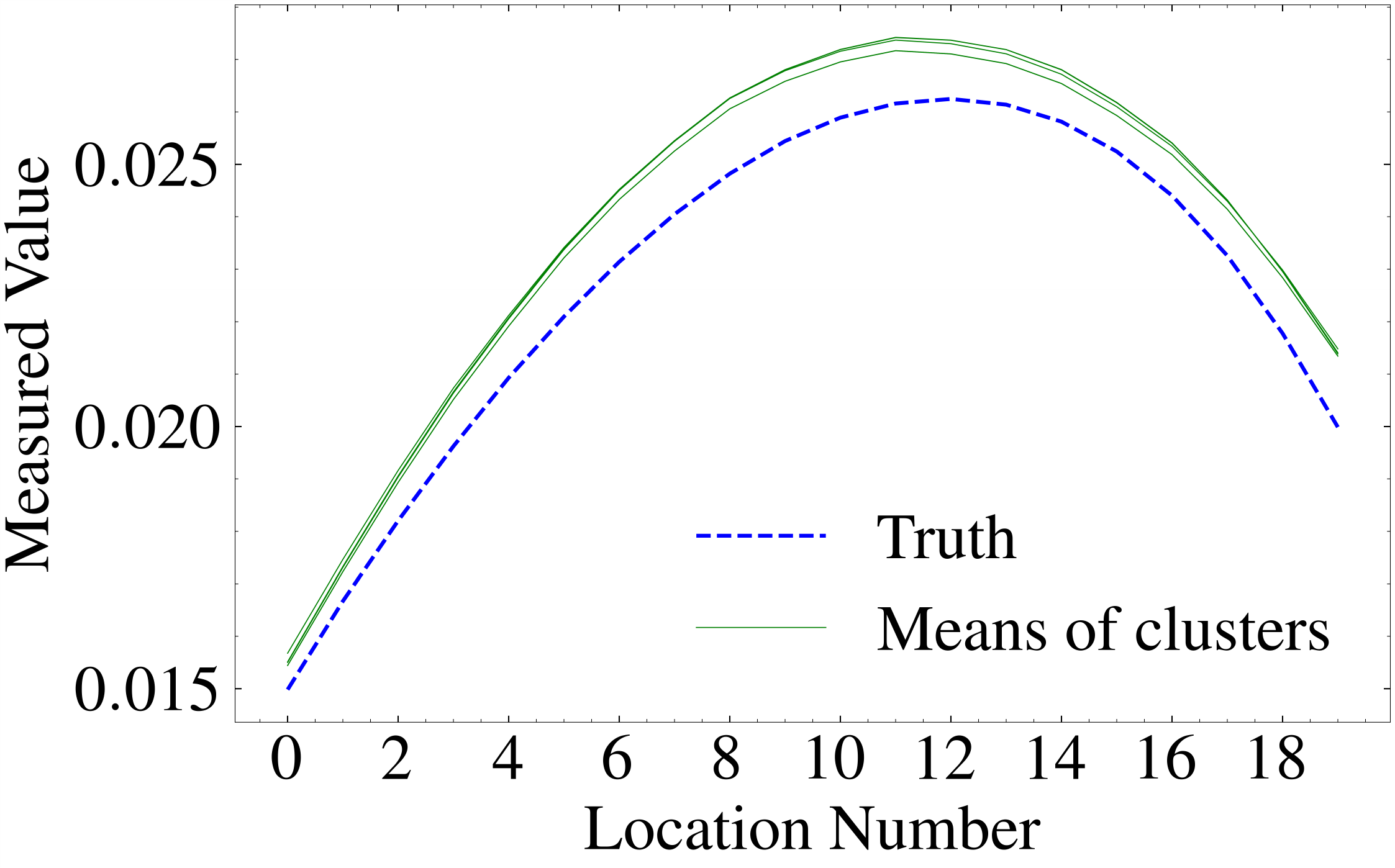}
	}
	\vskip 0.3 cm
	\caption{\emph{\small For multi-modal posteriors, the sample mean may not capture sufficient information. Instead, we solve the forward problem using cluster means of samples as parameters, and compare the solutions to the actual data. (a): The results of SMC-GM. (b): The results of SMC-pCN.
	}}	
	\label{fig:Example3Err}
\end{figure}

\paragraph{Numerical results} 
In the numerical experiment, we will compare the computational efficiency and accuracy of the two methods. Our algorithm demonstrates significant improvements in sampling efficiency. The time cost of SMC-GM is 92\,s, and that of pCN is 1441\,s. Since both have the same number of samples, we can conclude that the sampling efficiency of SMC-GM is 15.66 times that of SMC-pCN under our experimental setup.

Regarding accuracy, we assess whether SMC-GM accurately estimates the posterior measure of the unknown parameter \(u\). 

Note that the posterior measure \(\mu^\mathrm{d}\) is defined on a function space, and the covariance, being an operator that maps between function spaces, cannot be visualized directly. We model \(u\) as a Fourier series
\begin{align}
u(x) = \overline{u}(x) + \sum_{k=1}^\infty u_k e_k,\label{Fourier series}
\end{align}
where \(\overline{u}(x)\) is the mean function, \(e_k\) is the \(k\)-th basis function, and \(u_k\) is the \(k\)-th Fourier coefficient. Due to Theorem VI.16 in \cite{Reed1980}, the eigenfunctions of the prior covariance operator form a complete orthonormal basis. Thus, we can choose them as \(e_k\).

Firstly, in Figure \ref{fig:example1_mean}(b) and \ref{fig:example1_mean}(c), we present the posterior mean estimates of \(u\) derived from the SMC-GM and the SMC-pCN methods, respectively, which appear visually similar. We also include the background true function in Figure \ref{fig:example1_mean}(a) for comparison. The relative errors of the two methods are
\begin{align*}
err_{\text{\tiny GM}} = \frac{\|\overline{u}_{\text{\tiny GM}} - u_{\text{true}}\|_{L^2}}{\|u_{\text{true}}\|_{L^2}} = 0.0271, \quad
err_{\text{\tiny pCN}} = \frac{\|\overline{u}_{\text{\tiny pCN}} - u_{\text{true}}\|_{L^2}}{\|u_{\text{true}}\|_{L^2}} = 0.0254,
\end{align*}
where the \(L^2\) norm is defined as \(\|u\|_{L^2} := \sqrt{\int_{\Omega} u(x)^2 \, dx}\). Thus, the estimated posterior mean function of the parameter \(u\) obtained by SMC-GM is quantitatively similar to that obtained by SMC-pCN.

Secondly, we provide some discussions of the estimated posterior covariance functions. Using the samples from the posterior, we can estimate the posterior marginal densities of \(u_k\) defined in \eqref{Fourier series}. We uniformly selected 6 coefficients from the first 20 and plotted their density estimates in Figure \ref{fig:DensityOfCoef}. In each sub-figure, we use a red line to represent the results obtained by SMC-GM method and a blue line for the results obtained by SMC-pCN method. Finally, we can compute the total variation distance between two probability density functions
\begin{align*}d_{\text{\tiny TV}}(p_1, p_2) = \frac{1}{2} \int_{\mathbb{R}} |p_1(x) - p_2(x)| \, dx\end{align*}
to see if we have a good estimation, where the constant \(\frac{1}{2}\) appears in the definition to ensure that \(0 \leq d_{\text{\tiny TV}}(p_1, p_2) \leq 1\). The average total variation distance of the first 20 coefficients is
\begin{align}
\frac{1}{20} \sum_{i=1}^{20} d_{\text{\tiny TV}}(p^{SMC-GM}_i, p^{SMC-pCN}_i) = 0.15.
\end{align}
Thus, the error of the marginal density estimation by SMC-GM is small in terms of the average total variation distance.

\paragraph{Mesh independence}
Finally, we show the mesh independence of the SMC-GM method, a key property for a well-defined function space method. At each layer, the particle system determines the temperature of the next layer through the procedure described in \ref{SM-NumericalDetails}. As a result, differences between particle systems arising from different discretizations are reflected through the corresponding temperatures. 
 Figure \ref{fig_MeshIdpdt}(a) shows temperature plots for discrete dimensions \(N = \{400, 1600, 3600, 6400, 10000\}\), starting at 0 and increasing to 1. We observe that the temperatures show mesh independence.

For comparison, we also plot the temperature curves of the SMC-RW method, which lacks a well-defined definition in function space and is thus mesh-dependent. The curves in Figure \ref{fig_MeshIdpdt}(b) are distinct from one another. Furthermore, as the discrete dimension increases, the length of the curves increases as well, which will result in an infinite sampling time.

\subsection{Sparse measurement situation}
Then, we consider the same forward problem as at the beginning of this section, but with sparse measurement points, i.e., the number of observation points \(N_d\) is too small to determine a unique solution of the inverse problem. In this situation, the posterior measure is often multi-modal. Let \(N_d = 20\), and the measurement points \(\{x_i\}_{i=1}^{N_d}\) are taken at the coordinates \(\{(0.8, 0.2 + 0.03i)\}_{i=1}^{20}\). Denote the number of modals after clustering as \(N_{modal}\), then we can define an average \(l^2\) error
\begin{align*}
err:=\frac{1}{N_{modal}}\sum_{i=1}^{N_{modal}} \|\mathcal{F}u_i-d\|_{l^2},
\end{align*}
where the \(\| \cdot \|_{l^2}\) denote the usual \(l^2\) norm \(\|x\|_{l^2}=\sqrt{x^Tx}\).

\begin{table}[ht]
	\centering
	\renewcommand{\arraystretch}{1.2}
	\begin{tabular}{c|cccc}
		\Xhline{0.9pt}
		Method & Time cost & Particles number & Number of modals & Average error\\
		\hline
		SMC-GM & 5500.5 s& 60000 & 8& 0.0049\\
		\hline
		SMC-pCN & 8845.0 s& 6000 & 4 & 0.022\\
		\Xhline{0.9pt}
	\end{tabular}
	\vskip 0.1 cm
	\caption{Comparison of the performance of pCN and GM as transition kernel in SMC.}
	\label{tab:DarcySparseMeasure}
\end{table}

As shown in Table \ref{tab:DarcySparseMeasure}, we see that obtaining ten times the number of samples with the SMC-GM algorithm requires less time than with SMC-pCN. In other words, the efficiency of the former is thirteen times that of the latter. Considering the ability of the two algorithms to explore complex posterior distributions, we observe that SMC-GM detected 8 peaks, while SMC-pCN detected 4 peaks. We have depicted some of the cluster means in Figure \ref{fig:Example3Modals} (a)-(f). We also compare the data to the solutions of the forward problem using these means as parameters, as illustrated in Figure \ref{fig:Example3Err}. They all fit well, which also illustrates the multi-modal nature of the posterior distribution. In summary, compared to the SMC-pCN method, our proposed SMC-GM algorithm can obtain more samples in less time, detect more posterior peaks, and better align with the data when applied to the forward problem.
\section{Discussion}
In this paper, we propose a novel transition kernel, pCN-GM, which is derived from the Crank-Nicolson discretization of the Langevin system and is proven to be well-defined in infinite-dimensional spaces. By employing pCN-GM in the mutation step of SMC, we obtain the SMC-pCN-GM algorithm. As a theoretical foundation for SMC, we prove a convergence theorem under weaker conditions. Given the proven density of the Gaussian mixture measure, we propose an efficient approximation algorithm called SMC-GM and apply it to three typical inverse problems. All numerical results indicate that SMC-GM excels over SMC-pCN in terms of efficiency and multi-modal sampling capability, with comparable accuracy, and also exhibits mesh-independence. 

We have not provided an explicit expression for the error introduced by the Gaussian mixture approximation in the infinite-dimensional setting, which has been studied in the finite-dimensional context \cite{Zhang2020}. 
Further acceleration of SMC is an interesting direction that we plan to explore in future work. 
One possible approach is to introduce parallel MCMC schemes into the mutation step of SMC, for example via multi-level Monte Carlo methods~\cite{Wang2023}. 
Another possibility is to parallelize SMC itself~\cite{Verge2015}, by running two or more SMC samplers simultaneously and allowing state exchange between them.
\Acknowledgements{
We thank the reviewers for their valuable comments, which have greatly improved the quality of this paper. This work was partially supported by the National Natural Science Foundation of China (Grants No. 12322116, 12271428, 12326606, 12090021 and 12090020), the Fundamental and Interdisciplinary Disciplines Breakthrough Plan of the Ministry of Education of China (Grant No. JYB2025XDXM101), and the Tianyuan Fund for Mathematics of the National Natural Science Foundation of China (Grant No. 12426105).
}




\newpage
\begin{appendix}

\section{Four transition kernels and a note of  irreducibility}\label{SM-4Kernels}
In this section, we introduce four methods for the mutation step in SMC, each with a corresponding transition kernel \( Q(u, dv) \) and acceptance rate \( a(u, v) \). Here, \( u \) denotes the current state, and \( v \) denotes the next state. The general framework is the Metropolis-Hastings method, a renowned Markov chain Monte Carlo approach, which is described in Algorithm \ref{Alg:MG}.
\begin{algorithm}
	\caption{Metropolis-Hastings method}
	\label{Alg:MG}
	\begin{algorithmic}[]
		\STATE{\textbf{Input}: Specify the initial state \(u_0\) and the length of the Markov chain \(I\). }
		\STATE{\textbf{For} \(i=1,\cdots ,I\), \textbf{do}:}
		\STATE{\qquad Draw \(v_{i}\) from a transition kernel \(Q(u_{i-1},dv)\);}
		\STATE{\qquad Let \(u_{i}=v_{i}\) with probability \(a(u_i,v_i)\);}
		\STATE{\textbf{Output}: \(\{u_i\}_{i=1}^I\).}
	\end{algorithmic}
\end{algorithm}
Note that the Metropolis-Hastings method can be specified by \(Q(u,dv)\) and \(a(u,v)\). The acceptance rate function is always calculated by
\begin{align}
a(u,v)=\min\left\{1,\frac{\mu^\text{d}(du)Q(u,dv)}{\mu^\text{d}(dv)Q(v,du)}\right\}.\label{Apdx-auv}
\end{align}
Thus, we only provide examples of proposal distributions:
\begin{itemize}
	\item RW(Random walk): \begin{align*}Q_{\text{\tiny RW}}(u,dv)=u+\beta \mathcal{N}_{0,\mathcal{C}}(dv),\end{align*}
	\item pCN(preconditioned Crank-Nicolson): \begin{align*}Q_{\text{\tiny pCN}}(u,dv)=\sqrt{1-\beta^2}u+\beta  \mathcal{N}_{0,\mathcal{C}}(dv),\end{align*}
	\item pCN-GM(pCN based on Gaussian mixture approximation): \begin{align*}Q_{\text{\tiny pCN-GM}}(u,dv)=\sum_{i=1}^Mw_i m_i+\sqrt{1-\beta^2}\sum_{i=1}^Mw_i(u_i-m_i)+\beta\sum_{i=1}^Mw_i \mathcal{N}_{0,\mathcal{C}_i}(dv).\end{align*}
\end{itemize}
Here, $\beta$ represents the step size of the transition, $\mathcal{C}$ is the covariance operator of the prior, and the parameters of the Gaussian mixture measure, $\{w_i, m_i, \mathcal{C}_i\}_{i=1}^M$, should be determined to approximate the posterior with $\sum_{i=1}^M w_i \mathcal{N}(m_i, \mathcal{C}_i)$. The final proposal using in the mutation is given by
\begin{align*}Q_{\text{\tiny GM}}(u,dv)=\sum_{i=1}^Mw_i\mathcal{N}_{m_i,\mathcal{C}_i}(dv),\ a(u,v)=1,\end{align*}
which is an approximation of \(Q_{\text{\tiny pCN-GM}}\). The four transition kernels are closely related; \(Q_{\text{\tiny pCN}}\) is a generalization of \(Q_{\text{\tiny RW}}\) to the infinite-dimensional setting, \(Q_{\text{\tiny pCN-GM}}\) extends \(Q_{\text{\tiny pCN}}\) for multimodal distributions, and \(Q_{\text{\tiny GM}}\) is a special case of \(Q_{\text{\tiny pCN-GM}}\) (\(\beta=1\)). Note that the acceptance rate function defined by \eqref{Apdx-auv} may not be well-defined in the infinite-dimensional setting. For instance, the acceptance rate function derived from \(Q_{\text{\tiny RW}}\) is ill-defined, as discussed in Section \ref{SM-UniquenessOfpCNGM}, which outlines a necessary condition for the well-definedness of the acceptance rate function.

We then discuss whether the pCN algorithm is irreducible in the infinite-dimensional Hilbert space setting, and its implications for the SMC framework.
In the infinite-dimensional setting, the proposal of pCN is
\[P(u,dv)=\sqrt{1-\beta^2}u+\beta^2\xi,\quad \xi\sim \mu_{prior}=\mathcal{N}(0,\mathcal{C}), u,v\in \mathcal{H},\]
and $\mathcal{C}$ is a positive, symmetric, trace-class operator. In infinite-dimensional spaces, the Lebesgue measure does not exist, and two Gaussian measures with different means are mutually singular if the difference of their means lies outside the Cameron–Martin space \(\tilde{\mathcal{H}}\).
 The singularity of $\mu$ and $\nu$ means that there exists a measurable set $A$ such that $\mu(A)=1$ and $\nu(A)=0$, which is denoted by \(\mu\perp \nu\). Thus we have
\[\mathcal{N}(0,\mathcal{C})\perp \mathcal{N}(\sqrt{1-\beta^2}x,\beta^2\mathcal{C}),\quad\forall x\notin \tilde{\mathcal{H}} .\]
As a consequence of this singularity, for sets \(A\) with positive posterior measure \(\mu_{\text{post}}(A) > 0\), it generally holds that \(P(u,A)=0\) for \(\mu_{\text{post}}\)-almost all \(u\) \cite{Hairer2014}.
Therefore, under the standard finite-dimensional definition of irreducibility, the pCN transition kernel on infinite-dimensional state spaces fails to be irreducible. This phenomenon is not due to the properties of the pCN algorithm, but rather to the intrinsic nature of the infinite-dimensional setting.

Consequently, alternative approaches are required to establish ergodicity and related results. It is well known that for Markov chains with reversible transition kernels, a strictly positive spectral gap is equivalent to a form of geometric ergodicity~\cite{Roberts1997}. The uniqueness of the invariant distribution of pCN, together with the existence of an \(L^2\)-spectral gap and a Wasserstein spectral gap, is established in~\cite{Hairer2014}. Moreover, this spectral gap is dimension-independent, and laws of large numbers as well as central limit theorems are derived. Therefore, using pCN to sample in infinite-dimensional function spaces is feasible.

Finally, in SMC methods, the primary sampling mechanism is importance sampling, while MCMC steps are introduced to enhance particle diversity. Neither the theoretical analysis nor practical implementations of SMC require the Markov transitions to be run for long time horizons, and irreducibility of the underlying Markov kernels is therefore not a central requirement within the SMC framework. 

\section{From Langevin system to pCN-GM proposal}\label{SM-Langevin2pCNGM}
In this section, we detail the process of deriving the pCN-GM proposal from a system of Langevin equations. Consider the Langevin equation \cite{Andrew_review}
\begin{align*}
\frac{du}{ds} = -\mathcal{K}(\mathcal{L}u + \gamma D\Phi(u)) + \sqrt{2\mathcal{K}}\frac{db}{ds},
\end{align*}
where \( u \) represents velocity, \( b \) is a Brownian motion in \(\mathcal{H}\) with the covariance operator being the identity, \(\mathcal{L} = \mathcal{C}^{-1}\), and \(\mathcal{K} = \mathcal{C}\) acts as a preconditioner. Set \(\gamma = 0\) and consider a dynamical system of \( M \) particles with distinct distribution perturbations:
\begin{equation*}
\left\{
\begin{aligned}
\frac{du_1}{ds} &= -u_1 + \sqrt{2\mathcal{C}_1}\frac{db}{ds}, \\
\frac{du_2}{ds} &= -u_2 + \sqrt{2\mathcal{C}_2}\frac{db}{ds}, \\
&\vdots \\
\frac{du_M}{ds} &= -u_M + \sqrt{2\mathcal{C}_M}\frac{db}{ds}.
\end{aligned}
\right.
\end{equation*}
Let \( w_j \) denote mass, and we will examine the evolution of the average velocity
\begin{align*}
u = \frac{\sum_{j=1}^M w_j u_j}{\sum_{j=1}^M w_j}.
\end{align*}
By multiplying each equation by its respective mass and summing, we obtain
\begin{align*}
\sum_{j=1}^M w_j \frac{du_j}{ds} = -\sum_{j=1}^M w_j u_j + \sum_{j=1}^M w_j \sqrt{2\mathcal{C}_j} \frac{db}{ds}.
\end{align*}
Discretize the equation using the Crank-Nicolson method, and it follows that
\begin{align*}
\sum_{j=1}^M w_j \frac{v_j - u_j}{\delta} =& -\sum_{j=1}^M w_j \frac{u_j + v_j}{2} + \sum_{j=1}^M w_j \sqrt{2\mathcal{C}_j} \frac{b(s+\delta) - b(s)}{\delta}, \\
\sum_{j=1}^M w_j (v_j - u_j) =& -\delta \sum_{j=1}^M w_j \frac{u_j + v_j}{2} + \sum_{j=1}^M w_j \sqrt{2\mathcal{C}_j} \xi,\ \xi\sim \mathcal{N}(0,\delta I).
\end{align*}
Rearrange the terms to obtain the following equation:
\begin{align*}
\sum_{j=1}^M w_j \left(1 + \frac{1}{2}\delta \right) v_j =& \sum_{j=1}^M w_j \left(1 - \frac{1}{2}\delta\right) u_j + \sum_{j=1}^M w_j \sqrt{2\delta \mathcal{C}_j} \xi,\ \xi\sim \mathcal{N}(0,\delta I), \\
\left(2 + \delta \right) \sum_{j=1}^M w_j v_j =& \left(2 - \delta\right) \sum_{j=1}^M w_j u_j + \sqrt{8\delta} \sum_{j=1}^M w_j \xi_j,\ \xi_j\sim \mathcal{N}(0,\mathcal{C}_j), \\
\sum_{j=1}^M w_j v_j =& \frac{2 - \delta}{2 + \delta} \sum_{j=1}^M w_j u_j + \frac{\sqrt{8\delta }}{2 + \delta} \sum_{j=1}^M w_j \xi_j,\ \xi_j\sim \mathcal{N}(0,\mathcal{C}_j).
\end{align*}
Let \( \beta := \frac{\sqrt{8\delta }}{2 + \delta} \), then
\begin{align}
\sum_{j=1}^M w_j v_j = \sqrt{1 - \beta^2} \sum_{j=1}^M w_j u_j + \beta \sum_{j=1}^M w_j \xi_j,\ \xi_j\sim \mathcal{N}(0,\mathcal{C}_j).\label{Apdx pCN-GM0}
\end{align}

Divide by the total mass, but continue to denote it as \( w_j \), ensuring that \( \sum_{j=1}^M w_j = 1 \). Now, \(w_j\) has the same meaning as \(w_j\) in the main text. This approach can yield a mixed version of the pCN, i.e.,
\begin{align*}
v = \sqrt{1 - \beta^2} u + \beta \sum_{j=1}^M w_j \xi_j,\ \xi_j\sim \mathcal{N}(0,\mathcal{C}_j).
\end{align*}
In other words, we obtain a pCN method with mean-free Gaussian mixture measure. Subtract the mean from \eqref{Apdx pCN-GM0}, and we obtain that
\begin{align*}
\sum_{j=1}^M w_j (v_j - m_j) = \sqrt{1 - \beta^2} \sum_{j=1}^M w_j (u_j - m_j) + \beta \sum_{j=1}^M w_j \xi_j,\ \xi_j\sim \mathcal{N}(0,\mathcal{C}_j).
\end{align*}
The evolution of the average velocity is
\begin{align*}
v = \sqrt{1 - \beta^2} u + (1 - \sqrt{1 - \beta^2}) \sum_{j=1}^M w_j m_j + \beta \sum_{j=1}^M w_j \xi_j,\ \xi_j\sim \mathcal{N}(0,\mathcal{C}_j),
\end{align*}
which is precisely the pCN-GM transition kernel, corresponding to a Gaussian mixture measure \( \sum_{j=1}^M w_j \mathcal{N}(m_j, C_j) \).

\section{Some numerical details}\label{SM-NumericalDetails}
\textbf{Determine the temperatures}:  
We need to choose the temperatures \(h_j\in(0,1),1\leq i\leq J\) such that \(\sum_{j=1}^Jh_j=1\). If these temperatures are too small, the SMC will run slowly; if they are too large, the adjacent posteriors \(\mu_j\) and \(\mu_{j+1}\) will differ significantly, which is not conducive to importance sampling \cite{Dai2022}. We can use the effective sample size (ESS) to select appropriate temperatures:
\begin{align*}
\text{ESS} := \left(\sum_{i=1}^N w_i^2\right)^{-1},
\end{align*}
where \(w_i\) is the weight of the \(i\)-th particle. A small \(h_j\) can result in a large ESS. An excessively large value of \(h_j\) can cause the loss of a large number of samples during the resampling step, thereby leading to a small ESS. Thus, we can use a simple bisection method to find an \(h_j\) such that \(\text{ESS} > N_{\text{thresh}} := 0.6 N\) \cite{Beskos2015SC}.

\textbf{Settings of pCN transition kernel}: 
We should determine the step size parameter \( \beta \) in the proposal (\ref{pCNproposal}) of the main text to target average acceptance rates in a neighborhood of 0.2 \cite{Gelman1997,Beskos2009}. Let the average acceptance rate be denoted as \( \alpha^N_j \), where \( N \) represents the number of particles, and \( j \) is the index of the layer within the SMC algorithm. Then, we can employ an adaptive strategy as follows \cite{Beskos2015SC}:
\begin{align*}
\beta_{j+1} = \begin{cases}
2\beta_j, & \text{if } \alpha^N_j > 0.3, \\
0.5\beta_j, & \text{if } \alpha^N_j < 0.15, \\
\beta_j, & \text{if } 0.15 \leq \alpha^N_j \leq 0.3.
\end{cases}
\end{align*}
This update formula adjusts the step size of the next layer based on the average acceptance rate of each layer, thereby maintaining an average acceptance rate as close as possible to 0.2.

\section{Details of the main text}\label{SM-AllProofs}
In this section, we offer detailed proofs and some discussions for all the theorems and lemmas that were mentioned in the main text. These proofs are essential for a rigorous understanding of the mathematical framework and the theoretical underpinnings of the algorithms discussed.
\subsection{Discussion of total variation distance}\label{SM-TV}
In Section \ref{Sec.SMC-GM} we use the following distance:
\begin{align}
d(\mu,\nu)=\sup _{|f|_\infty\leq 1}\sqrt{\mathbb{E}_\omega \left|\int fd\mu-\int fd\nu\right|^2}.\label{Apdx-rTV}
\end{align}
Now we prove that \(d(\mu,\nu)\) is a random version of total variation distance defined by
\begin{align*}d_{\text{\tiny TV}}(\mu,\nu)=\int \left|\frac{d\mu}{d\eta}-\frac{d\nu}{d\eta}\right|d\eta,\end{align*}
where \(\eta\) is a measure satisfying \(\mu\ll \eta\) and \(\nu\ll \eta.\) 

If \(\mu\) and \(\nu\) are not random, it follows from formula \eqref{Apdx-rTV} that
\begin{align*}
d(\mu,\nu)=\sup _{|f|_\infty\leq 1}\left|\int fd\mu-\int fd\nu\right|.
\end{align*}
Let \(\eta = \frac{\mu + \nu}{2}\). Then, \(\mu \ll \eta\) and \(\nu \ll \eta\). Thus
\begin{align*}
d(\mu,\nu)=\sup _{|f|_\infty\leq 1}\left|
\int f\frac{d\mu}{d\eta}d\eta-\int f\frac{d\nu}{d\eta}d\eta
\right|=\sup _{|f|_\infty\leq 1}\left|
\int f\left(\frac{d\mu}{d\eta}-\frac{d\nu}{d\eta}\right)d\eta
\right|.
\end{align*}
It is obvious that
\begin{align*}
d(\mu,\nu)\leq \sup _{|f|_\infty\leq 1}
\int \left|f\left(\frac{d\mu}{d\eta}-\frac{d\nu}{d\eta}\right) \right|d\eta
\leq
\int \left|\frac{d\mu}{d\eta}-\frac{d\nu}{d\eta} \right|d\eta
=d_{{\text{\tiny TV}}}(\mu,\nu).
\end{align*}
On the other hand, let \(f_0=\text{sgn}(\frac{d\mu}{d\eta}-\frac{d\nu}{d\eta})\), where \(\text{sgn}(x)\) is the sign function. Now we have \(|f_0|_\infty\leq 1\) and
\begin{align*}
\left|
\int f_0\left(\frac{d\mu}{d\eta}-\frac{d\nu}{d\eta}\right)d\eta
\right|=\int \left|\frac{d\mu}{d\eta}-\frac{d\nu}{d\eta} \right|d\eta
=d_{{\text{\tiny TV}}}(\mu,\nu).
\end{align*}
Thus, \(d(\mu,\nu) = d_{\text{\tiny TV}}(\mu,\nu)\) holds for non-random measures \(\mu\) and \(\nu\), and we can regard \(d(\mu,\nu)\) as a random version of \(d_{\text{\tiny TV}}(\mu,\nu)\).
\subsection{Proof of Lemma \ref{Lemma1-3}(c) (error of operator \(L\))}\label{SM-ProofLemmaC}
Under Assumptions \ref{A1}, we have
\begin{align*}
d(L_j\mu_{j-1},L_j\mu)\leq 2\kappa_1^{-1}Z_j^{-1}d(\mu_{j-1},\mu).
\end{align*}
\begin{proof}
	Let $ g_j(v):=e^{-\Phi_j(u)} $, then it follows from Assumptions \ref{A1}(a) that $ |\kappa_1 g_j(v)|\leq 1$. For simplicity we denote $ g'=\kappa_1g_j $ with $ |g'|\leq 1$. 
	The aim is to measure the distance between $ L_j\mu_{j-1} $ and $ L_j \mu $. Hence we need to compute the difference between the integral of \(f\) with respect to two measures, i.e.,
	\begin{align*}
	(L_j\mu_{j-1})(f)-(L_j\mu)(f)&=
	\frac{\mu_{j-1}(fg_j)}{\mu_{j-1}(g_j)}-\frac{\mu(fg_j)}{\mu(g_j)}\\
	&=\frac{\mu_{j-1}(fg_j)}{\mu_{j-1}(g_j)}-\frac{\mu(fg_j)}{\mu_{j-1}(g_j)}
	+\frac{\mu(fg_j)}{\mu_{j-1}(g_j)}-\frac{\mu(fg_j)}{\mu(g_j)}\\
	&=\frac{\kappa_1^{-1}}{\mu_{j-1}(g_j)}[\mu_{j-1}( fg')-\mu( fg')]
	+\frac{\mu(fg_j)}{\mu(g_j)}\frac{\kappa_1^{-1}}{\mu_{j-1}(g_j)}[\mu( g')-\mu_{j-1}( g')].
	\end{align*}
	Take square of both sides, and using $ (a+b)^2\leq 2(a^2+b^2) $ we obtain
	\begin{align*}
	|(L_j\mu_{j-1})(f)-(L_j\mu)(f)|^2
	\leq&\frac{2\kappa_1^{-2}}{\mu_{j-1}(g_j)^2}|\mu_{j-1}( fg')-\mu( fg')|^2
	\\
	&+2\frac{\mu(fg_j)^2}{\mu(g_j)^2}\frac{\kappa_1^{-2}}{\mu_{j-1}(g_j)^2}|\mu( g')-\mu_{j-1}( g')|^2.
	\end{align*}
	By applying Assumptions \ref{A1}, we have
	\begin{align*}	|(L_j\mu_{j-1})(f)-(L_j\mu)(f)|^2
	\leq2\kappa_1^{-2}Z_j^{-2}|\mu_{j-1}( fg')-\mu(fg')|^2
	+2|f|_\infty^2\kappa_1^{-2}Z_j^{-2}|\mu(g')-\mu_{j-1}(g')|^2.\end{align*}
	Take the supremum of both sides for $|f|_\infty \leq 1$. Then, it follows that
	\begin{align*}
	&d(L_j\mu_{j-1},L_j\mu)\\&=\sup_{|f|_\infty\leq1}\sqrt{\mathbb{E}_\omega|(L_j\mu_{j-1})(f)-(L_j\mu)(f)|^2}\\
	&\leq \sup_{|f|_\infty\leq1}\sqrt{2\kappa_1^{-2}Z_j^{-2}\mathbb{E}_\omega|\mu_{j-1}( fg')-\mu(fg')|^2
		+2|f|_\infty^2\kappa_1^{-2}Z_j^{-2}\mathbb{E}_\omega|\mu(g')-\mu_{j-1}(g')|^2}\\
	&\leq \sqrt{2\kappa_1^{-2}Z_j^{-2}\sup_{|f|_\infty\leq1}\mathbb{E}_\omega|\mu_{j-1}( fg')-\mu(fg')|^2
		+2\kappa_1^{-2}Z_j^{-2}d(\mu,\mu_{j-1})^2}\\
	&= \sqrt{2\kappa_1^{-2}Z_j^{-2}d(\mu,\mu_{j-1})^2
		+2\kappa_1^{-2}Z_j^{-2}d(\mu,\mu_{j-1})^2}\\
	&= \frac{2}{\kappa_1Z_j}d(\mu,\mu_{j-1}),
	\end{align*}
	which completes the proof.
\end{proof}

\subsection{Proof and discussion of Theorem \ref{SMC approximation theorem} (convergence of SMC)}\label{SM-ConvergenceSMC}
\begin{proof}
	First we give an iterative error estimate: 
	\begin{align*}
d(\mu_{j+1}^N,\mu_{j+1})&=d(S^NP_{j+1}L_{j+1}\mu_j^N,\mu_{j+1})\\
	(\text{triangle inequality},\text{Lemma }\ref{Lemma1-3}(a))\quad
	&\leq \frac{1}{\sqrt{N}}+d(P_{j+1}L_{j+1}\mu_j^N,P_{j+1}\mu_{j+1})\\
	(\text{Lemma }\ref{Lemma1-3}(b))\quad
	&\leq\frac{1}{\sqrt{N}}+ d(L_{j+1}\mu_j^N,L_{j+1}\mu_{j}).\\
    (\text{Lemma }\ref{Lemma1-3}(c))\quad
	&\leq \frac{1}{\sqrt{N}}+2\kappa_1^{-1}Z_{j+1}^{-1}d(\mu_j^N,\mu_j)\\
	\end{align*}
	Notice that \(d(\mu_0^N,\mu_0)\leq \frac{1}{\sqrt{N}}\). 
	Iterating completes the proof, and the result is given by
\begin{align*}
d(\mu_J^N, \mu_J) \le \frac{1}{\sqrt{N}} \left(
1 + \sum_{j=1}^{J} \frac{2^j}{\kappa_1^j} \prod_{k=1}^{j}\frac{1}{Z_{J+1-k}}
\right).
\end{align*} 
\end{proof}
\subsection{Proof of Theorem \ref{mu_0dvPvdu is well defined} (the well-definedness of \(\mu_0(dv)Q(v,du)\))}
\label{SM-mu_0dvPvdu is well defined}
\begin{proof}
	The measure \(\mu(du)Q(u,dv)\) is defined on  \(\mathcal{H}\times\mathcal{H}\). Consider its characteristic function
	\begin{align*}
	&\int_{\mathcal{H}\times \mathcal{H}}e^{i(u,\xi)}e^{i(v,\eta)}\mu(du)Q(u,dv)\\
	=&\int_{\mathcal{H}} e^{i(u,\xi)} \mu(du)\int_{\mathcal{H}}e^{i(v,\eta)}Q(u,dv)\\
	=&\int_{\mathcal{H}} e^{i(u,\xi)} \mu(du)\sum_{j=1}^Mw_je^{i(\gamma u+(1-\gamma) m_j,\eta)-\frac{1}{2}(\eta,\beta^2\mathcal{C}_j\eta)}\\
	=&\sum_{j=1}^Mw_je^{i((1-\gamma) m_j,\eta)-\frac{1}{2}(\eta,\beta^2\mathcal{C}_j\eta)} \int_{\mathcal{H}} e^{i(u,\xi+\gamma \eta)} \mu(du)\\
	=&\sum_{j=1}^Mw_je^{i((1-\gamma) m_j,\eta)-\frac{1}{2}(\eta,\beta^2\mathcal{C}_j\eta)}
	e^{-\frac{1}{2}(\xi+\gamma \eta,\mathcal{C}(\xi+\gamma \eta))}\\
	=&\sum_{j=1}^Mw_j
	e^{i((1-\gamma) m_j,\eta)
		-\frac{\beta^2}{2}(\eta,\mathcal{C}_j\eta)
		-\frac{1}{2}(\xi,\mathcal{C}\xi)
		-\frac{\gamma^2}{2}(\eta,\mathcal{C}\eta)
		-\gamma(\eta,\mathcal{C}\xi)},
	\end{align*}
	which is the characteristic function of an \(M\)-component Gaussian mixture measure in \(\mathcal{H}\times\mathcal{H}\). The mean of \(j\)-th component is \((0,(1-\gamma) m_j)\), and the covariance operator is
	\begin{align*}\mathcal{V}_{j}=\left[\begin{array}{cc}
	\mathcal{C} & \gamma \mathcal{C} \\
	\gamma \mathcal{C} & \beta^2\mathcal{C}_j+\gamma^2\mathcal{C}
	\end{array}\right].\end{align*}
	In order to show that a Gaussian mixture measure is well-defined, we must show that each component is well-defined, that is, each \(\mathcal{V}_{j}\) is a positive definite, self-adjoint, and trace-class operator. 
	
	It is obvious that \(\mathcal{V}_{j}\) is self-adjoint. Then it follows from
	\begin{align*}((\xi,\eta),\mathcal{V}_{j}(\xi,\eta))=\frac{1}{2}(\eta,\beta^2\mathcal{C}_j\eta)
	+\frac{1}{2}(\xi+\gamma \eta,\mathcal{C}(\xi+\gamma \eta))\end{align*}
	that \(\mathcal{V}_{j}\) is positive. In order to show that \(\mathcal{V}_{j}\) is a trace-class operator, we need to compute the eigenpairs of \(\mathcal{V}_{j}\). Note that \(\mathcal{C}_j\) shares the same eigenfunctions, we claim that \(\mathcal{V}_{j}\) has eigenfunctions of the form \((\phi_i,t\phi_i)\), i.e.,
	\begin{align*}
	\left[\begin{array}{cc}
	\mathcal{C} & \gamma \mathcal{C} \\
	\gamma \mathcal{C} & \beta^2\mathcal{C}_j+\gamma^2\mathcal{C}
	\end{array}\right]\left[\begin{array}{cc}
	\phi_i\\t\phi_i
	\end{array}\right]=\left[\begin{array}{cc}
	(1+\gamma t)\lambda_{i}\phi_i\\
	(\gamma \lambda_{i}+\beta^2\lambda_{ji}t+\gamma^2\lambda_{i}t)\phi_i
	\end{array}\right]=\left[\begin{array}{cc}
	(1+\gamma t)\lambda_{i}\phi_i\\
	(\frac{\gamma}{t}\lambda_{i}+\beta^2\lambda_{ji}+\gamma^2\lambda_{i})t\phi_i
	\end{array}\right].
	\end{align*}
	Hence we have
	\begin{align}
	\lambda_{i}+\gamma t\lambda_{i}=\frac{\gamma}{t}\lambda_{i}+\beta^2\lambda_{ji}+\gamma^2\lambda_{i},\label{characteristic equation}
	\end{align}
	which is a  quadratic equation in \(t\). Note that \(\beta^2+\gamma^2=1\), the solution of equation \eqref{characteristic equation} is
	\begin{align}
	t^{\pm}=&\frac{\beta^2}{2\gamma\lambda_{i}}(\lambda_{ji}-\lambda_{i})
	\pm\sqrt{\frac{\beta^4}{4\gamma^2\lambda_{i}^2}(\lambda_{ji}-\lambda_{i})^2+1}.\label{sol of cf}
	\end{align}
	Thus \(\mathcal{V}_{j}\) have eigenfunction \([\phi_i,t^{\pm}\phi_i]\), and the corresponding eigenvalue \((1+\gamma t^{\pm})\lambda_{i}\). Note that
	\begin{align*}
	span\{[\phi_i,t^{\pm}\phi_i]_{i=1}^\infty\}=span\{[\phi_i,0]_{i=1}^\infty,[0,\phi_i]_{i=1}^\infty\}=\mathcal{H}\times\mathcal{H},
	\end{align*}
	which means we have found all eigenpairs. 
	
	Now we turn to prove that the operator \(\mathcal{V}_{j}\) is of trace class. Denote \( \ell^2 \) as the space of square-summable sequences. It follows from the equivalence among \(\{\mathcal{N}(m_j,\mathcal{C}_j)\}_{j=1}^M\) that \(\mathcal{C}^{-1/2}\mathcal{C}_j\mathcal{C}^{-1/2}-I\) is a Hilbert-Schmidt operator, which means \(\frac{\lambda_{ji}}{\lambda_i}-1\in \ell^2\). For simplicity, we denote \(l_{ji}=\frac{\lambda_{ji}}{\lambda_{i}}-1\in \ell^2\), and the solutions \eqref{sol of cf} become
	\begin{align*}
	t^\pm=&\frac{\beta^2}{2\gamma}l_{ji}
	\pm\sqrt{\frac{\beta^4}{4\gamma^2}l_{ji}^2+1}.
	\end{align*}
	Moreover, the sum of eigenvalues is given by
	\begin{align*}
	\text{tr}(\mathcal{V}_{j})=&\sum_i (1+\gamma t^+)\lambda_{i}+\sum_i (1+\gamma t^-)\lambda_{i}\\
	=&2\sum_i\lambda_{i}+\gamma\sum_i (t^++t^-)\lambda_{i}\\
	=&2\text{tr}(\mathcal{C})+\gamma\sum_i \frac{\beta^2}{\gamma}\lambda_{ji} \lambda_{i}\\
	=&2\text{tr}(\mathcal{C})+\beta^2\sum_i \lambda_{ji} \lambda_{i}<\infty,
	\end{align*}
	which completes the proof. 
\end{proof}

Using the same method, it can be deduced that \(\mu(dv)Q(v,du)\) is a Gaussian mixture measure. As a summary, we have that they are both \(M\)-component Gaussian mixture measures, i.e.,
\begin{align}
\mu_0(du)Q(u,dv)=&\sum_{j=1}^Mw_j\rho_{j}(du,dv),\quad \rho_{j}=\mathcal{N}\left(
m_{j},\mathcal{V}_{j}
\right),\\
\mu_0(dv)Q(v,du)=&\sum_{j=1}^Mw_j\rho_{j}'(du,dv),\quad \rho_{j}'=\mathcal{N}\left(
m_{j}',\mathcal{V}_{j}'
\right),
\end{align}
where the mean functions and covariance operators are as follows:
\begin{align}
m_{j}=&[0,(1-\gamma)) m_j],
&\mathcal{V}_{j}=\left[\begin{array}{cc}
\mathcal{C} & \gamma \mathcal{C} \\
\gamma \mathcal{C} & \beta^2\mathcal{C}_j+\gamma^2\mathcal{C}
\end{array}\right],\\
m_{j}'=&[(1-\gamma) m_j,0],
&\mathcal{V}_{j}'=\left[\begin{array}{cc}
\beta^2\mathcal{C}_j+\gamma^2\mathcal{C}& \gamma \mathcal{C} \\
\gamma \mathcal{C} & \mathcal{C}
\end{array}\right].
\end{align}
We claim that the eigenfunction of \(\mathcal{V}'_{j}\) is \([t\phi_i,\phi_i]^T\), then we have
\begin{align*}
\left[\begin{array}{cc}
\beta^2\mathcal{C}_j+\gamma^2\mathcal{C} & \gamma \mathcal{C} \\
\gamma \mathcal{C} & \mathcal{C}
\end{array}\right]\left[\begin{array}{cc}
t\phi_i\\\phi_i
\end{array}\right]=\left[\begin{array}{cc}
(\gamma \lambda_{i}+\beta^2\lambda_{ji}t+\gamma^2\lambda_{i}t)\phi_i\\
(1+\gamma t)\lambda_{i}\phi_i
\end{array}\right]=\left[\begin{array}{cc}
(\frac{\gamma}{t}\lambda_{i}+\beta^2\lambda_{ji}+\gamma^2\lambda_{i})t\phi_i\\
(1+\gamma t)\lambda_{i}\phi_i
\end{array}\right].
\end{align*}
Hence the eigenvalue \(\lambda_i\) should satisfy
\begin{align*}
\frac{\gamma}{t}\lambda_{i}+\beta^2\lambda_{ji}+\gamma^2\lambda_{i}=\lambda_{i}+\gamma t\lambda_{i}.
\end{align*}
which is the same as equation \eqref{characteristic equation}. Therefore, we find that the eigenvalues of \(\mathcal{V}'_{j}\) are identical to those of \(\mathcal{V}_{j}\). As a summary, it follows that
\begin{itemize}
	\item The eigenpairs of \(\mathcal{V}_{j}\) are \((1 + \gamma t^\pm)\lambda_i\) and \(\left[\phi_i, t^\pm \phi_i\right]\).
	\item The eigenpairs of \(\mathcal{V}_{j}'\) are \((1+\gamma t^\pm)\lambda_i\) and \([t^{\pm}\phi_i,\phi_i].\)
\end{itemize}
\subsection{Proof of Theorem \ref{HS} (the well-definedness of pCN-GM method)} 
\label{SM-well-defined pCN-GM}
\begin{proof}
	Our goal is to prove that \(\mathcal{V}^{-1/2}\mathcal{V}'\mathcal{V}^{-1/2}-I\) is a Hilbert-Schmidt operator, i.e., its eigenvalues \(\{\lambda_i\}_{i=1}^\infty \in \ell^2\). Thus, we compute the eigenvalues first. We only need to compute the eigenpairs of \(\mathcal{V}^{-1/2}\mathcal{V}'\mathcal{V}^{-1/2}\), i.e., to solve the following equation:
	\begin{align}
	\mathcal{V}'x=\eta \mathcal{V}x.\label{VV'V-I equation}
	\end{align}
	Note that \(\mathcal{V}\) and \(\mathcal{V}'\) have the same two-dimensional invariant subspace
	\begin{align*}
	\text{span}\{[\phi_i,t^{\pm}\phi_i]\}=\text{span}\{[\phi_i,0],[0,\phi_i]\}=\text{span}\{[t^{\pm}\phi_i,\phi_i]\}.
	\end{align*}
	Hence, we claim that the solution to \eqref{VV'V-I equation} is of the form \([a\phi_i, b\phi_i]\), and then substitute it back. It follows that
	\begin{align*}\left[\begin{array}{cc}
	\beta^2\mathcal{C}_{j_2}+\gamma^2\mathcal{C} & \gamma \mathcal{C} \\
	\gamma \mathcal{C} & \mathcal{C}
	\end{array}\right]\left[\begin{array}{cc}
	a\phi_i \\
	b\phi_i
	\end{array}\right]=\eta
	\left[\begin{array}{cc}
	\mathcal{C} & \gamma \mathcal{C} \\
	\gamma \mathcal{C} & \beta^2\mathcal{C}_{j_1}+\gamma^2\mathcal{C}
	\end{array}\right]\left[\begin{array}{cc}
	a\phi_i \\
	b\phi_i
	\end{array}\right].
	\end{align*}
	Currently, we have a system of equations
	\begin{align}
	(\beta^2\lambda_{j_2i}+\gamma^2\lambda_{i})a+\gamma \lambda_{i}b
	=\eta(\lambda_{i}a+\gamma \lambda_{i}b) \label{HS1},\\
	\gamma \lambda_{i}a+\lambda_{i}b=\eta
	[\gamma \lambda_{i}a
	+(\beta^2\lambda_{j_1i} \gamma^2\lambda_{i} )b].\label{HS2}
	\end{align} 
	It is apparent that \(b\neq 0\), thus we let \(b=\gamma \) for simplicity. Moreover, let
	\begin{align*}
	l_{j_1i}=\frac{\lambda_{j_1i}-\lambda_{i}}{\lambda_{i}},\quad
	l_{j_2i}'=\frac{\lambda_{j_2i}-\lambda_{i}}{\lambda_{i}},
	\end{align*}
	then we have \(l_{j_1i},l_{j_2i}'\in \ell^2\) due to the equivalence among the components of the Gaussian mixture measure. Then it follows from \eqref{HS1} that
	\begin{align}
	(\beta^2\lambda_{j_2i}+\gamma^2\lambda_{i})a
	+\gamma^2 \lambda_{i}
	=&\eta(\lambda_{i}a+\gamma^2 \lambda_{i}),\notag\\
	(\beta^2(\lambda_{j_2i}-\lambda_{i})+\lambda_{i})a
	+\gamma^2 \lambda_{i}
	=&\eta\lambda_{i}(a+\gamma^2 ),\notag\\
	(\beta^2l_{j_2i}'+1)a
	+\gamma^2
	=&\eta(a+\gamma^2 ),\notag\\
	(\beta^2l_{j_2i}'+1-\eta)a
	=&\gamma^2 (\eta-1).\label{HS3}
	\end{align}
	It follows from \eqref{HS2} that
	\begin{align}
	\lambda_{i}a+\lambda_{i}
	=&\eta[ \lambda_{i}a
	+(\beta^2\lambda_{j_1i}+\gamma^2\lambda_{i})],\notag\\
	\lambda_{i}( a+1)
	=&\eta[\lambda_{i}a
	+(\beta^2(\lambda_{j_1i}-\lambda_{i})+\lambda_{i})],\notag\\
	a+1
	=&\eta[a
	+(\beta^2l_{j_1i}+1)],\notag\\
	(1-\eta) a
	=&\eta(\beta^2l_{j_1i}+1)-1.\label{HS4}
	\end{align}
	Combining equations \eqref{HS3},\eqref{HS4} and eliminating \(a\), it follows that
	\begin{align*}
	(\beta^2l_{j_2i}'+1-\eta)[\eta(\beta^2l_{j_1i}+1)-1]
	=&-\gamma^2(\eta-1)^2,\\
	[\eta-(\beta^2l_{j_2i}'+1)][\eta(\beta^2l_{j_1i}+1)-1]
	=&\gamma^2(\eta-1)^2,
	\end{align*}
	By organizing this equation, we obtain
	\begin{align}
	(\beta^2l_{j_1i}+1-\gamma^2)\eta^2+[-(\beta^2l_{j_2i}'+1)(\beta^2l_{j_1i}+1)-1+2\gamma^2]\eta
	+\beta^2l_{j_2i}'+1-\gamma^2=&0,\notag\\
	\beta^2(l_{j_1i}+1)\eta^2-[\beta^4l_{j_1i}l_{j_2i}'+\beta^2(l_{j_1i}+l_{j_2i}')+2\beta^2]\eta+\beta^2(l_{j_2i}'+1)=&0,\notag\\
	(l_{j_1i}+1)\eta^2-[\beta^2l_{j_1i}l_{j_2i}'+(l_{j_1i}+l_{j_2i}')+2]\eta+l_{j_2i}'+1=&0.\label{HS5}
	\end{align}
	Denote the solutions of equation \eqref{HS5} as \(\eta_{i+},\eta_{i-}\). For simplicity we compute the sum and product  of the solutions. Recall that \(l_{j_1i},l_{j_2i}'\in \ell^2\) for \(1\leq j_1,j_2\leq M\), thus \(l_{j_1i},l_{j_2i}'\rightarrow 0,i\rightarrow \infty\).
	\begin{itemize}
		\item The sum is 
		\begin{align*}\eta_{i+}+\eta_{i-}
		=\frac{\beta^2l_{j_1i}l_{j_2i}'+(l_{j_1i}+l_{j_2i}')+2}
		{l_{j_1i}+1}
		\rightarrow 2,\ i\rightarrow \infty.\end{align*}
		Denote 
		\begin{align*}l_{3i}:=\eta_{i+}+\eta_{i-}-2
		=\frac{\beta^2 l_{j_1i} l_{j_2i}' + (l_{j_2i}' - l_{j_1i})}{l_{j_1i} + 1},\end{align*}
		then we have \(l_{3i}\in \ell^2 \). 
		\item The product is 
		\begin{align*}\eta_{i+}\eta_{i-}=\frac{l_{j_2i}'+1}{l_{j_1i}+1}\rightarrow 1.\end{align*}
		Denote 
		\begin{align*}l_{4i}:=\eta_{i+}\eta_{i-}-1=\frac{l_{j_2i}'-l_{j_1i}}{l_{j_1i}+1},\end{align*}
		then we have \(l_{2i}\in \ell^2\).
	\end{itemize}
	We need to show that \(\mathcal{V}^{-1/2}\mathcal{V}'\mathcal{V}^{-1/2}-I\) is a Hilbert-Schmidt operator, i.e. \(\{\eta_{i+}-1,\eta_{i-}-1\}_{i=1}^\infty\in \ell^2.\) Let's compute the sum in the subspace firstly:
	\begin{align*}
	(\eta_{i+}-1)^2+(\eta_{i-}-1)^2
	=&(\eta_{i+}^2+\eta_{i-}^2)+2-2(\eta_{i+}+\eta_{i-})\\
	=&(\eta_{i+}+\eta_{i-})^2-2(\eta_{i+}\eta_{i-}-1)-2(\eta_{i+}+\eta_{i-})\\
	=&(\eta_{i+}+\eta_{i-})(\eta_{i+}+\eta_{i-}-2)-2(\eta_{i+}\eta_{i-}-1),
	\end{align*}
	Using \(\eta_{i+}+\eta_{i-}=l_{3i}+2\) and \(\eta_{i+}\eta_{i-}=l_{4i}+1\), it follows that
	\begin{align*}
	(\eta_{i+}-1)^2+(\eta_{i-}-1)^2 
	=&l_{3i}^2+2l_{3i}-2l_{4i}\\
	=&l_{3i}^2+2 \left[\frac{\beta^2 l_{j_1i} l_{j_2i}'+(l_{j_2i}'-l_{j_1i})}
	{l_{j_1i}+1}-\frac{l_{j_2i}'-l_{j_1i}}{l_{j_1i}+1}\right]\\
	=&l_{3i}^2+2\beta^2\frac{l_{j_1i} l_{j_2i}'}{l_{j_1i}+1},
	\end{align*}
	which is summable because \(l_{j_1i}, l_{j_2i}', l_{3i} \in \ell^2\), and the proof is completed.
\end{proof}

\subsection{Discussion of Remark \ref{remark of pCN-GM} (uniqueness of pCN-GM proposal)} \label{SM-UniquenessOfpCNGM}
\begin{proof}
	This proof process is similar to that described in Section \ref{SM-well-defined pCN-GM}. Let \(\beta^2+\gamma^2=I\), and we will show that \(I\) must be \(1\) to ensure the equivalence. The condition is that \(\mathcal{V}^{-1/2}\mathcal{V}'\mathcal{V}^{-1/2}-I\) is a Hilbert-Schmidt operator, i.e., its eigenvalue \(\lambda_i \in \ell^2\). To compute the eigenvalue, we solve
	\begin{align*}
	\mathcal{V}'x = \eta \mathcal{V}x.
	\end{align*}
	We claim that the solution is of the form \([a\phi_i, b\phi_i]\). It follows that
	\begin{align*}\left[\begin{array}{cc}
	\beta^2\mathcal{C}_{j_2}+\gamma^2\mathcal{C} & \gamma \mathcal{C} \\
	\gamma \mathcal{C} & \mathcal{C}
	\end{array}\right]\left[\begin{array}{cc}
	a\phi_i \\
	b\phi_i
	\end{array}\right] = \eta
	\left[\begin{array}{cc}
	\mathcal{C} & \gamma \mathcal{C} \\
	\gamma \mathcal{C} & \beta^2\mathcal{C}_{j_1}+\gamma^2\mathcal{C}
	\end{array}\right]\left[\begin{array}{cc}
	a\phi_i \\
	b\phi_i
	\end{array}\right].\end{align*}
	Currently, we have a system of equations:
	\begin{align}
	(\beta^2\lambda_{j_2i}+\gamma^2\lambda_{i})a+\gamma \lambda_{i}b
	=\eta(\lambda_{i}a+\gamma \lambda_{i}b) \label{HS1'},\\
	\gamma \lambda_{i}a+\lambda_{i}b=\eta
	[\gamma \lambda_{i}a
	+(\beta^2\lambda_{j_1i} \gamma^2\lambda_{i} )b].\label{HS2'}
	\end{align} 
	It is apparent that \(b\neq 0\), thus we let \(b=\gamma \) for simplicity. Moreover, let
	\begin{align*}
	l_{j_1i}=\frac{\lambda_{j_1i}-\lambda_{i}}{\lambda_{i}}, \quad l_{j_2i}'=\frac{\lambda_{j_2i}-\lambda_{i}}{\lambda_{i}},
	\end{align*}
	then \(l_{j_1i}, l_{j_2i}' \in \ell^2\) due to the equivalence among components of the Gaussian mixture measure. From this step, the calculation process differs from Section  \ref{SM-well-defined pCN-GM}. It follows from equation \eqref{HS1'} that
	\begin{align}
	(\beta^2\lambda_{j_2i}+\gamma^2\lambda_{i})a
	+\gamma^2 \lambda_{i}
	=&\eta(\lambda_{i}a+\gamma^2 \lambda_{i}),\notag\\
	(\beta^2(\lambda_{j_2i}-\lambda_{i})+I\lambda_{i})a
	+\gamma^2 \lambda_{i}
	=&\eta\lambda_{i}(a+\gamma^2 ),\notag\\
	(\beta^2 l_{j_2i}' + I) a + \gamma^2 =& \eta (a + \gamma^2 ), \notag \\
	(\beta^2 l_{j_2i}' + I - \eta) a =& \gamma^2 (\eta - 1). \label{HS3'}
	\end{align}
	It follows from equation \eqref{HS2'} that
	\begin{align}
	\lambda_{i}a+\lambda_{i}
	=&\eta[ \lambda_{i}a
	+(\beta^2\lambda_{j_1i}+\gamma^2\lambda_{i})],\notag\\
	\lambda_{i}( a+1)
	=&\eta[\lambda_{i}a
	+(\beta^2(\lambda_{j_1i}-\lambda_{i})+I\lambda_{i})],\notag\\
	a + 1 =& \eta[a + (\beta^2 l_{j_1i} + I)], \notag \\
	(1 - \eta) a =& \eta(\beta^2 l_{j_1i} + I) - 1. \label{HS4'}
	\end{align}
	Combining equations \eqref{HS3'} and \eqref{HS4'}, we then eliminate \(a\) and find that
	\begin{align*}
	(\beta^2 l_{j_2i}' + I - \eta)[\eta(\beta^2 l_{j_1i} + I) - 1] =& -\gamma^2 (\eta - 1)^2, \\
	(\eta - (\beta^2 l_{j_2i}' + I))[\eta(\beta^2 l_{j_1i} + I) - 1] =& \gamma^2 (\eta - 1)^2.
	\end{align*}
	By organizing this equation, we obtain
	\begin{align}
	(\beta^2 l_{j_1i} + I - \gamma^2) \eta^2 + [-(\beta^2 l_{j_2i}' + I)(\beta^2 l_{j_1i} + I) - 1 + 2\gamma^2] \eta + \beta^2 l_{j_2i}' + I - \gamma^2 =& 0, \notag \\
	\beta^2 (l_{j_1i} + 1) \eta^2 - [\beta^4 l_{j_1i} l_{j_2i}' + \beta^2 (l_{j_1i} + l_{j_2i}') + I^2 + 1 - 2\gamma^2] \eta + \beta^2 (l_{j_2i}' + 1) =& 0, \notag \\
	(l_{j_1i} + 1) \eta^2 - \left[\beta^2 l_{j_1i} l_{j_2i}' + (l_{j_1i} + l_{j_2i}') + \frac{I^2 + 1 - 2\gamma^2}{\beta^2}\right] \eta + l_{j_2i}' + 1 =& 0, \notag \\
	(l_{j_1i} + 1) \eta^2 - \left[\beta^2 l_{j_1i} l_{j_2i}' + (l_{j_1i} + l_{j_2i}') + \frac{(I - 1)^2}{\beta^2} + 2\right] \eta + l_{j_2i}' + 1 =& 0. \label{HS5'}
	\end{align}
	Denote the solutions to equation \eqref{HS5'} as \(\eta_1,\eta_2\). We still compute the sum and product of the solutions first. 
	\begin{itemize}
		\item The sum is \(\eta_1+\eta_2
		=\frac{\beta^2ll'+(l+l')+2+\frac{(I-1)^2}{\beta^2}}
		{l+1}\rightarrow2+\frac{(I-1)^2}{\beta^2} \), and let \begin{align*}l_1:=\eta_1+\eta_2-2=\frac{\beta^2ll'+(l'-l)+\frac{(I-1)^2}{\beta^2}}{l+1}\rightarrow \frac{(I-1)^2}{\beta^2}.\end{align*}
		\item The product is \(\eta_1\eta_2=\frac{l'+1}{l+1}\rightarrow 1\), and let \begin{align*}l_2:=\eta_1\eta_2-1=\frac{l'-l}{l+1}\in \ell^2.\end{align*}
	\end{itemize}
	
	Denote the solutions to equation \eqref{HS5'} as \(\eta_{i+},\eta_{i-}\). We still compute the sum and product of the solutions first. 
	\begin{itemize}
		\item The sum is \begin{align*}\eta_{i+}+\eta_{i-}
		=\frac{\beta^2 l_{j_1i} l_{j_2i}'+(l_{j_1i}+l_{j_2i}')+2+\frac{(I-1)^2}{\beta^2}}
		{l_{j_1i}+1}\rightarrow2+\frac{(I-1)^2}{\beta^2} ,\end{align*}
		then let \begin{align*}l_{3i}:=\eta_{i+}+\eta_{i-}-2=\frac{\beta^2 l_{j_1i} l_{j_2i}'+(l_{j_2i}'-l_{j_1i})+\frac{(I-1)^2}{\beta^2}}{l_{j_1i}+1}\rightarrow \frac{(I-1)^2}{\beta^2}.\end{align*}
		\item The product is \begin{align*}\eta_{i+}\eta_{i-}=\frac{l_{j_2i}'+1}{l_{j_1i}+1}\rightarrow 1,\end{align*}then let \begin{align*}l_{4i}:=\eta_{i+}\eta_{i-}-1=\frac{l_{j_2i}'-l_{j_1i}}{l_{j_1i}+1}\in \ell^2.\end{align*}
	\end{itemize}
	
	To see if \(\mathcal{V}^{-1/2}\mathcal{V}'\mathcal{V}^{-1/2}-I\) is a Hilbert-Schmidt operator, let's compute the sum of \((\eta-1)^2\) in the subspace firstly:
	\begin{align*}
	(\eta_{i+}-1)^2+(\eta_{i-}-1)^2
	=&(\eta_{i+}+\eta_{i-})(\eta_{i+}+\eta_{i-}-2)-2(\eta_{i+}\eta_{i-}-1)\\
	=&l_{3i}^2+2l_{3i}-2l_{4i}\\
	=&l_{3i}^2+2 \left[\frac{\beta^2 l_{j_1i} l_{j_2i}'+(l_{j_2i}'-l_{j_1i})}
	{l_{j_1i}+1}-\frac{l_{j_2i}'-l_{j_1i}}{l_{j_1i}+1}\right]\\
	=&l_{3i}^2+2\beta^2\frac{l_{j_1i} l_{j_2i}'}{l_{j_1i}+1}.
	\end{align*}
	If \(I=1\), then \(l_{3i} \in \ell^2\) and the above sum is finite because \(l_{j_1i}, l_{j_2i}' \in \ell^2\). If \(I \neq 1\), then \(l_{3i} \rightarrow \frac{(I-1)^2}{\beta^2} \neq 0\), the above sum is infinite.
	
	In conclusion, \(\mathcal{V}^{-1/2}\mathcal{V}'\mathcal{V}^{-1/2}-I\) is a Hilbert-Schmidt operator if and only if \(\beta^2 + \gamma^2 = 1\). Consequently, if we want to use a positive-coefficient linear combination of the current state \(u\) and a Gaussian mixture, then the only choice is our pCN-GM proposal.
\end{proof}

\subsection{Proof of Theorem \ref{GMM approximate post} (denseness of Gaussian mixture measures)} \label{SM-GMM approximate post}
\begin{proof}
	We use the total variation distance
	\begin{align*}
	d(\mu,\nu)=\sup _{|f|_\infty\leq 1}\left|\int fd\mu-\int fd\nu\right|.
	\end{align*}
	First we can approximate \(\int fd\mu_{post}
	=\int f(u)\frac{1}{Z}e^{-\Phi(u)}\mu_0(du)\) by replacing the likelihood with a finite-dimensional version. Let \(Z_N=\int e^{-\Phi(u_N)}d\mu_0^N(u_N)\), it follows that
	\begin{align*}
	|Z_N-Z|=&\left|\int e^{-\Phi(u_N)}d\mu_0^N(u_N)-\int e^{-\Phi(u)}d\mu_0(u)\right|\\
	\leq&\left|\int e^{-\Phi(u_N)}d\mu_0^N(u_N)-\int e^{-\Phi(u_N)}d\mu_0(u)\right|
	+\int |e^{-\Phi(u_N)}- e^{-\Phi(u)}|d\mu_0(u)\\
	=&\int |e^{-\Phi(u_N)}- e^{-\Phi(u)}|d\mu_0(u)\rightarrow 0, N\rightarrow \infty.
	\end{align*}
	Here we need \(\Phi\) to be continuous. Then we can measure the distance between \(\mu_{post}\) and the posterior derived from the finite-dimensional likelihood. Using Dominated Convergence Theorem, we have
	\begin{align}
	&\left|\int fd\mu_{post}-\int f(u)\frac{1}{Z_N}e^{-\Phi(u_N)}\mu_0(du)\right|\nonumber\\
	=&\left| \int f(u)\frac{1}{Z}e^{-\Phi(u)}\mu_0(du)-\int f(u)\frac{1}{Z_N}e^{-\Phi(u_N)}\mu_0(du)\right|\nonumber\\
	\leq &\int |f(u)|\left|\frac{1}{Z}e^{-\Phi(u)}-\frac{1}{Z_N}
	e^{-\Phi(u_N)}\right|\mu_0(du)\nonumber\\
	\leq &\int |f(u)|\left|
	\frac{1}{Z}(e^{-\Phi(u)}
	-e^{-\Phi(u_N)})
	+\left(\frac{1}{Z}
	-\frac{1}{Z_N}\right)e^{-\Phi(u_N)}
	\right|\mu_0(du)\nonumber\\
	\leq &|f|_\infty\epsilon_1\label{estimate1}.
	\end{align}
	By decomposing the prior into its first \(N\) dimensions and the remaining dimensions, we obtain
	\begin{align*}
	\int f(u)\frac{1}{Z_N}e^{-\Phi(u_N)}\mu_0(du)
	=&\iint f(u)\frac{1}{Z_N}e^{-\Phi(u_N)}\mu_0^N(du_N)\mu_0^\perp(du_\perp)\\
	=&\iint f(u)\frac{1}{Z_N}e^{-\Phi(u_N)}p(u_N)\lambda(du_N)\mu_0^\perp(du_\perp),
	\end{align*}
	where \(\lambda\) is the Lebesgue measure on \(\mathbb{R}^N\). Note that \(Z_N=\int e^{-\Phi(u_N)}p(u_N)\lambda(du_N)\), hence \(\frac{1}{Z_N}e^{-\Phi(u_N)}p(u_N)\) is a density function defined on \(\mathbb{R}^N\), which can be approximated by a Gaussian mixture density function. More precisely, denote the sequence of Gaussian mixture density functions as follows:
	\begin{align*}
	g_k(u_N)=&\sum_{i=1}^{I_k} w_{ki}p(u_N;m_{ki},\Sigma_{ki}),\\
	p(u;m,\Sigma)=&\frac{1}{(2\pi)^{N/2}|\det (\Sigma)|^{1/2}}\exp\left\{-\frac{1}{2}(u-m)^T\Sigma^{-1}(u-m)\right\}.
	\end{align*}
	Then it follows from Theorem 5.d in \cite{Nguyen2020}  that the density function can be approximated by \(g_k\), i.e.,\begin{align*}\lim_{k\rightarrow \infty}g_k(u_N) = \frac{1}{Z_N}e^{-\Phi(u_N)}p(u_N), a.e.\end{align*}
	According to the Dominated Convergence Theorem, we have
	\begin{align*}
	&\lim_{k\rightarrow \infty}\iint f(u)\left[g_k(u_N)-\frac{1}{Z_N}e^{-\Phi(u_N)}p(u_N)\right]\lambda(du_N)\mu_0^\perp(du_\perp)\\
	=&\iint \lim_{k\rightarrow \infty}f(u)\left[g_k(u_N)-\frac{1}{Z_N}e^{-\Phi(u_N)}p(u_N)\right]\lambda(du_N)\mu_0^\perp(du_\perp)\\
	=&0.
	\end{align*}
	Hence there exists a \(g_K\) such that
	\begin{align}
	\iint f(u)\left[g_K(u_N)-\frac{1}{Z_N}e^{-\Phi(u_N)}p(u_N)\right]\lambda(du_N)\mu_0^\perp(du_\perp)<\epsilon.\label{GMM_K}
	\end{align}
	Construct a measure \(\tilde{\mu}_{post}^N\) in the first \(N\) dimensions based on the density function \(g_K\) and the Lebesgue measure, and construct the product measure of \(\tilde{\mu}_{post}^N\) and \(\mu_0^\perp\), i.e.,
	\begin{align*}
	\tilde{\mu}_{post}^N(du_N):=&g_K(u_N)\lambda (du_N),\\
	\tilde{\mu}_{post}(du):=&\tilde{\mu}_{post}^N(du_N)\times \mu_0^\perp(du_\perp).
	\end{align*}
	Then it follows from \eqref{GMM_K} that
	\begin{align}
	\left|\int f(u)\frac{1}{Z_N}e^{-\Phi(u_N)}\mu_0(du)-\int  f(u)\tilde{\mu}_{post}(du)\right|<\epsilon_2.
	\label{estimate2}
	\end{align}
	Combining the formulas \eqref{estimate1} and \eqref{estimate2}, we can deduce that
	\begin{align*}
	&\left|\int fd\mu_{post}-\int fd\tilde{\mu}_{post}\right|\\
	\leq &\left|
	\int fd\mu_{post}
	-\int f(u)\frac{1}{Z_N}e^{-\Phi(u_N)}\mu_0(du)
	\right|+\left|
	\int f(u)\frac{1}{Z_N}e^{-\Phi(u_N)}\mu_0(du)
	-\int fd\tilde{\mu}_{post}
	\right|\\
	<&|f|_\infty\epsilon_1+\epsilon_2.
	\end{align*}
	Note that \(\epsilon_1\) and \(\epsilon_2\) do not depend on \(f\), therefore we can take the supremum of \(f\). Consequently, it follows that
	\begin{align*}d(\mu_{\text{post}}, \tilde{\mu}_{\text{post}}) < \epsilon_1 + \epsilon_2.\end{align*}Now we find a Gaussian mixture measure, \(\tilde{\mu}_{\text{post}}\), which approximates the posterior measure. Furthermore, since we only modified the first \(N\) dimensions when constructing \(\tilde{\mu}_{\text{post}}\), this measure is naturally equivalent to the prior measure.
\end{proof}

\subsection{Proof of Theorem \ref{theorem4} (error of Gaussian mixture approximation)}\label{SM-ErrOfGMM}
The proof of the first inequality \eqref{GMneq1} is straightforward, i.e., 
\begin{align*}d(Q,P\mu_j^N ) \leq d(Q,\mu_j) + d(P\mu_{j},P\mu_j^N)\leq d(Q,\mu_{j})+d(\mu_{j},\mu_j^N).\end{align*}
To prove the second inequality \eqref{GMneq2}, let $X_1,\dots,X_N$ be samples obtained by SMC and define the empirical measure
\begin{align*}
\mu^N := \frac1N\sum_{i=1}^N \delta_{X_i}.
\end{align*}
Let \(\rho=\mathcal{N}(0,\mathcal{C}_0)\), and define the Gaussian mixture
\begin{align*}
Q=\mu^N * \mathcal\rho
=
\frac1N\sum_{i=1}^N \bigl(\delta_{X_i} * \mathcal N(0,\mathcal{C}_0)\bigr) 
=\frac1N\sum_{i=1}^N \mathcal N(X_i,\mathcal{C}_0),
\end{align*}
where the convolution of two finite Borel measures is defined by
\begin{align*}
(\nu_1*\nu_2)(A)
:=
\int \nu_1(A-x)\,\nu_2(dx).
\end{align*}
Difference is estimated by
\begin{align*}
d'(Q,\mu)\leq& d'(Q,\mu*\rho)+d'(\mu*\rho,\mu).
\end{align*}
The first term is bounded by
\begin{align}
d'(\mu^N*\rho,\mu*\rho)
&=\sup_{|f|_\infty\le1,\|D^2f\|\leq 1}\sqrt{
	\mathbb E\Big[
	\int fd(\mu^N*\rho)-\int fd(\mu*\rho)
	\Big]^2} \notag\\
&=\sup_{|f|_\infty\le1,\|D^2f\|\leq 1}\sqrt{
	\mathbb E\Big[
	\int (f*\rho)d(\mu^N-\mu*\rho)
	\Big]^2} \notag\\
&\le d'(\mu^N,\mu) \label{neq1}
\end{align}
where $\rho(-A)=\rho(A)$ for all $ A\in \mathcal{B}(\mathcal H)$, \(|f*\rho(x)|=\left|\int f(x-y)\rho(dy)\right|\leq 1\), and \( D^2(f*\rho)=D^2f*\rho\) are employed.

Using once more the identity $\rho(-A)=\rho(A)$ for every $A\in\mathcal{B}(\mathcal{H})$, the second term is bounded by
\begin{align*}
d'(\mu*\rho,\mu)
=&\sup_{|f|_\infty\le1,\|D^2f\|\leq 1}\left|
\int (f*\rho-f)\mu(dx)\right|\\
=&\sup_{|f|_\infty\le1,f\in C_\infty}\left|
\int_{x\in \mathcal{H}} \left(\int_{y\in \mathcal{H}} f(x-y)\rho(dy)-f\right)\mu(dx)\right|\\
=&\sup_{|f|_\infty\le1,\|D^2f\|\leq 1}\left|
\int_{x\in \mathcal{H}} \int_{y\in \mathcal{H}} \left(f(x+y)-f(x)\right)\rho(dy)\mu(dx)\right|.
\end{align*}

Using a Taylor expansion together with $\mathbb E_{Y\sim \rho}Y=0$, we obtain
\begin{align*}
d'(\mu*\rho,\mu)
=&\sup_{|f|_\infty\le1,\|D^2f\|\leq 1}\left|
\int_{x\in \mathcal{H}} \int_{y\in \mathcal{H}} \frac{1}{2}(D^2f(x+ty)y,y)\rho(dy)\mu(dx)\right|\\
\leq &\frac{1}{2}
\int \|y\|^2\rho(dy)=\frac{1}{2} \text{Tr}(\mathcal{C}_0).
\end{align*}
We choose $\rho = \mathcal{N}(0, \mathcal{C}_0)$ such that $\mathrm{Tr}(\mathcal{C}_0) \le 2\epsilon \, d(\mu^N, \mu)$, which implies that
\begin{align}
d'(\mu * \rho, \mu) \le \epsilon d'(\mu^N, \mu). \label{neq2}
\end{align}
Combining estimates \eqref{neq1} and \eqref{neq2} then yields the desired result.

\begin{remark}It is worth noting that a smaller $\epsilon$ in the error bound of this theorem should not be interpreted as always preferable. In practice, taking $\epsilon$ too small may cause the effective sample size to decay rapidly, leading to a breakdown of the sampling procedure. From a theoretical perspective, in the limiting case $\epsilon \to 0$, the mutation step becomes completely ineffective, in the sense that it does not improve sample diversity at all. Moreover, since at each layer the re-weighting and resampling steps tend to reduce sample diversity, the overall effect is a progressive loss of diversity across iterations. In the extreme case, the final SMC output may degenerate to a single sample drawn from the prior.
\end{remark}

\begin{remark}The Gaussian mixture measure employed in the proof is numerically implementable, provided that an appropriate covariance operator \(\mathcal{C}\) is selected. However, most existing strategies for choosing such a covariance operator have been developed for finite-dimensional settings \cite{Markley2010,Halbe2013}. Therefore, in the numerical experiments, we adopt the Gaussian mixture estimation method proposed in~\cite{Feng2018}, which is designed specifically for the infinite-dimensional setting. In the following, we demonstrate that, with an appropriate choice of the covariance operator \(\mathcal{C}\), the two Gaussian mixture constructions can nonetheless exhibit similar numerical behavior.

First, we illustrate the similarity between the Gaussian mixtures estimated by the two methods. Specifically, we consider high-dimensional function samples and construct Gaussian mixtures using both methods. We examine the densities of the resulting Fourier coefficients of samples drawn from the two Gaussian mixtures. The same quantities as in Section~3.2 are used for comparison. Here, we use ``KDE'' to denote the Gaussian mixture construction employed in the proof, since it corresponds to a kernel density estimation approach.
\begin{center}
\captionsetup{aboveskip=0.1pt}
\begin{minipage}{0.3\linewidth}
    \centering
    \includegraphics[width=\linewidth]{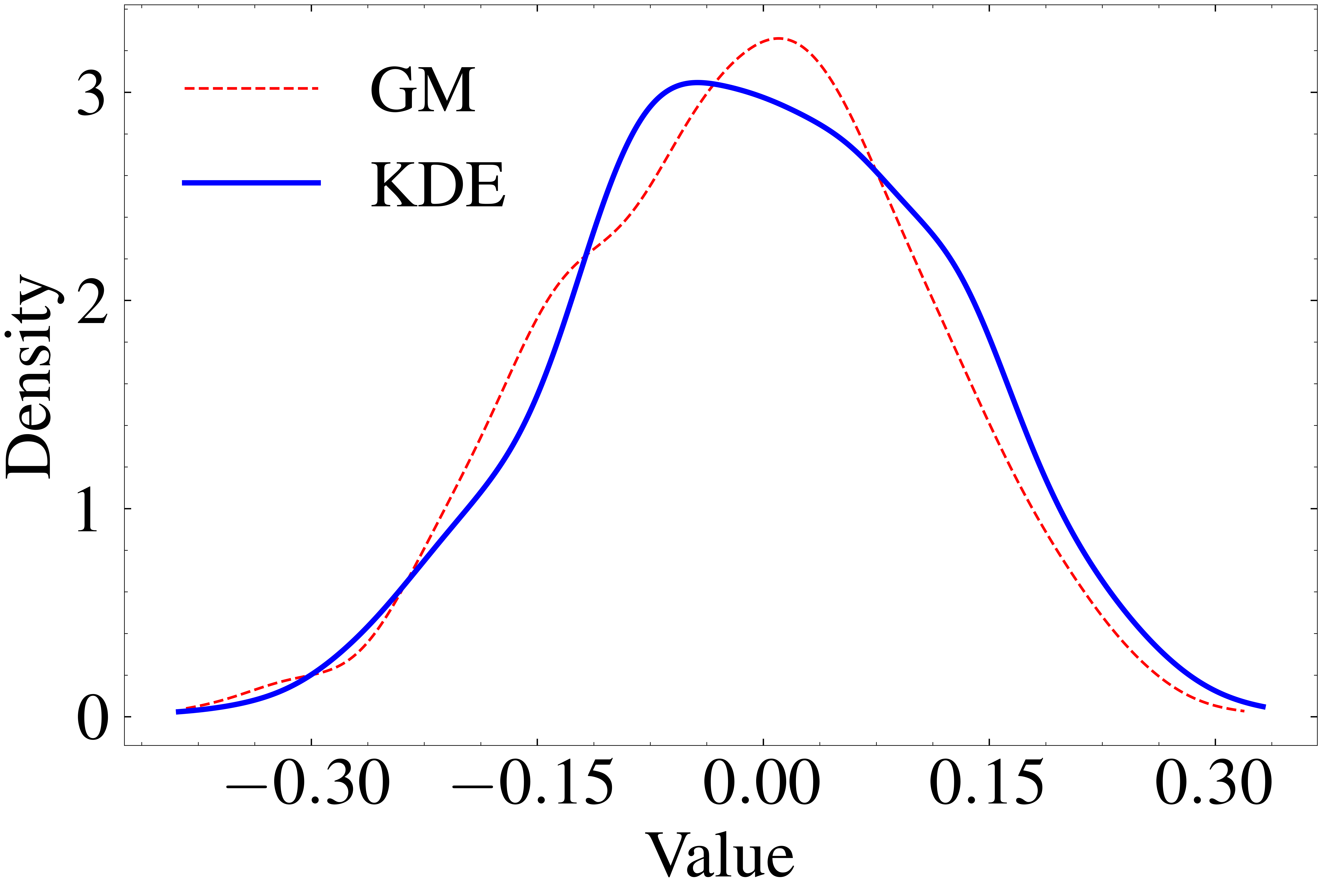}
    \small (a) Density of \(u_1\)
\end{minipage}
\hfill
\begin{minipage}{0.3\linewidth}
    \centering
    \includegraphics[width=\linewidth]{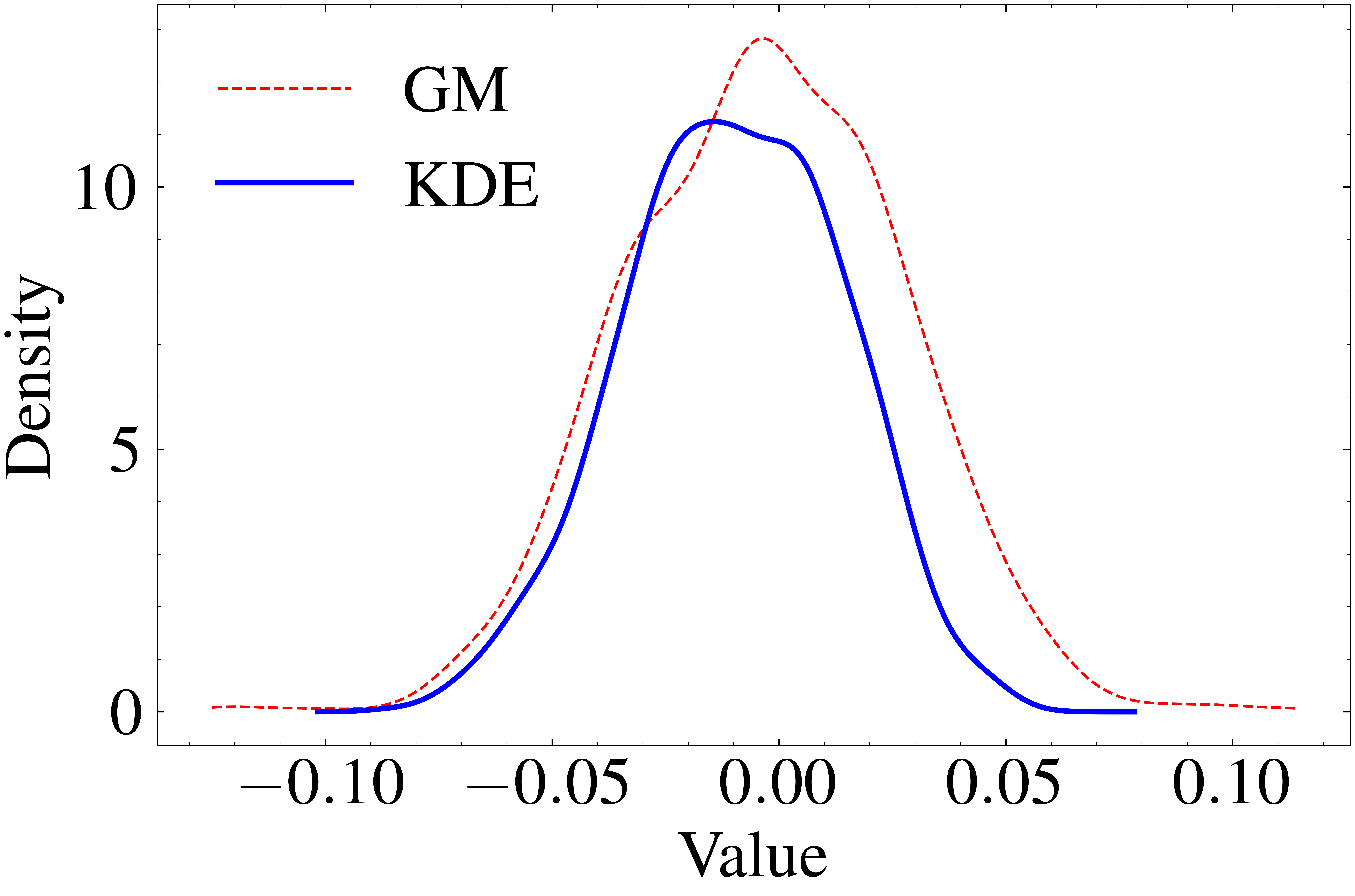}
    \small (a) Density of \(u_4\)
\end{minipage}
\hfill
\begin{minipage}{0.3\linewidth}
    \centering
    \includegraphics[width=\linewidth]{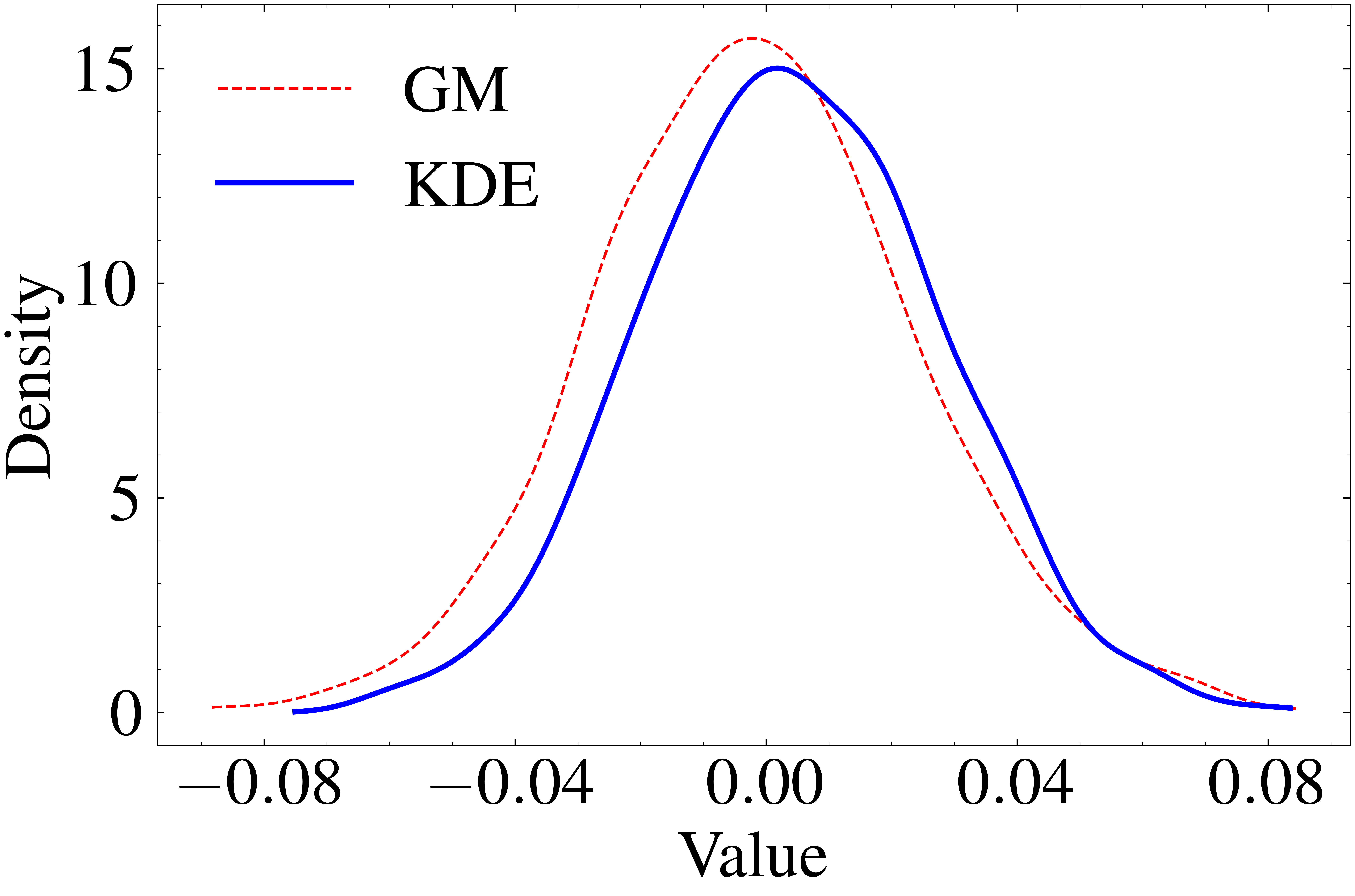}
    \small (a) Density of \(u_7\)
\end{minipage}

\vspace{1ex}  

\begin{minipage}{0.3\linewidth}
    \centering
    \includegraphics[width=\linewidth]{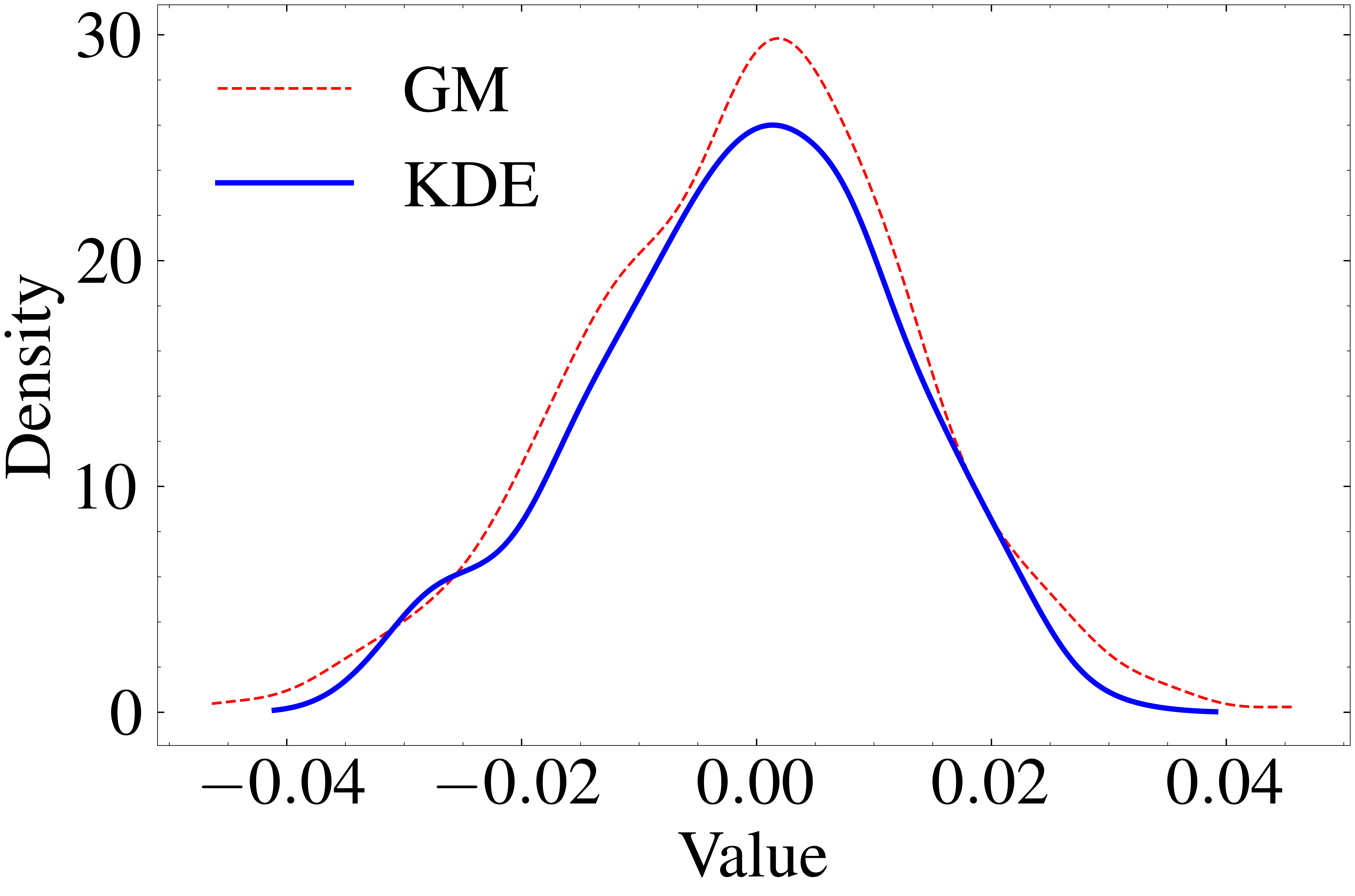}
    \small (d) Density of \(u_{10}\)
\end{minipage}
\hfill
\begin{minipage}{0.3\linewidth}
    \centering
    \includegraphics[width=\linewidth]{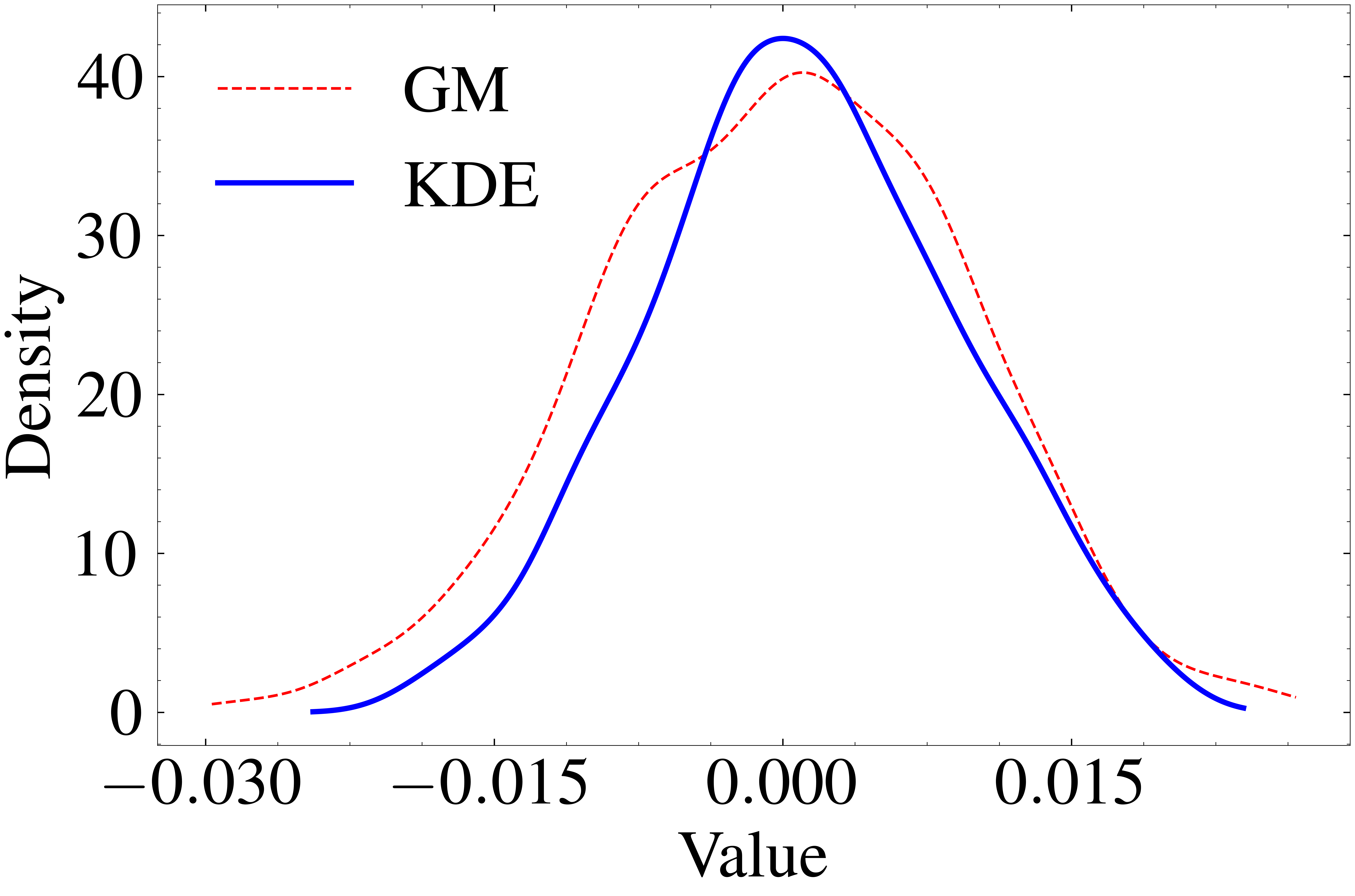}
    \small (e) Density of \(u_{13}\)
\end{minipage}
\hfill
\begin{minipage}{0.3\linewidth}
    \centering
    \includegraphics[width=\linewidth]{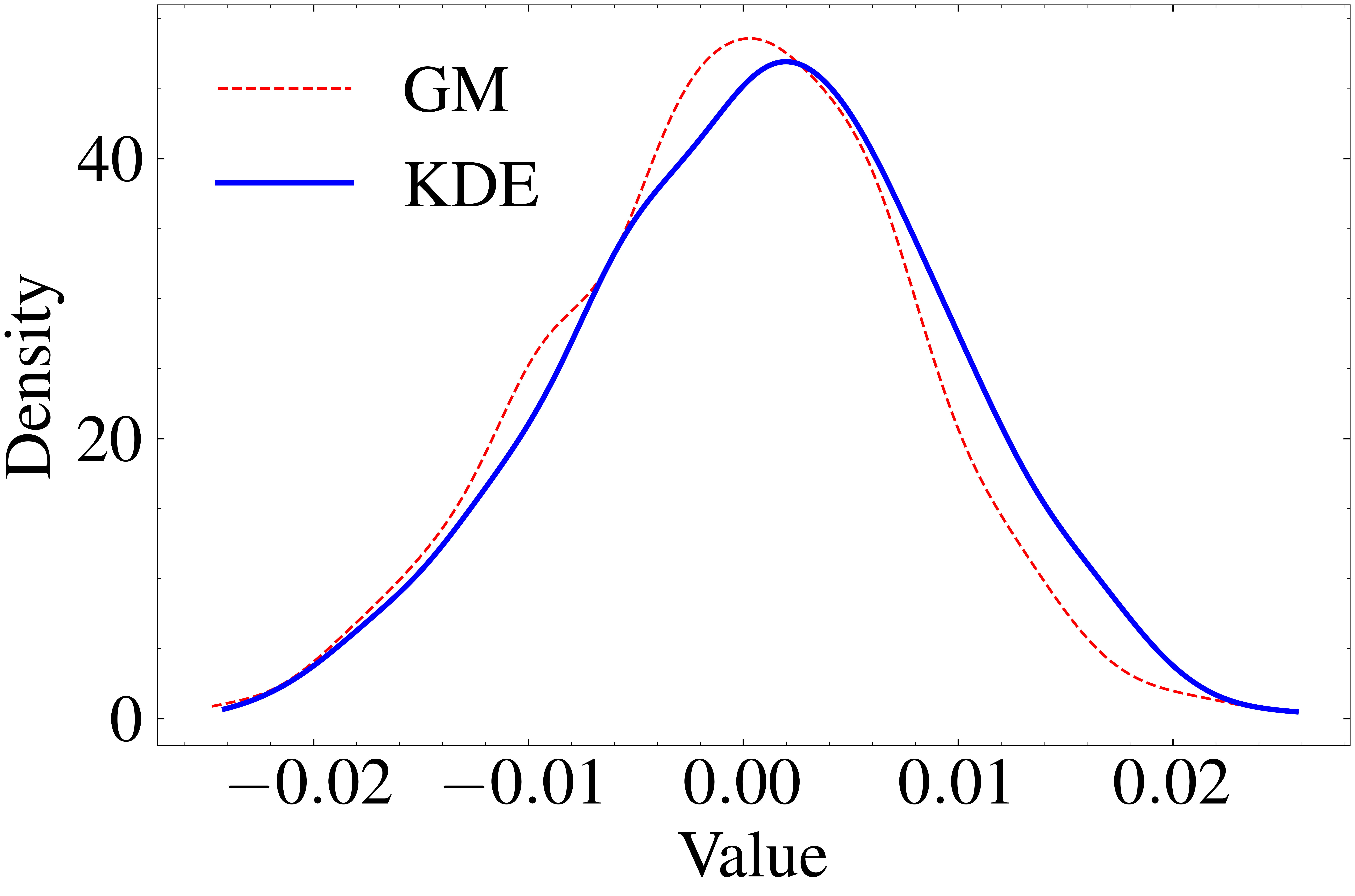}
    \small (f) Density of \(u_{16}\)
\end{minipage}
\par\vspace{0.7em}
\small Figure D.1: Densities of the 1st, 4th, 7th, 10th, 13th, 16th Fourier coefficients.
\end{center}

Second, we employ both Gaussian mixture constructions as the mutation step in SMC and compare the resulting estimates. For the construction used in the proof, we take the sample covariance matrix, multiplied by a scaling parameter \(s\), as the covariance matrix of each Gaussian component. The resulting posterior mean, together with the true value, is shown below.
\begin{center}
\captionsetup{aboveskip=0.1pt}
\begin{minipage}{0.45\linewidth}
    \centering
    \includegraphics[width=\linewidth]{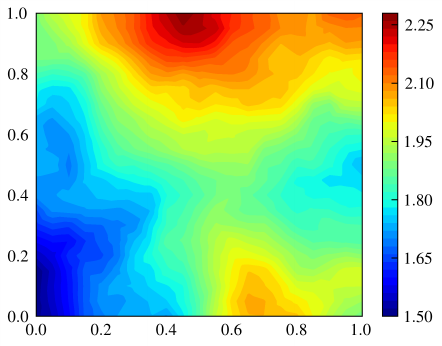}
    \small (a) Background truth
\end{minipage}
\hfill
\begin{minipage}{0.45\linewidth}
    \centering
    \includegraphics[width=\linewidth]{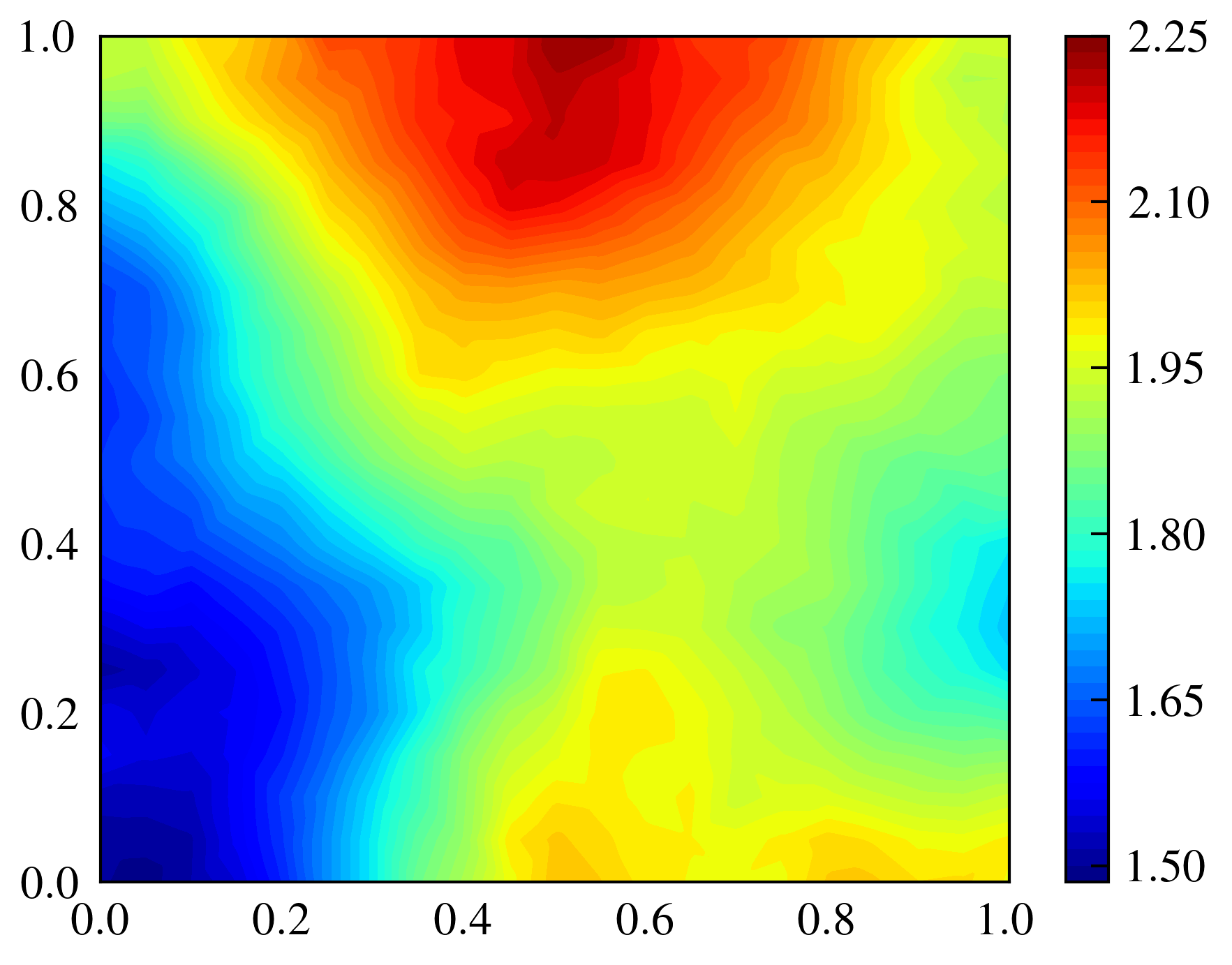}
    \small (b) Mean of SMC-GM
\end{minipage}
\par\vspace{0.7em}
\small Figure D.2: Performance of SMC with KDE estimation.
\end{center}
Numerical results indicate that the KDE-based estimation also achieves sufficient mixing and yields results comparable to those obtained with SMC-GM and SMC-pCN, provided that a suitable scaling parameter $s= 0.04$ is selected. The relative error between the sample mean and the truth is $0.0296$. For reference, the relative errors reported in the main text for SMC-pCN and SMC-GM are $0.0271$ and $0.0254$, respectively. Therefore, both methods can achieve good SMC sampling performance when appropriate parameters are chosen.
\end{remark}

\color{black}

\end{appendix}
\end{document}